\documentclass{amsart}
\usepackage{amssymb}

\usepackage{amsmath,amsthm,amssymb,latexsym,xcolor}
\usepackage[notref,notcite]{showkeys}
\usepackage{epsfig}
\usepackage{enumerate}
\usepackage{mathrsfs}

\newcounter{personaln}

\input tcilatex

\begin{document}
\title[Linear-fractional BGW and RW's]{Sterile versus prolific individuals pertaining to linear-fractional
Bienaym\'{e}-Galton-Watson trees}
\author{Thierry Huillet$^{\left( 1\right) *}$, Servet Mart\'{i}nez$^{\left( 2\right)
}$}
\address{$^{\left( 1\right) }$Laboratoire de Physique Th\'{e}orique et
Mod\'{e}lisation\\
CY Cergy Paris University, CNRS UMR-8089\\
Site de Saint Martin\\
2 avenue Adolphe-Chauvin\\
95302 Cergy-Pontoise, FRANCE\\
E-mail: thierry.huillet@cyu.fr\\
$^{\left( 2\right) }$Centro de Modelamiento Matem\'{a}tico,\\
Facultad de Ciencias F\'{i}sicas y Matem\'{a}ticas,\\
Universidad de Chile\\
Beauchef 851, Edificio Norte Piso 7\\
CP 837 0456\\
Santiago CHILE\\
Email: smartine@dim.uchile.cl}
\maketitle

\begin{abstract}
In a Bienaym\'{e}-Galton-Watson process for which there is a positive
probability for individuals of having no offspring, there is a subtle
balance and dependence between the sterile nodes (the dead nodes or leaves)
and the prolific ones (the productive nodes) both at and up to the current
generation. We explore the many facets of this problem, especially in the
context of an exactly solvable linear-fractional branching mechanism at all
generation. Eased asymptotic issues are investigated. Relation of this
special branching process to skip-free to the left and simple random walks'
excursions is then investigated. Mutual statistical information on their
shapes can be learnt from this association.\newline



$^{*}$ corresponding author
\end{abstract}

\section{Introduction and summary of the results}

Bienaym\'{e}-Galton-Watson (BGW) branching processes have for long been a
milestone in the understanding of multiplicative cascade phenomena [Good
(1949) and Otter (1949)], starting with the extinction of family names in
population dynamics, in the second half of the 19-th century. See [Kendall,
(1966) and Jagers, (2020)] for historical background. BGW processes also
appear as crude models for the spread of rumors, gravitational clustering
[Sheth, (1996)], cosmic-ray cascades and the proliferation of free neutrons
in nuclear fission, [Harris, (1963)].

In a BGW branching tree process for which (our background hypothesis
throughout this paper): 
\begin{eqnarray*}
&&\text{(i) there is a positive probability for individuals of having no
offspring,} \\
&&\text{(ii) the offspring number has all its moments,}
\end{eqnarray*}
there is a subtle balance and dependence between the present (the state of
the population at some current generation) and the past (for example the
cumulated number of nodes before this current generation). In both cases,
the population, as a tree, can be split into two main types of individuals:
the sterile (the dead leaves) and the prolific ones (the living nodes), both
at and up to the currently observed generation number. Specific interest
into these two types of individuals goes back at least to [R\'{e}nyi,
(1959)].

By making extensive use of the apparatus of iterated generating functions,
we explore the balance and dependence between the two. We consider the many
facets of this problem, especially in the context of the linear-fractional
(LF) branching mechanism (obeying (i)(ii)), leading to explicit computations
at all generation as a result of the stability of LF transformations under
composition. Such explicit results, although being specific, illuminate the
known asymptotic results concerning more general branching processes, and,
on the other hand, may bring insight into less investigated aspects of the
theory of branching processes. In places, we will then make allusions to
other important branching mechanisms. Depending on the BGW process being
(sub-)critical, supercritical, or nearly supercritical, asymptotic issues
are investigated when the generation number goes to infinity. This is part
of a general program of understanding on how the full past interact with the
present in a BGW process. The organization of the paper is as follows:

- \textit{Section} $2$: We recall some basics on the current population size
of BGW processes satisfying (i)(ii), including some well-known limit laws,
conditions of criticality and conditionings on criticalities. BGW processes
with two-parameters LF branching mechanisms were used in the extinction of
family names problem, [Steffenson, (1930, 1933)]; they deserve specific
interest being an exactly solvable model. The discrete LF BGW process is
indeed an important particular case being amenable to explicit computations,
including - the one of its extinction probability,- the one of the important
fixed point parameter $\tau $ defined in (\ref{f6}) and used in the
conditioning of non-critical BGW processes to critical ones and in the law
of the total progeny - the one of its $n$-step transition matrix leading to
a tractable potential theory - the exact expression of a Kolmogorov constant
- an explicit limit law of the $Q$-process. It is furthermore embeddable in
a continuous-time binary BGW process. We mention a Harris derivation of its
`stationary distribution' [Harris, 1963]. We derive the joint law of the
sterile versus prolific individuals at a current generation $n$.

- \textit{Section} $3$: We are concerned with the interaction between the
past and the present of BGW processes, starting with the total progeny of
one or more founders (forests of trees).

The analysis of the total progeny marginal is first introduced in Section
3.1 and developed in Section 3.5. Section 3.6 considers the critical case.

We next derive a general recurrence for the joint law of the sterile and
prolific individuals at and up to generation $n$, starting from a single
founder; see (\ref{f16}) of Section 3.2. We particularize this recurrence to
study:

- the joint law of the cumulated number of sterile against the cumulated
number of prolific individuals in the BGW tree. This is introduced in
Section 3.3 and developed in Section 3.7.

- the joint law of the cumulated number of sterile against the current
number of prolific individuals in the BGW tree; see Section 3.4. We show
through a simple example that answering the question: do the number of ever
dead in a population outnumber or not the currently living is a delicate
question that can be answered positively or negatively in the LF case,
depending on the range of its two parameters.

- the joint law of the total number of leaves (sterile individuals) against
the total progeny, obtaining the limiting large deviation function in
Section 3.7.

- \textit{Section }$4$: Two relations of BGW processes to random walks
(RW's) are next investigated. The first concerns general skip-free to the
left RW's related to BGW processes through a discrete version of the
Lamperti time-change theorem known for continuous-state branching processes,
[Lamperti, (1967)]. It connects the full past and the present of BGW
processes. Making use of the scale function [Marchal, (2001)], we derive the
law of the maximum (width) of a BGW till its first extinction, concomitantly
with the law of the supremum of the RW till its first hitting time of $0$.
Both are jointly finite (or not) depending on the criticality (or not) of
the BGW process. We derive an explicit computation in the LF and binary
fission cases.

The second concerns simple random walks (SRW's) with nearest-neighbours
moves. We revisit and highlight the Harris construction of a geometric BGW
tree nested inside a SRW excursion with no holding probability, [Harris,
1952]. A node of the tree at some height $h-1$ branches when it reveals
minima of the SRW at $h$. We propose an extension of this construction to
SRW's with holding probabilities (whose sample paths shows highlands and
valleys), the nested tree having now a LF branching mechanism. This relation
allows to derive useful mutual information on both the tree and the random
walk: in particular, we show that the first hitting time of the SRW is
approximately twice the total progeny of its nested tree whose time to
extinction is approximately the height of the SRW. The width of the nested
BGW process is the largest size of the SRW valleys.

%Let $A(r)$ be the area of a $d$-dimensional sphere of radius $r>0$. We can
%interpret $N_{n}$ as the number of external nodes uniformly distributed on
%the $n$-spheres of area $a_{n}:=A(nr)$ or $A(r/n)$, $n\geq 1$, which yields%
%a geometric vision of the BGW process along nested spheres, either ex- or%
%implosive: an onion model. [Nils Berglund: http://images.math.cnrs.fr/La
%probabilit\'{e} d'extinction d'une esp\`{e}ce menac\'{e}e]

%Consider an infinite step-function as a broken-line interface separating
%solid at the bottom to liquid on the top. Suppose the step-function is
%bounded above at height $0$. Initially, an infinite sea at level $1$ covers
%the whole solid structure. When the level of the sea is lowered by one unit,
%many nested sub-seas (lakes) will appear each of which can subsequently
%reveal additional sub-sub-seas or disappear being not deep enough....the
%problem (for this construction to be amenable to a BGW process) is that
%sub-seas have to be rescaled to the mother sea before envisaging the
%appearance of next sub-sub-seas.

%Let $X$ be a random set of points of the plane (Poisson in[ Sheth, (1996)]).
%Tag one of its particular point and center a ball of radius $r>0$ on it. If
%the ball does not catch other points of $X$, the process stops (extinction).
%If it does, center a ball of radius $r$ on each of the first generation
%points (the offspring) and iterate the process. The problem is that there
%may be common points caught by first (and subsequent) generation points: the
%mutual independence of the second (and subsequent) offspring fail in general.

\section{Generalities on Bienaym\'{e}-Galton-Watson (BGW) branching processes
}

We start with generalities on such BGW processes, avoiding the case
displaying finite-time explosion [Sagitov and Lindo (2015)].

\subsection{Current population size: the probability generating function
(p.g.f.) approach}

Consider a discrete-time BGW branching process [Harris (1963); Athreya and
Ney (1972)] whose reproduction law is given by the probability law $\mathbf{P%
}\left( M=m\right) =:\pi _{m}$, $m\geq 0$ for the number $M$ of offspring
per capita. At each generation, each individual $i$ alive, independently of
another, generates $M_{i}$ offspring, with $M_{i}\overset{d}{=}M$ (in
distribution). Unless specified otherwise, we assume $\pi _{0}>0$ so that
the process can go extinct. We let $\phi \left( z\right) =\mathbf{E}\left(
z^{M}\right) =\sum_{m\geq 0}\pi _{m}z^{m}$ be the p.g.f. of $M$, so with $%
\phi \left( 1\right) =1$ and $\phi \left( z\right) $ has convergence radius $%
z_{*}\in \left( 1,\infty \right] $ ($M$ has all its moments finite and
geometric tails). In such cases, $\phi \left( z\right) $ is analytic in the
open disk $\left| z\right| <z_{*}$ of the complex plane. The latter
assumptions also guarantee the existence of two real fixed points to the
equation $\phi \left( z\right) =z$, one of which being $1$ (a double fixed
point if $\phi \left( 1\right) =\phi ^{\prime }\left( 1\right) =1$)\emph{.}

BGW processes are therefore primarily concerned with\emph{\ }asexual
organisms who die while giving birth.\emph{\ }As such, BGW processes are
birth and death Markov processes with non-overlapping generations. However,
an individual dying while giving birth to a single offspring, so with
probability (w.p.) $\pi _{1}$, can also be interpreted as an individual
whose lifetime is delayed by one unit, thereby generating overlapping
generations. According to this new interpretation, each individual $i$ of a
BGW of some generation may equivalently (in law) be considered as splitting
after a geometric lifetime with success probability $1-\pi _{1}$ while
giving birth to $M_{i}$ offspring, given $M_{i}\neq 1$.

With $N_{n}\left( 1\right) $\ the number of individuals alive at generation $%
n$\ given $N_{0}\left( 1\right) =1$, we have 
\begin{equation*}
\phi _{n}\left( z\right) :=\mathbf{E}\left( z^{N_{n}\left( 1\right) }\right)
=\phi ^{\circ n}\left( z\right) ,
\end{equation*}
where $\phi ^{\circ n}\left( z\right) $\ is the $n$-th composition of $\phi
\left( z\right) $\ with itself, \footnote{%
Throughout this work, a p.g.f. will therefore be a function $\phi $ which is
absolutely monotone on $\left( 0,1\right) $ with all nonnegative derivatives
of any order there, obeying $\phi \left( 1\right) \leq 1.$ The defective
case $\phi \left( 1\right) <1$ will appear only marginally$.$}.
Equivalently, $\phi _{n}\left( z\right) $ obeys (as from a recursion from
the root) the branching property 
\begin{equation}
\phi _{n+1}\left( z\right) =\phi \left( \phi _{n}\left( z\right) \right)
=\phi _{n}\left( \phi \left( z\right) \right) \text{, }\phi _{0}\left(
z\right) =z.  \label{1.1}
\end{equation}
Similarly, if $N_{n}\left( i\right) $\ is the number of individuals alive at
generation $n$\ given there are $N_{0}=i$\ independent founders, we clearly
get 
\begin{equation}
\mathbf{E}\left( z^{N_{n}\left( i\right) }\right) :=\mathbf{E}\left(
z^{N_{n}}\mid N_{0}=i\right) =\phi _{n}\left( z\right) ^{i}.  \label{1.1.a}
\end{equation}
We shall also let 
\begin{equation}
\tau _{i,j}=\inf \left( n\geq 1:N_{n}=j\mid N_{0}=i\right) ,  \label{1.1.b}
\end{equation}
the first hitting time of state $j\neq i$\ given $N_{0}=i\neq 0$.\newline

If $\phi \left( 1\right) =1$ (the regular case), depending on $\mu :=\mathbf{%
E}\left( M\right) \leq 1$ (i.e. the (sub-)critical case ) or $\mu >1$
(supercritical case): the process $N_{n}\left( 1\right) $ goes eventually
extinct with probability $1$ or goes eventually extinct with probability $%
\rho _{e}<1$ where $\rho _{e}$ is the smallest fixed point solution in $%
\left[ 0,1\right] $ to $\phi \left( z\right) =z$, respectively (state $%
\left\{ 0\right\} $ is absorbing). In the latter case, the distribution of
the time to extinction $\tau _{1,0}$ is given by $\mathbf{P}\left( \tau
_{1,0}\leq 0\right) =0$ and 
\begin{equation*}
\mathbf{P}\left( \tau _{1,0}\leq n\right) =\mathbf{P}\left( N_{n}\left(
1\right) =0\right) =\phi _{n}\left( 0\right) ,\text{ }n\geq 1,
\end{equation*}
and the process explodes with complementary probability $\overline{\rho }%
_{e}:=1-\rho _{e}$, but not in finite time and $\tau _{1,0}=\infty $.
Clearly also, if there are $i$ independent founders instead of simply $1$, 
\begin{equation*}
\mathbf{P}\left( \tau _{i,0}\leq n\right) =\mathbf{P}\left( N_{n}\left(
i\right) =0\right) =\phi _{n}\left( 0\right) ^{i}.
\end{equation*}
Note that $\tau _{1,0}\overset{d}{=}H\left( 1\right) $ where $H\left(
1\right) $ is the height of $\left\{ N_{n}\left( 1\right) \right\} $.
Similarly, $\tau _{i,0}\overset{d}{=}H\left( i\right) $ where $H\left(
i\right) $ is the height of $\left\{ N_{n}\left( i\right) \right\} .$\newline

\textbf{Remark (}$N_{n}\left( 1\right) $ is positively correlated\textbf{):}
With $n_{1},n>0,$ letting $\phi _{n_{1},n_{1}+n}\left( z_{1},z_{2}\right) :=%
\mathbf{E}\left( z_{1}^{N_{n_{1}}\left( 1\right) }z_{2}^{N_{n_{1}+n}\left(
1\right) }\right) ,$ the branching property states that 
\begin{equation*}
\phi _{n_{1},n_{1}+n}\left( z_{1},z_{2}\right) =\phi _{n_{1}}\left(
z_{1}\phi _{n}\left( z_{2}\right) \right) .
\end{equation*}
Differentiating twice with respect to $z_{1}$ and then $z_{2}$ and
evaluating the result at $\left( 1,1\right) $ yields the (non-stationary)
autocovariance$:$%
\begin{eqnarray}
\text{Cov}\left( N_{n_{1}}\left( 1\right) ,N_{n_{1}+n}\left( 1\right)
\right) &=&\mathbf{E}\left( N_{n}\left( 1\right) \right) \mathbf{E}\left(
N_{n_{1}}\left( 1\right) ^{2}\right) -\mathbf{E}\left( N_{n_{1}}\left(
1\right) \right) \mathbf{E}\left( N_{n_{1}+n}\left( 1\right) \right)  \notag
\\
&=&\sigma ^{2}\left( N_{n_{1}}\left( 1\right) \right) \mu ^{n}>0\text{ if }%
\mu \neq 1  \notag \\
&=&\sigma ^{2}\left( N_{n_{1}}\left( 1\right) \right) =n_{1}\sigma ^{2}>0%
\text{ if }\mu =1
\end{eqnarray}
where $\sigma ^{2}=\sigma ^{2}\left( M\right) $ is the variance of $M$ and $%
\sigma ^{2}\left( N_{n_{1}}\left( 1\right) \right) =\sigma ^{2}\frac{\mu
^{n_{1}-1}\left( \mu ^{n_{1}}-1\right) }{\mu -1}$ the variance of $%
N_{n_{1}}\left( 1\right) $. The autocorrelation follows as 
\begin{equation*}
\text{Corr}\left( N_{n_{1}}\left( 1\right) ,N_{n_{1}+n}\left( 1\right)
\right) =\left\{ 
\begin{array}{c}
\frac{\text{Cov}\left( N_{n_{1}}\left( 1\right) ,N_{n_{1}+n}\left( 1\right)
\right) }{\sigma \left( N_{n_{1}}\left( 1\right) \right) \sigma \left(
N_{n_{1}+n}\left( 1\right) \right) } \\ 
\frac{\sqrt{\mu ^{n}}\sqrt{\left| \mu ^{n_{1}}-1\right| }}{\sqrt{\left| \mu
^{n_{1}+n}-1\right| }}\in \left( 0,1\right) \text{ if }\mu \neq 1 \\ 
=\left( 1+n/n_{1}\right) ^{-1/2}\in \left( 0,1\right) \text{ if }\mu \neq 1.%
\text{ }
\end{array}
\right.
\end{equation*}
Note $\lim_{n_{1}\rightarrow \infty }$Corr$\left( N_{n_{1}}\left( 1\right)
,N_{n_{1}+n}\left( 1\right) \right) =1$ ($=\mu ^{n/2}$) if $\mu \geq 1$ ($%
\mu <1$)$.$ $\blacksquare $

\subsection{Some limit laws in the regular case $\phi \left( 1\right) =1$
(see [Harris, 1963])}

We recall that $\phi \left( z\right) $ is assumed to have a convergence
radius $z_{*}\in \left( 1,\infty \right] .$

- subcritical case $\mu =\phi ^{\prime }\left( 1\right) <1$ ($\rho =1$)$:$
as $n\rightarrow \infty $%
\begin{equation}
N_{n}\left( 1\right) \mid N_{n}\left( 1\right) >0\overset{d}{\rightarrow }%
N_{\infty },  \label{f1}
\end{equation}
with 
\begin{equation}
\mathbf{E}\left( z^{N_{n}\left( 1\right) \mid N_{n}\left( 1\right)
>0}\right) =\frac{\phi _{n}\left( z\right) -\phi _{n}\left( 0\right) }{%
1-\phi _{n}\left( 0\right) }\rightarrow \phi _{\infty }\left( z\right) =%
\mathbf{E}\left( z^{N_{\infty }}\right) ,  \label{f2}
\end{equation}
solving the Schr\"{o}der functional equation [Hoppe, (1980)]: 
\begin{equation}
1-\phi _{\infty }\left( \phi \left( z\right) \right) =\mu \left( 1-\phi
_{\infty }\left( z\right) \right) .  \label{SFE}
\end{equation}

- critical case $\mu =\phi ^{\prime }\left( 1\right) =1$ ($\rho =1$)$:$ as $%
n\rightarrow \infty $%
\begin{equation}
\frac{N_{n}\left( 1\right) }{n}\mid N_{n}\left( 1\right) >0\overset{d}{%
\rightarrow }E,  \label{f3}
\end{equation}
where $E\overset{d}{\sim }$Exp$\left( 1\right) .$

- supercritical case $\mu =\phi ^{\prime }\left( 1\right) >1$ ($\rho <1$)$:$
as $n\rightarrow \infty $%
\begin{equation}
\mu ^{-n}N_{n}\left( 1\right) \overset{d}{\rightarrow }W\geq 0,  \label{f4}
\end{equation}
where the Laplace-Stieltjes transform (LST) $\phi _{W}\left( \lambda \right)
=\mathbf{E}e^{-\lambda W}$ of $W$ solves the Poincar\'{e}-Abel functional
equation 
\begin{equation}
\phi _{W}\left( \mu \lambda \right) =\phi \left( \phi _{W}\left( \lambda
\right) \right) ,  \label{f5}
\end{equation}
having mass $\rho $ at $W=0$. With complementary probability $\overline{\rho 
}=1-\rho $, the support of $W$ is the half-line. This results from the fact
that $\mu ^{-n}N_{n}\left( 1\right) $ is a non-negative martingale.

\subsection{The transition matrix approach and conditionings via Doob's
transforms}

A Bienaym\'{e}-Galton-Watson process is a time-homogeneous Markov chain with
denumerable state-space $\Bbb{N}_{0}:=\left\{ 0,1,...\right\} ;$ [see Woess,
(2009)]. Its stochastic irreducible transition matrix is $P$, with entries $%
P\left( i,j\right) =\left[ z^{j}\right] \phi \left( z\right) ^{i}=\mathbf{P}%
\left( N_{1}\left( i\right) =j\right) $ (with $\left[ z^{j}\right] \phi
\left( z\right) ^{i}$\ meaning the $z^{j}$-coefficient of the p.g.f. $\phi
\left( z\right) ^{i}$). State $\left\{ 0\right\} $ is absorbing and so $%
P\left( 0,j\right) =\delta _{0,j}$. When there is explosion as in the
supercritical cases, an interesting problem arises when conditioning $%
\left\{ N_{n}\right\} $ either on extinction or on explosion. This may be
understood by transformations of paths as follows:\newline

- \emph{Regular supercritical BGW process conditioned on extinction:} The
harmonic column vector $\mathbf{h}$, solution to $P\mathbf{h}=\mathbf{h}$,
is given by its coordinates $h\left( i\right) =\rho _{e}^{i}$, $i\geq 0$,
because $\sum_{j\geq 0}P\left( i,j\right) \rho _{e}^{j}=\phi \left( \rho
_{e}\right) ^{i}=\rho _{e}^{i}$. Letting $D_{\mathbf{h}}:=$diag$\left(
h\left( 0\right) ,h\left( 1\right) ,...\right) $, introduce the stochastic
matrix $P_{\mathbf{h}}$ given by a Doob transform [Norris (1998) and Rogers
and Williams (1994), p. $327$)]: $P_{\mathbf{h}}=D_{\mathbf{h}}^{-1}PD_{%
\mathbf{h}}$ or $P_{\mathbf{h}}\left( i,j\right) =h\left( i\right)
^{-1}P\left( i,j\right) h\left( j\right) =P\left( i,j\right) \rho _{e}^{j-i}$%
, $i,j\geq 0$. Note $h\left( N_{n}\left( i\right) \right) =\rho
_{e}^{N_{n}\left( i\right) }$ is a martingale because $\mathbf{E}\left(
h\left( N_{n}\left( i\right) \right) \right) =\phi _{n}\left( \rho
_{e}\right) ^{i}=\rho _{e}^{i}=h\left( i\right) =h\left( N_{0}\left(
i\right) \right) $. Then $P_{\mathbf{h}}$ is the transition matrix of $%
N_{1}\left( i\right) $ conditioned on almost sure extinction, with $P_{%
\mathbf{h}}^{n}\left( i,j\right) =h\left( i\right) ^{-1}P^{n}\left(
i,j\right) h\left( j\right) =\rho _{e}^{j-1}P^{n}\left( i,j\right) $ giving
the $n$-step transition matrix of the conditioned process. Equivalently,
when conditioning $N_{n}\left( 1\right) $ on almost sure extinction, one is
led to a regular subcritical BGW process with modified Harris-Sevastyanov
branching mechanism $\widetilde{\phi }_{0}\left( z\right) :=\phi \left( \rho
_{e}z\right) /\rho _{e}$, satisfying $\widetilde{\phi }_{0}\left( 1\right)
=1 $ and $\widetilde{\phi }_{0}^{\prime }\left( 1\right) =\phi ^{\prime
}\left( \rho _{e}\right) <1$. Indeed, $\widetilde{\phi }_{0}\left( z\right)
=\sum_{j\geq 0}P_{\mathbf{h}}\left( 1,j\right) z^{j}$. Upon iterating, we
get the composition rule $\widetilde{\phi }_{0,n}\left( z\right) =\phi
_{n}\left( \rho _{e}z\right) /\rho _{e}$. See [Klebaner et al. (2007), pp.
47-53]. \newline

- \emph{Regular supercritical BGW process conditioned on almost sure
explosion:}\textbf{\ }Similarly, when conditioning $\left\{ N_{n}\left(
1\right) \right\} $ on almost sure explosion, one is led to an explosive
supercritical BGW process with new Harris-Sevastyanov branching mechanism $%
\widetilde{\phi }_{\infty }\left( z\right) :=\left[ \phi \left( \rho _{e}+%
\overline{\rho }_{e}z\right) -\rho _{e}\right] /\overline{\rho }_{e}$,
satisfying $\widetilde{\phi }_{\infty }\left( 0\right) =0$ (all individuals
of the modified process are productive) and $\widetilde{\phi }_{\infty
}\left( 1\right) =\left( \phi \left( 1\right) -\rho _{e}\right) /\overline{%
\rho }_{e}=1$. Upon iterating, we get the composition rule $\widetilde{\phi }%
_{\infty }^{\circ n}\left( z\right) =\left[ \phi _{n}\left( \rho _{e}+%
\overline{\rho }_{e}z\right) -\rho _{e}\right] /\overline{\rho }_{e}$. With
probability $1$, this process drifts to $\infty $ in infinite time if $\phi
\left( 1\right) =1$.\newline

- \emph{Regular supercritical BGW process conditioned on never hitting }$%
\left\{ 0,\infty \right\} $\emph{:}\textbf{\ }BGW processes are unstable in
that they cannot reach a proper stationary distribution, being attracted
either at $\left\{ 0\right\} $ ($\mu \leq 1$) or at $\left\{ 0,\infty
\right\} $, ($\mu >1$). The following selection of paths reveals a proper
stationary measure in the supercritical case. Let $\overline{P}$ be a
substochastic matrix obtained from $P$ while removing its first row and
column. The largest eigenvalue (spectral radius) of $\overline{P}$ is $%
\gamma =\phi ^{\prime }\left( \rho _{e}\right) <1$. The corresponding
positive right (column) eigenvector $\mathbf{u}$ obeys $\overline{P}\mathbf{u%
}=\gamma \mathbf{u}$ with $u\left( i\right) =i\rho _{e}^{i-1}$, $i\geq 1$,
because $\sum_{j\geq 1}\overline{P}\left( i,j\right) j\rho _{e}^{j-1}=\phi
^{\prime }\left( \rho _{e}\right) i\phi \left( \rho _{e}\right)
^{i-1}=\gamma i\rho _{e}^{i-1}$. Conditioning $\left\{ N_{n}\left( 1\right)
\right\} $ on never hitting $\left\{ 0,\infty \right\} $ in the remote
future is given by the $Q$-process with stochastic transition matrix $%
Q=\gamma ^{-1}D_{\mathbf{u}}^{-1}\overline{P}D_{\mathbf{u}}$ or $Q\left(
i,j\right) =\gamma ^{-1}u\left( i\right) ^{-1}P\left( i,j\right) u\left(
j\right) =\gamma ^{-1}\rho _{e}^{j-i}i^{-1}\overline{P}\left( i,j\right) j$, 
$i,j\geq 1$ [see Lambert (2010) and Sagitov and Lindo (2015), Section $6$].
The modified Lamperti-Ney branching mechanism [Lamperti and Ney, (1968)] of
the $Q$-process has p.g.f. 
\begin{equation*}
\widetilde{\phi }_{Q}\left( z\right) :=\gamma ^{-1}\sum_{j\geq 1}\overline{P}%
\left( 1,j\right) j\rho _{e}^{j-1}z^{j}=z\phi ^{\prime }\left( z\rho
_{e}\right) /\phi ^{\prime }\left( \rho _{e}\right) .
\end{equation*}
The $Q$-process $\left\{ \widetilde{N}_{n}\left( 1\right) \right\} $ has an
invariant probability mass function (up to a normalization $K$) given by the
Hadamard product 
\begin{equation}
\mathbf{P}\left( \widetilde{N}_{\infty }\left( 1\right) =i\right) =Kv\left(
i\right) u\left( i\right) \text{, }i\geq 1,  \label{f6}
\end{equation}
where $\mathbf{v}$, with entries $v\left( i\right) $, obeys $\mathbf{v}%
^{\prime }\overline{P}=\gamma \mathbf{v}^{\prime },$ as a positive left
eigenvector\footnote{%
Here, a boldface variable, say $\mathbf{x}$, will represent a column-vector
so that its transpose, say $\mathbf{x}^{\prime }$, will be a row-vector.}.
Recalling $\overline{P}\left( i,j\right) =\left[ z^{j}\right] \phi \left(
z\right) ^{i}$, the generating function $v\left( z\right) =\sum_{i\geq
1}v\left( i\right) z^{i}$ of the $v\left( i\right) $'s obeys the Abel's
functional equation 
\begin{equation}
v\left( \phi \left( z\right) \right) -v\left( \phi \left( 0\right) \right)
=\gamma v\left( z\right) .  \label{f6a}
\end{equation}

Note that the Lamperti-Ney branching mechanism $\widetilde{\phi }_{Q}\left(
z\right) =z\phi ^{\prime }\left( z\rho _{e}\right) /\phi ^{\prime }\left(
\rho _{e}\right) $ is obtained as the composition of the Harris-Sevastyanov $%
\widetilde{\phi }_{0}\left( z\right) $ with the branching mechanism $\phi
_{SB}\left( z\right) =z\phi ^{\prime }\left( z\right) /\phi ^{\prime }\left(
\rho _{e}\right) $, the one of a size-biased version of $M$: $\widetilde{%
\phi }_{Q}\left( z\right) =\widetilde{\phi }_{0}\left( \phi _{SB}\left(
z\right) \right) $; see [Klebaner et al. (2007)]$.$\newline

- \emph{Regular supercritical or subcritical BGW processes conditioned to be
critical:}\textbf{\ }We end up with a last conditioning leading to a
critical BGW tree with mean offspring number $\mu _{c}=1.$ Let $\phi $,
regular, obey: $\phi $ has convergence radius $z_{*}>1$ (possibly $%
z_{*}=\infty $) and $\pi _{0}>0$. For such $\phi $'s, the unique positive
real root to the equation 
\begin{equation}
\phi \left( \tau \right) -\tau \phi ^{\prime }\left( \tau \right) =0,
\label{f8}
\end{equation}
exists, with $\rho _{e}=1<\tau <z_{*}$ if $\mu <1$ ($\phi \left( \tau
\right) >1$), $\tau =1$ if $\mu =1$ and $\rho _{e}<\tau <1<z_{*}$ if $\mu >1$
($\phi \left( \tau \right) <1$). In both cases, $\phi ^{\prime }\left( \tau
\right) <1.$\newline

Start with a supercritical branching process ($\mu >1$) and consider a
process whose modified branching mechanism is $\widetilde{\phi }_{c}\left(
z\right) =\phi \left( \tau z\right) /\phi \left( \tau \right) $, satisfying $%
\widetilde{\phi }_{c}\left( 1\right) =1$ and $\widetilde{\phi }_{c}^{\prime
}\left( 1\right) =:\mu _{c}=1$, the one of a critical branching process with
mean $1$ offspring distribution and variance: $\sigma _{c}^{2}=\tau ^{2}\phi
^{\prime \prime }\left( \tau \right) /\phi \left( \tau \right) $. Upon
iterating, we get the composition rule $\widetilde{\phi }_{c}^{\circ
n}\left( z\right) =\phi _{\tau }^{\circ n}\left( \tau z\right) /\tau $ where 
$\phi _{\tau }\left( z\right) =\phi \left( z\right) /\phi ^{\prime }\left(
\tau \right) $ is a scaled version of $\phi \left( z\right) $. Note $\phi
_{\tau }\left( 1\right) =1/\phi ^{\prime }\left( \tau \right) >1$ and $\phi
_{\tau }\left( \tau \right) =\tau $. The transition matrix $P_{c}$ of the
critical process is given by its entries 
\begin{equation*}
P_{c}\left( i,j\right) =\left[ z^{j}\right] \widetilde{\phi }_{c}\left(
z\right) ^{i}=\frac{\tau ^{j}}{\phi \left( \tau \right) ^{i}}P\left(
i,j\right) =\frac{\tau ^{j-i}}{\phi ^{\prime }\left( \tau \right) ^{i}}%
P\left( i,j\right) .
\end{equation*}
This transformation kills the supercritical paths to only select the
critical ones.

Similarly, starting with a subcritical branching process ($\mu <1$) and
considering a process whose modified branching mechanism (as a p.g.f.) is $%
\widetilde{\phi }_{c}\left( z\right) =\phi \left( \tau z\right) /\phi \left(
\tau \right) $, satisfying $\widetilde{\phi }_{c}\left( 1\right) =1$ and $%
\widetilde{\phi }_{c}^{\prime }\left( 1\right) =:\mu _{c}=1$, the one of a
critical branching process. Upon iterating, we get the composition rule $%
\widetilde{\phi }_{c}^{\circ n}\left( z\right) =\phi _{\tau }^{\circ
n}\left( \tau z\right) /\tau $ where $\phi _{\tau }\left( z\right) =\phi
\left( z\right) /\phi ^{\prime }\left( \tau \right) $ is a scaled version of 
$\phi \left( z\right) $. Note again $\phi _{\tau }\left( 1\right) =1/\phi
^{\prime }\left( \tau \right) >1$.

This transformation creates critical paths from the subcritical ones.\newline

The large-$n$ asymptotic properties of the above processes requires the
evaluation of the large-$n$ iterates of a p.g.f. There are classes of
discrete branching processes for which the $n$-step p.g.f. $\phi _{n}\left(
z\right) $\ of $N_{n}\left( 1\right) $\ (but also the `tilded' ones of their
conditioned versions) is exactly computable, thereby making the above
computations concrete and somehow explicit.

This is the case for the LF p.g.f. $\phi \left( z\right) =\pi _{0}+\overline{%
\pi }_{0}\frac{\pi z}{1-\overline{\pi }z}$\ for which, assuming $\rho
_{e}:=\pi _{0}/\overline{\pi }<1$\ (the super-criticality condition, see
below):

- Almost sure extinction 
\begin{equation*}
\bullet \text{ }\widetilde{\phi }_{0}\left( z\right) :=\phi \left( \rho
_{e}z\right) /\rho _{e}=\overline{\pi }+\pi \frac{\overline{\pi }_{0}z}{%
1-\pi _{0}z}.
\end{equation*}

- Immortal individuals 
\begin{equation*}
\bullet \text{ }\widetilde{\phi }_{\infty }\left( z\right) :=\left[ \phi
\left( \rho _{e}+\overline{\rho }_{e}z\right) -\rho _{e}\right] /\overline{%
\rho }_{e}=\frac{z\pi /\overline{\pi }_{0}}{1-\left( 1-\pi /\overline{\pi }%
_{0}\right) z}.
\end{equation*}

- $Q$-process 
\begin{equation*}
\bullet \text{ }\widetilde{\phi }_{Q}\left( z\right) =z\phi ^{\prime }\left(
z\rho _{e}\right) /\phi ^{\prime }\left( \rho _{e}\right) =z\left( \frac{%
\overline{\pi }_{0}}{1-\pi _{0}z}\right) ^{2}.
\end{equation*}
To compute the probability mass function (\ref{f6}), we first need to solve (%
\ref{f6a}), or equivalently 
\begin{equation*}
v\left( \phi \left( z\right) \right) -1=\gamma v\left( z\right) ,
\end{equation*}
while imposing $v\left( \pi _{0}\right) =1$\ and $v\left( 0\right) =0,$\ $%
v\left( \rho _{e}\right) =1/\left( 1-\gamma \right) >1$. Recall $\gamma
=\phi ^{\prime }\left( \rho _{e}\right) =\pi /\overline{\pi }_{0}<1$. The
solution (satisfying $v\left( \pi _{0}\right) =1$) is found to be 
\begin{equation*}
v\left( z\right) =1-\frac{1}{\log m}\log \left[ \frac{a-z}{a-\pi _{0}}\frac{%
1-\pi _{0}}{1-z}\right] ,
\end{equation*}
with $m>1$, $1>a>\rho _{e}>\pi _{0}.$\ It diverges at $z=a$. The condition $%
v\left( \rho _{e}\right) =1/\left( 1-\gamma \right) $\ yields 
\begin{equation*}
a=\frac{\rho _{e}\overline{\pi }_{0}-\overline{\rho }_{e}\pi _{0}m^{\gamma
/\left( 1-\gamma \right) }}{\overline{\pi }_{0}-\overline{\rho }_{e}\pi
_{0}m^{\gamma /\left( 1-\gamma \right) }}\in \left( \rho _{e},1\right) .
\end{equation*}
The condition $v\left( 0\right) =0$\ yields 
\begin{equation*}
m=\frac{a\overline{\pi }_{0}}{a-\pi _{0}}>1.
\end{equation*}
Consequently, with $\sum_{i\geq 1}v_{i}=\infty ,$\ the left eigenvector of $%
\overline{P}$\ associated to the eigenvalue $\gamma $\ is given by 
\begin{equation*}
v_{i}=\left[ z^{i}\right] v\left( z\right) =\frac{1}{\log m}\frac{a^{-i}-1}{i%
}\text{, }i\geq 1
\end{equation*}
and, up to the finite normalization factor $K=\left( a-\rho _{e}\right) 
\overline{\rho }_{e}\log m/\left( 1-a\right) $, with $i\geq 1,$\ recalling $%
u\left( i\right) =i\rho _{e}^{i-1},$%
\begin{equation*}
\bullet \text{ }\mathbf{P}\left( \widetilde{N}_{\infty }\left( 1\right)
=i\right) =Kv\left( i\right) u\left( i\right) =\frac{K}{\log m}\left(
a^{-i}-1\right) \rho _{e}^{i-1}=\frac{\left( a-\rho _{e}\right) \overline{%
\rho }_{e}}{\rho _{e}\left( 1-a\right) }\left[ \left( \rho _{e}/a\right)
^{i}-\rho _{e}^{i}\right] ,\text{ }
\end{equation*}
is the explicit invariant probability mass of this $Q$-process. It decays
asymptotically geometrically at rate $\rho _{e}/a.$

- Forced criticality: with 
\begin{equation*}
\bullet \text{ }\tau =\frac{-\pi _{0}\overline{\pi }+\sqrt{\pi _{0}\overline{%
\pi }_{0}\pi \overline{\pi }}}{\overline{\pi }\left( \pi -\pi _{0}\right) }
\end{equation*}
\ the explicit radical solution to the quadratic equation $\phi \left( \tau
\right) -\tau \phi ^{\prime }\left( \tau \right) =0$\ in the LF case, 
\begin{eqnarray*}
\widetilde{\phi }_{c}\left( z\right) &=&\phi \left( \tau z\right) /\phi
\left( \tau \right) =\frac{1}{\phi \left( \tau \right) }\left( \pi _{0}+%
\overline{\pi }_{0}\frac{\pi \tau z}{1-\overline{\pi }\tau z}\right) \\
&=&P_{0}+\overline{P}_{0}\frac{Pz}{1-\overline{P}z}\text{ with }P_{0}=\frac{%
\pi _{0}}{\phi \left( \tau \right) }\text{ and }\overline{P}=\overline{\pi }%
\tau .
\end{eqnarray*}

General branching processes conditioned as above are still branching
processes. Except for the negative binomial [NB$\left( 2,\pi _{0}\right) $]
p.g.f. $\widetilde{\phi }_{Q}\left( z\right) $, the `tilded' p.g.f.'s of LF
branching mechanisms are again LF ones. Iterating such `tilded' LF p.g.f.'s
yield again LF mechanisms.

\subsection{The linear-fractional model}

We shall deal with the following regular LF case with two parameters $\pi
_{0},\pi \in \left( 0,1\right) $ (unless otherwise specified), as a
zero-inflated geometric p.g.f.: 
\begin{equation}
\bullet \text{ }\phi \left( z\right) =\pi _{0}+\overline{\pi }_{0}\frac{\pi z%
}{1-\overline{\pi }z}=1-\frac{1}{\overline{\pi }/\overline{\pi }_{0}+\pi /%
\overline{\pi }_{0}\left( 1-z\right) ^{-1}},  \label{Mlaw}
\end{equation}
for which $\pi _{m}=\mathbf{P}\left( M=m\right) =\overline{\pi }_{0}\pi 
\overline{\pi }^{m-1}$, $m\geq 1$ ($\overline{\pi }_{0}=1-\pi _{0}$ and $%
\overline{\pi }=1-\pi $)$.$ This distribution has mean $\mu :=\mathbf{E}%
\left( M\right) =\overline{\pi }_{0}/\pi $ and variance $\sigma ^{2}:=\sigma
^{2}\left( M\right) =\overline{\pi }_{0}\left( \overline{\pi }+\pi
_{0}^{{}}\right) /\pi ^{2}$ and, alternatively, 
\begin{equation*}
\bullet \text{ }\phi \left( z\right) =\frac{1+\left( 1-z\right) \left( 
\overline{\pi }-\overline{\pi }_{0}\right) /\pi }{1+\left( 1-z\right) 
\overline{\pi }/\pi }=\frac{\pi _{0}+z\left( \pi -\pi _{0}\right) }{1-z%
\overline{\pi }}.
\end{equation*}
[Athreya-Ney (1972), p. 22] suggest that one could bound an arbitrary
generating function $\phi $, $\phi ^{\prime \prime }\left( 1\right) <\infty $
between two LF generating functions. Linear-fractional branching mechanisms
is one of some rare p.g.f.'s which is stable under composition [Sagitov and
Lindo (2015); Grosjean and Huillet (2017)].

For this reproduction model, 
\begin{equation*}
\mu =\phi ^{\prime }\left( 1\right) =\overline{\pi }_{0}/\pi .
\end{equation*}
The non-trivial ($\neq 1$) solution to $\phi \left( \rho \right) =\rho $ is 
\begin{equation}
\bullet \text{ }\rho =\pi _{0}/\overline{\pi },  \label{rho}
\end{equation}
with $\rho =:\rho _{e}<1$, the extinction probability, if $\mu >1$ ($a<1$).
If $\mu >1$ ($\mu <1$), this BGW process is supercritical (subcritical, with 
$\rho >1$). It is critical when $\mu =1$.

It has mode at the origin if and only if $\pi _{0}>\pi _{1}$, else $\pi
_{0}>\pi /\left( 1+\pi \right) $. Otherwise, the mode is at $1$. Note here $%
\pi _{1}=\overline{\pi }_{0}\pi $, relevant in the non-overlapping
interpretation of this process. Given some individual produces offspring
(with probability $\overline{\pi }_{0}$), the number of offspring is
geometrically distributed with success probability $\pi $. This branching
mechanism model was considered by [Steffenson, (1930, 1933)] in the
extinction of family surnames problem; see [Kendall, (1966)] for historical
background. In the 1920 United-States census of white males with $\rho
_{e}\sim 0.860$ as the probability of the termination of the male line of
descent from a new-born male, the data fits the facts fairly well using $%
\widehat{\pi }_{0}$ $=0.481$ and $\widehat{\overline{\pi }}=0.559$ ($\mu
=1.163>1$; $\sigma =1.633$). The mode is at the origin. The probability of
having more than $l$ offspring is $\mathbf{P}\left( M>l\right) =\overline{%
\pi }_{0}\pi ^{l+1}/\overline{\pi }$ which, for these values of $\left( \pi
_{0},\overline{\pi }\right) $, yields $0.0325$ if $l=3.$\newline

\textbf{Remark }(geometric infinite-divisibility): Let $M^{\prime }$\ be a
random variable (r.v.) obtained as a Geo$\left( \nu \right) $\ sum of i.i.d.
Bernoulli$\left( p\right) $\ r.v.'s, $p,\nu \in \left( 0,1\right) $. Its
p.g.f. reads 
\begin{equation*}
\phi _{M^{\prime }}^{{}}\left( z\right) :=\mathbf{E}z^{M^{\prime }}=\frac{%
\nu \left( q+pz\right) }{1-\overline{\nu }\left( q+pz\right) }=\frac{\nu q}{%
1-\overline{\nu }q}+\frac{p}{1-\overline{\nu }q}\frac{\frac{\nu }{1-%
\overline{\nu }q}z}{1-\frac{\overline{\nu }p}{1-\overline{\nu }q}z}.
\end{equation*}
It can be put under the form (\ref{Mlaw}) if 
\begin{eqnarray*}
\pi _{0} &=&\frac{\nu q}{1-\overline{\nu }q}\text{ and }\pi =\frac{\nu }{1-%
\overline{\nu }q} \\
q &=&\frac{\pi _{0}}{\pi }\text{ and }\nu =\frac{\pi -\pi _{0}}{\overline{%
\pi }_{0}},
\end{eqnarray*}
so only if $\pi _{0}<\pi $. Under this condition therefore, is the r.v. $M$\
whose law is defined in (\ref{Mlaw}) interprets as a Bernoulli-thinning of a
Geo$\left( \nu \right) -$distributed r.v.. This will be the case if $\pi
_{0}<\pi /\left( 1+\pi \right) <\pi $\ (when the mode of $M$\ is away from $%
0 $\ at $1$). $M^{\prime }$\ is easily shown to be infinitely-divisible
(else compound-Poisson) with clusters' having Fisher's log-series
distribution.

When $\pi _{0}>\pi ,$\ the general LF p.g.f. can be put under the compound
Geo$_{0}$\ form 
\begin{equation*}
\phi \left( z\right) =\frac{\pi _{0}-z\left( \pi _{0}-\pi \right) }{1-z%
\overline{\pi }}=\frac{\pi }{1-\overline{\pi }\psi \left( z\right) },
\end{equation*}
for some well-defined LF clusters' p.g.f. 
\begin{equation*}
\psi \left( z\right) =\frac{1}{\overline{\pi }}\frac{\pi _{0}-\pi -z\left(
\pi _{0}-\pi -\pi \overline{\pi }\right) }{\pi _{0}-z\left( \pi _{0}-\pi
\right) }.
\end{equation*}
In that case, $M$\ is an independent random Geo$_{0}\left( \pi \right) $\
sum of i.i.d. clusters with LF sizes. It is thus geometrically-infinitely
divisible (hence infinitely divisible or compound-Poisson). $\blacksquare $%
\newline

From (\ref{Mlaw}), $\phi \left( z\right) =\left( \alpha z+\beta \right)
/\left( \gamma z+\delta \right) $ is an homography (M\"{o}bius transform)
encoded by the matrix 
\begin{equation*}
A=\left[ 
\begin{array}{ll}
\alpha & \beta \\ 
\gamma & \delta
\end{array}
\right] =\left[ 
\begin{array}{ll}
\pi -\pi _{0} & \pi _{0} \\ 
-\overline{\pi } & 1
\end{array}
\right] \text{.}
\end{equation*}
$A$ is invertible because $\left| A\right| :=\alpha \delta -\beta \gamma =%
\overline{\pi }_{0}\pi =\pi _{1}\neq 0$. Diagonalization of $A$ (with row
sum $\pi $) yields $A=SDS^{-1}$ with 
\begin{equation*}
S=\left[ 
\begin{array}{ll}
1 & 1 \\ 
1 & \overline{\pi }/\pi _{0}
\end{array}
\right] \text{, }S^{-1}=\frac{1}{\overline{\pi }-\pi _{0}}\left[ 
\begin{array}{ll}
\overline{\pi } & -\pi _{0} \\ 
-\pi _{0} & \pi _{0}
\end{array}
\right] \text{, }D=\left[ 
\begin{array}{ll}
\pi & 0 \\ 
0 & \overline{\pi }_{0}
\end{array}
\right]
\end{equation*}
and 
\begin{equation*}
A^{n}=\frac{1}{\overline{\pi }-\pi _{0}}\left[ 
\begin{array}{ll}
\overline{\pi }\pi ^{n}-\pi _{0}\overline{\pi }_{0}^{n} & \pi _{0}\overline{%
\pi }_{0}^{n}-\pi _{0}\pi ^{n} \\ 
\overline{\pi }\pi ^{n}-\overline{\pi }\overline{\pi }_{0}^{n} & \overline{%
\pi }\overline{\pi }_{0}^{n}-\pi _{0}\pi ^{n}
\end{array}
\right] =\left[ 
\begin{array}{ll}
\alpha _{n} & \beta _{n} \\ 
\gamma _{n} & \delta _{n}
\end{array}
\right]
\end{equation*}
is the homography-matrix associated to $\phi _{n}\left( z\right) =\phi
^{\circ n}\left( z\right) =\left( \alpha _{n}z+\beta _{n}\right) /\left(
\gamma _{n}z+\delta _{n}\right) $. The matrix $A^{n}$ has row sum $\pi ^{n}$
with $\alpha _{n}+\beta _{n}=\gamma _{n}+\delta _{n}=\pi ^{n}$ translating
that $\phi _{n}\left( z\right) $ is a p.g.f.. Among the sequences $\left(
\alpha _{n},\beta _{n},\gamma _{n},\delta _{n}\right) $, only two of them
are therefore independent.

An alternative representation of $\phi _{n}$ is 
\begin{equation}
\bullet \text{ }\phi _{n}\left( z\right) =1-\frac{1}{b_{n}+a_{n}\left(
1-z\right) ^{-1}},  \label{R2}
\end{equation}
with, by recurrence, 
\begin{equation}
a_{n}=a^{n}\text{, }b_{n}=b\left( 1+a+...+a^{n-1}\right) =b\frac{a^{n}-1}{a-1%
},\text{ if }a\neq 1,  \label{R3}
\end{equation}
\begin{equation*}
a_{n}=1\text{, }b_{n}=bn\text{, if }a=1\text{ (critical case)}
\end{equation*}
and $a=\pi /\overline{\pi }_{0}=1/\mu $ and $b=\overline{\pi }/\overline{\pi 
}_{0}$ ($a+b=1/\overline{\pi }_{0}$, $a_{n}+b_{n}=\left( a^{n}\pi _{0}/%
\overline{\pi }_{0}-b\right) /\left( a-1\right) >0$).

It can be checked that the following relations between $\left(
a_{n},b_{n}\right) $ and $\left( \alpha _{n},\beta _{n},\gamma _{n},\delta
_{n}\right) $ hold: 
\begin{eqnarray*}
\alpha _{n} &=&\overline{\pi }_{0}^{n}\left( 1-b_{n}\right) ;\text{ }\beta
_{n}=\overline{\pi }_{0}^{n}\left( a_{n}+b_{n}-1\right) , \\
\gamma _{n} &=&-\overline{\pi }_{0}^{n}b_{n};\text{ }\delta _{n}=\overline{%
\pi }_{0}^{n}\left( a_{n}+b_{n}\right) ,
\end{eqnarray*}
\begin{equation*}
a_{n}=\overline{\pi }_{0}^{-n}\left( \alpha _{n}+\beta _{n}\right) =\left(
\pi /\overline{\pi }_{0}\right) ^{n}\text{; }b_{n}=-\overline{\pi }%
_{0}^{-n}\gamma _{n}.
\end{equation*}

From the second representation of $\phi _{n}\left( z\right) ,$ 
\begin{equation*}
\mathbf{P}\left( N_{n}\left( 1\right) >0\right) =1-\phi _{n}\left( 0\right) =%
\mathbf{P}\left( \tau _{1,0}>n\right) =1/\left( b_{n}+a_{n}\right)
\end{equation*}
and 
\begin{eqnarray*}
\phi _{n}^{\prime }\left( 1\right) &=&\mathbf{E}N_{n}\left( 1\right) =\frac{1%
}{a_{n}}=\mu ^{n},\text{ } \\
\sigma ^{2}\left( N_{n}\left( 1\right) \right) &=&\frac{2b_{n}+a_{n}-1}{%
a_{n}^{2}}=\sigma ^{2}\frac{\mu ^{n-1}\left( \mu ^{n}-1\right) }{\mu -1}.
\end{eqnarray*}
The probability of non-extinction (survival) at generation $n$\ is: 
\begin{eqnarray*}
\bullet \text{ }\mathbf{P}\left( \tau _{1,0}>n\right) &=&1/\left(
b_{n}+a_{n}\right) \sim a^{-n}/\left( 1+b/\left( a-1\right) \right) \text{
if }a>1 \\
&=&\left( 1-\frac{\overline{\pi }}{\pi _{0}}\right) \mu ^{n}\text{ (}\mu <1%
\text{) (subcritical regime).}
\end{eqnarray*}
\begin{eqnarray*}
\bullet \text{ }\mathbf{P}\left( \tau _{1,0}>n\right) &=&1/\left(
b_{n}+a_{n}\right) \sim 1/\left( bn\right) \text{ if }a=1 \\
&=&\frac{\pi }{\overline{\pi }}n^{-1}\text{ (}\mu =1\text{) (critical
regime).}
\end{eqnarray*}
\begin{eqnarray*}
\bullet \text{ }\mathbf{P}\left( \tau _{1,0}>n\right) &=&\left( \left(
1-a\right) /b\right) /\left[ 1-a^{n}\left( 1-\left( 1-a\right) /b\right)
\right] \text{ if }a<1 \\
&=&\overline{\rho }_{e}/\left[ 1-\rho _{e}\mu ^{-n}\right] \text{ (}\mu >1%
\text{) (supercritical regime).}
\end{eqnarray*}
The time to extinction of the subcritical LF BGW has geometric tails (rapid
extinction), whereas the time to extinction of the critical LF BGW has
power-law tails with index $1$\ (slow extinction). In both cases, extinction
is almost sure ($\rho _{e}=1$).

In the supercritical regime, $\mathbf{P}\left( \tau _{1,0}>n\right)
\rightarrow \left( 1-a\right) /b=\overline{\rho }_{e}$, the first-order
correcting term being geometrically small.

From the exact expression of $\mathbf{P}\left( \tau _{1,0}>n\right) $, as in
[Garcia-Millan R. et al., 2015 and Corral \'{A}. et al., 2016], we observe
the following finite-size scaling law in the \textit{slightly supercritical}
regime for which $\mu =1+x/n,$\ $x>0$\ and $\rho _{e}\sim 1-2\left( \mu
-1\right) /\sigma _{c}^{2}$, [$\sigma _{c}^{2}=2\pi _{0}/\pi $, the critical
variance of $M$\ when $\mu =1$, see (\ref{e4a}) below:\ 
\begin{equation*}
\bullet \text{ }n\mathbf{P}\left( \tau _{1,0}>n\right) \rightarrow r\left(
x\right) :=\frac{1}{\sigma _{c}^{2}}\frac{2xe^{x}}{e^{x}-1}\text{ as }%
n\rightarrow \infty .
\end{equation*}
As in the strictly critical regime, the time to extinction has power-law
tails with index $1$, but with a non-constant asymptotic rate $r\left(
x\right) $.\newline

\textbf{Remark }(transition matrix powers): With $\pi _{n}\left( 0\right) =%
\mathbf{P}\left( N_{n}\left( 1\right) =0\right) =\mathbf{P}\left( \tau
_{1,0}\leq n\right) =\phi _{n}\left( 0\right) $, ($\overline{\pi }_{n}\left(
0\right) =1-\pi _{n}\left( 0\right) $), the p.g.f. $\phi _{n}\left( z\right) 
$\ at step $n$\ can be put under the form (similar to \ref{Mlaw} when $n=1$) 
\begin{equation*}
\phi _{n}\left( z\right) =\pi _{n}\left( 0\right) +\overline{\pi }_{n}\left(
0\right) g_{n}\left( z\right) ,
\end{equation*}
where $g_{n}\left( z\right) $\ is the p.g.f. of a geometric distribution
with failure probability $\overline{\pi }_{n}=1-\frac{a_{n}}{a_{n}+b_{n}}$.
The law of $N_{n}\left( 1\right) $\ has mode at the origin if and only if $%
\pi _{n}\left( 0\right) >\overline{\pi }_{n}\left( 0\right) \pi _{n}$, else $%
\overline{\pi }_{n}\left( 0\right) <1/\left( 1+\pi _{n}\right) $. Otherwise,
the mode is at $1$. The Fa\`{a}-di-Bruno formula [see Comtet, (1970)] allows
for an explicit expression of the step$-n$\ transition probability
(involving $i$\ founders) 
\begin{equation*}
P^{n}\left( i,j\right) =\mathbf{P}\left( N_{n}\left( i\right) =j\right)
=\left[ z^{j}\right] \phi _{n}\left( z\right) ^{i},
\end{equation*}
resulting from the composition of the binomial p.g.f. $\left( \pi _{n}\left(
0\right) +\overline{\pi }_{n}\left( 0\right) z\right) ^{i}$\ with a
geometric one $g_{n}\left( z\right) $. For all $n\geq 1$, we get $%
P^{n}\left( 0,j\right) =\delta _{0,j}$, $P^{n}\left( i,0\right) =\pi
_{n}\left( 0\right) ^{i}$, $i\geq 1$\ and 
\begin{equation*}
\bullet \text{ }P^{n}\left( i,j\right) =\overline{\pi }_{n}^{j}\sum_{k=1}^{i%
\wedge j}\binom{i}{k}\binom{j-1}{k-1}\left( \frac{\overline{\pi }_{n}\left(
0\right) \pi _{n}}{\overline{\pi }_{n}}\right) ^{k}\pi _{n}\left( 0\right)
^{i-k}\text{, }i,j\geq 1.
\end{equation*}
As a result, we obtained a closed-form expression of the Green kernel of the
LF model:

\begin{equation*}
G_{i,j}\left( u\right) :=\sum_{n\geq 0}u^{n}P^{n}\left( i,j\right) .
\end{equation*}
Note from the above expressions of the survival probabilities $\mathbf{P}%
\left( \tau _{1,0}>n\right) =1-\pi _{n}\left( 0\right) $\ that, whatever the
regime, 
\begin{equation*}
\bullet \text{ }G_{i,0}\left( 1\right) :=\sum_{n\geq 0}P^{n}\left(
i,0\right) =\sum_{n\geq 0}\pi _{n}\left( 0\right) ^{i}=\infty .
\end{equation*}
In particular, $G_{0,0}\left( 1\right) =\infty $\ and state $0$\ is visited
infinitely often. $\blacksquare $\newline

\textbf{Remark }(embedding)\textbf{:}\textit{\ }Let $f\left( z\right)
=A\left( z-1\right) +B/2\left( z-1\right) ^{2}$, $B>0$\ ($A=f^{\prime
}\left( 1\right) $, $B=f^{\prime \prime }\left( 1\right) $)$.$\ Consider the
p.g.f. $\phi _{t}\left( z\right) =\mathbf{E}\left( z^{N_{t}\left( 1\right)
}\right) $\ of a continuous-time branching process $N_{t}\left( 1\right) $\
with binary branching mechanism $f\left( z\right) .$\ We have $\partial
_{t}\phi _{t}\left( z\right) =f\left( \phi _{t}\left( z\right) \right) ;$\ $%
t\geq 0,$\ $\phi _{0}\left( z\right) =z$, whose solution when $A\neq 0$\ is 
\begin{equation*}
\bullet \text{ }\phi _{t}\left( z\right) =1-\frac{1}{\frac{B}{2A}\left(
1-e^{-At}\right) +e^{-At}\left( 1-z\right) ^{-1}},
\end{equation*}
with $\phi _{t}\left( z\right) \rightarrow 1-\frac{2A}{B}$\ as $t\rightarrow
\infty $\ if $A>0$\ (the supercritical case). Recall $\phi _{t+s}\left(
z\right) =\phi _{t}\left( \phi _{s}\left( z\right) \right) $, $s,t\geq 0$,
as a semi-group. With $A=-\log a$\ and $B=-2\left( b\log a\right) /(1-a)>0$\
and $t=n,$\ this is (\ref{R2}) showing that the discrete-time Markov chain
with LF branching mechanism is embeddable in the continuous-time branching
process with binary fission. $\blacksquare $\newline

Coming back to the discrete-time setting, we have:

- In the subcritical case when $\mu <1$ ($a>1$), with $a_{n}\rightarrow 0$
and $b_{n}/a_{n}\rightarrow b/\left( 1-a\right) $, we have 
\begin{eqnarray*}
\mathbf{E}\left( z^{N_{n}\left( 1\right) }\mid N_{n}\left( 1\right)
>0\right) &=&\frac{\phi _{n}\left( z\right) -\phi _{n}\left( 0\right) }{%
1-\phi _{n}\left( 0\right) }=\frac{z}{\left( 1-z\right) b_{n}/a_{n}+1} \\
&\rightarrow &\phi _{\infty }\left( z\right) =\frac{z}{1-\left( 1-z\right)
b/\left( a-1\right) },
\end{eqnarray*}
the p.g.f. of a geometric r.v. with mean $b/\left( a-1\right) =\overline{\pi 
}/\left( \pi -\overline{\pi }_{0}\right) $, solving the associated
Schr\"{o}der functional equation (\ref{SFE})$.$ The reciprocal of the mean
is $K=\left( \pi -\overline{\pi }_{0}\right) /\overline{\pi }$, the
Kolmogorov constant for which $\mu ^{-n}\mathbf{P}\left( \tau
_{1,0}>0\right) \underset{n\rightarrow \infty }{\rightarrow }K$. This
constant is thus explicit in the LF case.

- In the supercritical case when $\mu >1$ ($a<1$), with $\mu
^{-n}b_{n}/a_{n}=b_{n}\rightarrow b/\left( 1-a\right) ,$%
\begin{eqnarray*}
\mathbf{E}\left( e^{-\lambda \mu ^{-n}N_{n}\left( 1\right) }\mid N_{n}\left(
1\right) >0\right) &=&\frac{e^{-\lambda \mu ^{-n}}}{\left( 1-e^{-\lambda \mu
^{-n}}\right) b_{n}/a_{n}+1} \\
&\rightarrow &\frac{1}{1+\lambda b/\left( 1-a\right) },
\end{eqnarray*}
the LST of an exponential distribution with mean $b/\left( 1-a\right) =%
\overline{\pi }/\left( \overline{\pi }_{0}-\pi \right) $; see [Yaglom,
(1947)]$.$

- When $\pi =\pi _{0}$ ($b=1$), $\phi \left( z\right) =\frac{\pi _{0}}{1-%
\overline{\pi }_{0}z},$ the p.g.f. of a shifted to the left by one unit, say
Geo$_{0}\left( \pi _{0}\right) $, distribution. Here, $\mu =\phi ^{\prime
}\left( 1\right) =\overline{\pi }_{0}/\pi _{0}=1/a$ with $a$ a non-trivial
solution to $\phi \left( a\right) =a$. If $\mu >1$ ($\pi _{0}<1/2$), the BGW
is supercritical and $a=\pi _{0}/\overline{\pi }_{0}=\rho _{e}<1$, the
extinction probability. If $\mu <1$ ($\pi _{0}>1/2$), this BGW is
subcritical and 
\begin{equation*}
\mathbf{E}\left( z^{N_{n}\left( 1\right) }\mid N_{n}\left( 1\right)
>0\right) \rightarrow \frac{z}{1-\left( 1-z\right) /\left( a-1\right) },
\end{equation*}
the p.g.f. of a geometric r.v. with mean $1/\left( a-1\right) =\overline{\pi 
}_{0}/\left( 2\pi _{0}-1\right) .$

- In the critical case when $\overline{\pi }=\pi _{0},$ with $\phi \left(
z\right) =\frac{\left( 1-2\overline{\pi }\right) z+\overline{\pi }}{1-%
\overline{\pi }z},$ the matrix $A$ is 
\begin{equation*}
A=\left[ 
\begin{array}{ll}
\alpha & \beta \\ 
\gamma & \delta
\end{array}
\right] =\left[ 
\begin{array}{ll}
\pi -\overline{\pi } & \overline{\pi } \\ 
-\overline{\pi } & 1
\end{array}
\right] ,
\end{equation*}
with the double eigenvalue $\pi $. We get 
\begin{equation*}
A^{n}=\left[ 
\begin{array}{ll}
\alpha _{n} & \beta _{n} \\ 
\gamma _{n} & \delta _{n}
\end{array}
\right] =\pi ^{n-1}\left[ 
\begin{array}{ll}
\pi -n\overline{\pi } & n\overline{\pi } \\ 
-n\overline{\pi } & \pi +n\overline{\pi }
\end{array}
\right] ,
\end{equation*}
still with row sum $\pi ^{n}.$ Alternatively, 
\begin{equation*}
\phi _{n}\left( z\right) =1-\frac{1}{b_{n}+a_{n}\left( 1-z\right) ^{-1}},
\end{equation*}
with 
\begin{equation*}
a_{n}=1\text{; }b_{n}=-\overline{\pi }_{0}^{-n}\gamma _{n}=n\overline{\pi }%
/\pi .
\end{equation*}
Note 
\begin{equation*}
1-\phi _{n}\left( 0\right) =\mathbf{P}\left( N_{n}\left( 1\right) >0\right) =%
\mathbf{P}\left( \tau _{1,0}>n\right) =1/\left( b_{n}+a_{n}\right) =1/\left(
1+n\overline{\pi }/\pi \right) .
\end{equation*}
In the critical case, there is almost sure (a.s.) extinction but the time to
extinction is slow with power-law tails of order $1/n$.

- When $\pi _{0}=0$ (immortal individuals), 
\begin{equation*}
\phi \left( z\right) =\frac{\pi z}{1-\overline{\pi }z}=1-\frac{1}{\overline{%
\pi }+\pi \left( 1-z\right) ^{-1}},
\end{equation*}
the p.g.f. of a proper geometric distribution with failure probability $\pi $%
. In that case, 
\begin{equation*}
\phi _{n}\left( z\right) =\phi ^{\circ n}\left( z\right) =1-\frac{1}{%
b_{n}+a_{n}\left( 1-z\right) ^{-1}},
\end{equation*}
with 
\begin{equation*}
a_{n}=a^{n}\text{, }b_{n}=b\left( 1+a+...+a^{n-1}\right) ,
\end{equation*}
and $a=\pi $ and $b=\overline{\pi }.$ Here, $\mu =\phi ^{\prime }\left(
1\right) =1/\pi >1$ and $\phi \left( 0\right) =0:$ the model is strictly
supercritical, with $N_{n}\left( 1\right) >0\overset{d}{\rightarrow }\infty $%
, corresponding to explosion with probability $1$ ($\rho _{e}=0$). Note
indeed $\phi _{n}\left( 0\right) =\mathbf{P}\left( N_{n}\left( 1\right)
=0\right) =0$, for all $n\geq 0.$

The representation (\ref{R2}) of $\phi _{n}\left( z\right) $ is useful
because 
\begin{eqnarray*}
\frac{1-\phi _{n}\left( z\right) }{1-z} &=&\frac{1}{a_{n}+b_{n}\left(
1-z\right) }=\sum_{k\geq 0}\mathbf{P}\left( N_{n}\left( 1\right) >k\right)
z^{k}\text{ with} \\
\mathbf{P}\left( N_{n}\left( 1\right) >k\right) &=&\left[ z^{k}\right] \frac{%
1-\phi _{n}\left( z\right) }{1-z}=\frac{1}{a_{n}+b_{n}}\left( \frac{b_{n}}{%
a_{n}+b_{n}}\right) ^{k} \\
\mathbf{P}\left( N_{n}\left( 1\right) =k\right) &=&\frac{a_{n}}{\left(
a_{n}+b_{n}\right) ^{2}}\left( \frac{b_{n}}{a_{n}+b_{n}}\right) ^{k-1}\text{%
, }k\geq 1;\text{ }\mathbf{P}\left( N_{n}\left( 1\right) =0\right) =1-\frac{1%
}{a_{n}+b_{n}}.
\end{eqnarray*}

- When $n$ is large and in the subcritical case $a>1$ ($\mu <1$), 
\begin{eqnarray*}
\mathbf{P}\left( N_{n}\left( 1\right) =0\right) &\sim &1-a^{-n}\frac{a-1}{%
a+b-1}\text{ and } \\
\mathbf{P}\left( N_{n}\left( 1\right) =k\right) &\sim &a^{-n}\left( \frac{a-1%
}{a+b-1}\right) ^{2}\left( \frac{b}{a+b-1}\right) ^{k-1}\text{ if }k\geq 1,
\end{eqnarray*}
with $\mathbf{P}\left( N_{n}\left( 1\right) =k\mid N_{n}\left( 1\right)
>0\right) \sim \frac{a-1}{a+b-1}\left( \frac{b}{a+b-1}\right) ^{k-1}.$ Note
the decay rate $b/\left( a+b-1\right) =\overline{\pi }/\pi _{0}<1$.

- In the supercritical case $a<1$ ($\mu >1$), with $\frac{a+b-1}{1-a}=\rho
_{e},$%
\begin{eqnarray*}
\mathbf{P}\left( N_{n}\left( 1\right) =0\right) &\sim &\frac{a+b-1}{1-a}%
\left( 1-\frac{1-a}{b}a^{n}\right) \text{ and} \\
\mathbf{P}\left( N_{n}\left( 1\right) =k\right) &\sim &\left( \frac{1-a}{b}%
\right) ^{2}a^{n}\left( 1-\frac{1-a}{b}a^{n}\right) ^{k}\text{ if }k\geq 1,
\end{eqnarray*}
with $\mathbf{P}\left( N_{n}\left( 1\right) =k\mid N_{n}\left( 1\right)
>0\right) \sim \frac{1-a}{b}a^{n}\left( 1-\frac{1-a}{b}a^{n}\right) ^{k}.$
In the range $k=x/a^{n}$, we recover the Yaglom large-$n$ limiting Exp$%
\left( \frac{1-a}{b}\right) $ density for $\mu ^{-n}N_{n}\left( 1\right)
\mid N_{n}\left( 1\right) >0.$\newline

Observing finally 
\begin{equation*}
\phi _{n}\left( z\right) =\frac{\alpha _{n}}{\gamma _{n}}-\frac{1}{\gamma
_{n}}\frac{\alpha _{n}\delta _{n}-\beta _{n}\gamma _{n}}{\gamma _{n}z+\delta
_{n}},
\end{equation*}
this is also ($\gamma _{n}<0$) 
\begin{equation}
\mathbf{P}\left( N_{n}\left( 1\right) =0\right) =\beta _{n}/\delta _{n}\text{%
; }\mathbf{P}\left( N_{n}\left( 1\right) =k\right) =-\frac{\alpha _{n}\delta
_{n}-\beta _{n}\gamma _{n}}{\gamma _{n}\delta _{n}}\left( -\frac{\gamma _{n}%
}{\delta _{n}}\right) ^{k};\text{ }k\geq 1,  \label{coeffs}
\end{equation}
in terms of $\left( \alpha _{n},\beta _{n},\gamma _{n},\delta _{n}\right) $
resulting from the diagonalization of $A.$

\subsection{Frequency spectrum (linear-fractional model)}

Stricto sensu, the BGW process $\left\{ N_{n}\left( 1\right) \right\} $\ has
no non-trivial invariant measure, translating that BGW processes are very
unstable. Letting $n\rightarrow \infty $\ in (\ref{1.1}), we get $\phi
_{\infty }\left( z\right) \equiv \rho \wedge 1$, $z\in \left[ 0,1\right) $, $%
\phi _{\infty }\left( 1\right) =1,$\ as the only p.g.f. solution to $\phi
_{\infty }\left( z\right) =\phi _{\infty }\left( \phi \left( z\right)
\right) $. Let $\overline{P}$\ be the matrix obtained after deleting the
first row and column of $P$. Looking for a positive solution to $\mathbf{%
\varphi }^{\prime }=\mathbf{\varphi }^{\prime }\overline{P}$\ with $\mathbf{%
\varphi }:=\left( \varphi _{1},\varphi _{2},...\right) ^{\prime }$\ the
column vector of the asymptotic occupation states', $\varphi \left( z\right)
:=\sum_{i\geq 1}z^{i}\varphi _{i}$\ must solve $\varphi \left( \phi \left(
0\right) \right) =1$\ and 
\begin{equation*}
\varphi \left( z\right) =\varphi \left( \phi \left( z\right) \right) -1\text{%
, }z\in \left[ 0,\rho \wedge 1\right) ,
\end{equation*}
the Abel's functional equation.

If $\phi \left( z\right) $\ is the LF branching mechanism, with $\varphi
\left( \rho \wedge 1\right) =\infty $\ if $\rho \neq 1$, we get ($\mu =%
\overline{\pi }_{0}/\pi ,$\ $\rho =\pi _{0}/\overline{\pi }$) 
\begin{equation*}
\bullet \text{ }\varphi \left( z\right) =\left\{ 
\begin{array}{c}
1+\frac{1}{\log \left( 1/\mu \right) }\log \left( \frac{\rho -z}{\rho -\pi
_{0}}\frac{1-\pi _{0}}{1-z}\right) \text{ if }\mu <1 \\ 
\frac{\pi }{\overline{\pi }}\frac{z}{1-z}\text{ if }\mu =1 \\ 
1-\frac{1}{\log \mu }\log \left( \frac{\rho -z}{\rho -\pi _{0}}\frac{1-\pi
_{0}}{1-z}\right) \text{ if }\mu >1,
\end{array}
\right.
\end{equation*}
corresponding, with $i\geq 1$, to [see Harris (1963), p. 28] 
\begin{equation*}
\bullet \text{ }\varphi _{i}=\left\{ 
\begin{array}{c}
\frac{1}{\log \left( 1/\mu \right) }\frac{1-\rho ^{-i}}{i}\text{ if }\mu <1%
\text{ (}\rho >1\text{)} \\ 
\frac{\pi }{\overline{\pi }}\text{ if }\mu =1\text{ (}\rho =1\text{)} \\ 
\frac{1}{\log \mu }\frac{\rho ^{-i}-1}{i}\text{ if }\mu >1\text{ (}\rho
=\rho _{e}<1).
\end{array}
\right.
\end{equation*}
[Harris, (1963)] interprets $\left( \varphi _{0}=\infty ,\varphi _{i},i\geq
1\right) $ as the stationary measure of the BGW process, so long as $0\cdot
\infty =0$ is forced. He also mentions that this stationary measure is
unique (up to a multiplicative constant) when $\mu =1$; this is not the case
in general when $\mu >1,$ see [Kingman, (1965)]. Putting aside the question
of unicity, in all cases, there exists a left eigenvector $\mathbf{\varphi }$
of the substochastic matrix $\overline{P}$, associated to the eigenvalue $1$%
, obeying $\sum_{i\geq 1}\varphi _{i}=\infty $, always. The interpretation
of the $\varphi _{i}$'s as a `stationary distribution' of the LF BGW process
remains obscure to us.

\subsection{The joint law of the sterile and prolific individuals at
generation $n$}

Suppose there is a single founder. Let $\left( N_{n}^{0}\left( 1\right)
,N_{n}^{1}\left( 1\right) \right) $\ be the number of (sterile, prolific)
individuals, descending from the founder and present at generation $n$, so
with $N_{n}^{0}\left( 1\right) +N_{n}^{1}\left( 1\right) =N_{n}\left(
1\right) $. The sterile individuals are the leaves of the tree. By prolific
individuals, we mean those individuals having at least one descendant. Let 
\begin{equation}
\phi _{n}\left( z,z_{0},z_{1}\right) =\mathbf{E}\left( z^{N_{n}\left(
1\right) }z_{0}^{N_{n}^{0}\left( 1\right) }z_{1}^{N_{n}^{1}\left( 1\right)
}\right) .  \label{f9}
\end{equation}
Clearly, $\phi _{n}\left( z,z_{0},z_{1}\right) =\phi ^{\circ n}\left( \pi
_{0}zz_{0}+\overline{\pi }_{0}zz_{1}\right) $\ with $\phi ^{\circ 0}\left(
z\right) =z,$\ leading to (as from a recursion from the root) 
\begin{equation}
\bullet \text{ }\phi _{n+1}\left( z,z_{0},z_{1}\right) =\phi \left( \phi
_{n}\left( z,z_{0},z_{1}\right) \right) ;\text{ }\phi _{0}\left(
z,z_{0},z_{1}\right) =\pi _{0}zz_{0}+\overline{\pi }_{0}zz_{1}.  \label{f10}
\end{equation}
So, with $n\geq 1,$%
\begin{equation}
\phi _{n}\left( z,z_{0},z_{1}\right) =\frac{\alpha _{n}z+\beta _{n}}{\gamma
_{n}z+\delta _{n}}\mid _{z=\pi _{0}zz_{0}+\overline{\pi }_{0}zz_{1}}=\phi
_{n}\left( \pi _{0}zz_{0}+\overline{\pi }_{0}zz_{1}\right) ,  \label{f11}
\end{equation}
and 
\begin{equation}
\mathbf{E}\left( z_{0}^{N_{n}^{0}\left( 1\right) }z_{1}^{N_{n}^{1}\left(
1\right) }\mid N_{n}\left( 1\right) =k\right) =\frac{\left[ z^{k}\right]
\phi _{n}\left( z,z_{0},z_{1}\right) }{\left[ z^{k}\right] \phi _{n}\left(
z,1,1\right) }.  \label{f12}
\end{equation}
Also, 
\begin{eqnarray*}
\mathbf{E}\left( z_{0}^{N_{n}^{0}\left( 1\right) /N_{n}\left( 1\right)
}\right) &=&\phi _{n}\left( z_{0}^{-1},z_{0},1\right) =\phi ^{{{}^{\circ }}%
n}\left( \pi _{0}+\overline{\pi }_{0}/z_{0}\right) \\
&=&1-\frac{1}{b_{n}+\frac{a_{n}}{\overline{\pi }_{0}}\left( 1-1/z_{0}\right)
^{-1}},
\end{eqnarray*}
giving the law of the ratio $N_{n}^{0}\left( 1\right) /N_{n}\left( 1\right)
. $ Furthermore, 
\begin{equation*}
\mathbf{P}\left( N_{n}^{0}\left( 1\right) =0\right) =\phi _{n}\left(
1,0,1\right) =\phi ^{\circ n}\left( \overline{\pi }_{0}\right)
\end{equation*}
is the probability that generation $n$ shows no leaves. Consistently, 
\begin{equation*}
\mathbf{P}\left( N_{n}^{1}\left( 1\right) =0\right) =\phi _{n}\left(
1,1,0\right) =\phi ^{\circ n}\left( \pi _{0}\right) =\phi ^{\circ n+1}\left(
\pi _{0}\right) =\mathbf{P}\left( N_{n+1}\left( 1\right) =0\right)
\end{equation*}
is the probability that at generation $n$ there are no prolific individuals.
We refer to Section $3.7$ for some conclusions which can be drawn.

\section{Sterile vs prolific: joint laws of the past and the present for BGW
processes}

We start with recalling the law of the total progeny of a BGW process, which
is relative to its past. We then investigate the joint laws of
sterile/prolific individuals, relative to the past and the present. The
special case of the LF model is developed.

\subsection{Total progeny: the past}

Let $\overline{N}_{n}\left( i\right) $ be the cumulated number of
individuals (nodes) in the BGW tree up to generation $n$, starting from $i$
founders. From the recursion from the preceding step 
\begin{equation*}
\overline{N}_{n+1}\left( 1\right) \overset{d}{=}\overline{N}_{n}\left(
1\right) +\sum_{i=1}^{N_{n}\left( 1\right) }M_{i}.
\end{equation*}
From the recursion from the root, we have 
\begin{equation*}
\overline{N}_{n+1}\left( 1\right) \overset{d}{=}1+\sum_{m=1}^{M}\overline{N}%
_{n}^{\left( m\right) }\left( 1\right) .
\end{equation*}
With $\Phi _{n}\left( z\right) =\mathbf{E}\left( z^{\overline{N}_{n}\left(
1\right) }\right) ,$ $\Phi _{n}\left( 0\right) =0$ and $\Phi _{0}\left(
z\right) =z$, therefore [see Harris, (1963)], 
\begin{equation}
\Phi _{n+1}\left( z\right) =z\phi \left( \Phi _{n}\left( z\right) \right)
\label{f13}
\end{equation}
and 
\begin{equation*}
\mathbf{E}\left( z^{\overline{N}_{n}\left( i\right) }\right) =\Phi
_{n}\left( z\right) ^{i}.
\end{equation*}
In the (sub)-critical cases, the size of the BGW tree is finite and $\Phi
\left( z\right) =z\phi \left( \Phi \left( z\right) \right) $ is the
functional equation solving $\Phi \left( z\right) =\mathbf{E}\left( z^{%
\overline{N}_{\infty }\left( 1\right) }\right) $ with $\Phi \left( 1\right)
=1$ and $\Phi \left( 0\right) =0.$ [See the Section $3.5$ for additional
information]. In the supercritical case, the size of the BGW tree is finite
only with probability $\rho _{e}$ and $\Phi \left( z\right) =z\phi \left(
\Phi \left( z\right) \right) $ is the functional equation solving $\Phi
\left( z\right) =\mathbf{E}\left( z^{\overline{N}_{\infty }\left( 1\right)
}\right) $ with $\Phi \left( 1\right) =\rho _{e}$ and $\Phi \left( 0\right)
=0.$

With $\Phi _{n}\left( z,\overline{z}\right) =\mathbf{E}\left( z^{N_{n}\left(
1\right) }\overline{z}^{\overline{N}_{n}\left( 1\right) }\right) $ and $\Phi
_{0}\left( z,\overline{z}\right) =z\overline{z},$ the joint p.g.f. of $%
\left( N_{n}\left( 1\right) ,\overline{N}_{n}\left( 1\right) \right) $ 
\begin{eqnarray*}
\Phi _{n+1}\left( z,\overline{z}\right) &=&\Phi _{n}\left( \phi \left( z%
\overline{z}\right) ,\overline{z}\right) \\
\Phi _{n+1}\left( z,\overline{z}\right) &=&\overline{z}\phi \left( \Phi
_{n}\left( z,\overline{z}\right) \right) ,
\end{eqnarray*}
where the first recursion is from the preceding step, while the second is
from the root, [see Pakes, (1971)]. The process $\left( N_{n}\left( 1\right)
,\overline{N}_{n}\left( 1\right) \right) $ is a bivariate Markov chain whose
marginals are Markovian.

Consider now the disjoint set of nodes $\left( N_{n}\left( 1\right) ,%
\overline{N}_{n-1}\left( 1\right) \right) $, rather than looking at $%
\overline{N}_{n}\left( 1\right) $.

$\bullet $ Defining $\Psi _{n}\left( z,\overline{z}\right) =\mathbf{E}\left(
z^{N_{n}\left( 1\right) }\overline{z}^{\overline{N}_{n-1}\left( 1\right)
}\right) =\Phi _{n}\left( z/\overline{z},\overline{z}\right) $, with $\Psi
_{0}\left( z,\overline{z}\right) =z,$%
\begin{eqnarray*}
\Psi _{n+1}\left( z,\overline{z}\right) &=&\Psi _{n}\left( \overline{z}\phi
\left( z\right) ,\overline{z}\right) \\
\Psi _{n+1}\left( z,\overline{z}\right) &=&\overline{z}\phi \left( \Psi
_{n}\left( z,\overline{z}\right) \right) .
\end{eqnarray*}
Defining the `marked' p.g.f. $\phi _{\overline{z}}\left( z\right) :=%
\overline{z}\phi \left( z\right) $, 
\begin{equation*}
\Psi _{n}\left( z,\overline{z}\right) =\phi _{\overline{z}}^{\circ n}\left(
z\right) ,
\end{equation*}
where the iteration is on $z$. Note $\phi _{\overline{z}}\left( 1\right) =%
\overline{z}<1$ is the p.g.f. of a r.v. variable assigning mass $1-\overline{%
z}$ to $\infty $. We have 
\begin{equation}
\bullet \text{ }\frac{\left[ \overline{z}^{k}\right] \Psi _{n}\left( z,%
\overline{z}\right) }{\left[ \overline{z}^{k}\right] \Psi _{n}\left( 1,%
\overline{z}\right) }=\mathbf{E}\left( z^{N_{n}\left( 1\right) }\mid 
\overline{N}_{n-1}\left( 1\right) =k\right) ,  \label{f14}
\end{equation}

\begin{equation*}
\bullet \text{ }\frac{\Psi _{n}\left( z,1\right) -\Psi _{n}\left( z,0\right) 
}{1-\Psi _{n}\left( 1,0\right) }=\mathbf{E}\left( z^{N_{n}\left( 1\right)
}\mid \overline{N}_{n-1}\left( 1\right) >0\right) ,
\end{equation*}
where 
\begin{eqnarray*}
\Psi _{n}\left( z,1\right) &=&\phi ^{\circ n}\left( z\right) ,\Psi
_{n}\left( z,0\right) =\mathbf{E}\left( z^{N_{n}\left( 1\right) }1_{%
\overline{N}_{n-1}\left( 1\right) =0}\right) \text{ and} \\
\Psi _{n}\left( 1,0\right) &=&\mathbf{P}\left( \overline{N}_{n-1}\left(
1\right) =0\right) .
\end{eqnarray*}

When $\phi \left( z\right) $ is a LF branching mechanism encoded by $A$,

\begin{equation*}
A=\left[ 
\begin{array}{ll}
\pi -\pi _{0} & \pi _{0} \\ 
-\overline{\pi } & 1
\end{array}
\right] \rightarrow A_{\overline{z}}=\left[ 
\begin{array}{ll}
\overline{z}\pi -\overline{z}\pi _{0} & \overline{z}\pi _{0} \\ 
-\overline{\pi } & 1
\end{array}
\right] .
\end{equation*}
$A_{\overline{z}}^{n}$ yields in principle the expression of $\Psi
_{n}\left( z,\overline{z}\right) =\phi _{\overline{z}}^{\circ n}\left(
z\right) .$

In the supercritical case, we have $\Phi \left( 1\right) =\mathbf{P}\left( 
\overline{N}\left( 1\right) <\infty \right) =\rho _{e}$, the extinction
probability of $N_{n}\left( 1\right) $. It obeys 
\begin{equation}
\rho _{e}=\phi \left( \rho _{e}\right) .  \label{e3}
\end{equation}
with $\rho _{e}=1$ if $\mu \leq 1$, $\overline{\rho }_{e}=1-\rho _{e}>0$ if $%
\mu >1.$ If $\mu <1$, $m=\mathbf{E}\overline{N}\left( 1\right) =1/\left(
1-\mu \right) <\infty $, otherwise if $\mu \geq 1$, $m=\infty $. In the
critical case when $\mu =1$, $\overline{N}\left( 1\right) <\infty $ with
probability $1$ but $m=\mathbf{E}\overline{N}\left( 1\right) =\infty $ as a
result of $\overline{N}\left( 1\right) $ displaying heavy tails. In the
supercritical case when $\mu >1$, $m=\infty $ because with some positive
probability $\overline{\rho }_{e}$, the tree is a giant tree with infinitely
many nodes or branches (one more node than branches in a tree corresponding
to the root).

Whenever one deals with a supercritical situation with $\rho _{e}=\Phi
\left( 1\right) <1$, defining the p.g.f. of $\overline{N}_{\infty }\left(
1\right) =\overline{N}_{\infty }\left( 1\right) \mid \overline{N}_{\infty
}\left( 1\right) <\infty $ to be 
\begin{equation*}
\widetilde{\Phi }\left( z\right) =\frac{\Phi \left( z\right) -\Phi \left(
1\right) }{1-\Phi \left( 1\right) },
\end{equation*}
we have 
\begin{equation*}
\widetilde{\Phi }\left( z\right) =z\widetilde{\phi }_{\infty }\left( 
\widetilde{\Phi }\left( z\right) \right) \text{ and }\widetilde{\phi }%
_{\infty }\left( z\right) =\frac{\phi \left( z\right) -\rho _{e}}{1-\rho _{e}%
},
\end{equation*}
where $\overline{\phi }\left( z\right) $ is the modified subcritical
branching mechanism with mean $\widetilde{\phi }_{\infty }^{\prime }\left(
1\right) =\phi ^{\prime }\left( \rho _{e}\right) <1$. Conditioning a
supercritical tree on being finite is amenable to a subcritical tree problem
so with extinction probability $1$.

But this requires the computation of $\rho _{e}$ which can be quite involved
in general (although explicit in the LF case).

Indeed however (with $\left[ z^{k}\right] f\left( z\right) $ denoting the
coefficient in front of $z^{k}$ in the power-series expansion of $f\left(
z\right) $ at $0$), by Lagrange inversion formula 
\begin{equation}
\overline{\pi }_{k}=\left[ z^{k}\right] \Phi \left( z\right) =\mathbf{P}%
\left( \overline{N}\left( 1\right) =k\right) =\frac{1}{k}\left[
z^{k-1}\right] \phi \left( z\right) ^{k},  \label{e4}
\end{equation}
so that 
\begin{equation*}
\rho _{e}=\mathbf{P}\left( \overline{N}\left( 1\right) <\infty \right)
=\sum_{k\geq 1}\mathbf{P}\left( \overline{N}\left( 1\right) =k\right)
=\sum_{m\geq 0}\frac{1}{m+1}\left[ z^{m}\right] \phi \left( z\right) ^{m+1},
\end{equation*}
is the power series expansion of the extinction probability $\rho _{e}$ in
the supercritical case. There is an estimate of $\rho _{e}$ when the BGW
process is nearly supercritical ($\mu $ slightly above $1$). Let $\overline{%
\rho }_{e}=1-\rho _{e}$ be the survival probability and $f\left( z\right)
=\phi \left( z\right) -z$, with 
\begin{equation*}
f\left( 1\right) =0,f^{\prime }\left( 1\right) =\mu -1\text{ and }%
f^{^{\prime \prime }}\left( 1\right) =\mathbf{E}\left( M\left( M-1\right)
\right) =\sigma ^{2}+\mu ^{2}-\mu \underset{\mu \sim 1^{+}}{\sim }\sigma
_{c}^{2},
\end{equation*}
where $\sigma _{c}^{2}$ is the variance of $M$ at criticality. We have 
\begin{equation*}
\rho _{e}=\phi \left( \rho _{e}\right) \Leftrightarrow f\left( 1-\overline{%
\rho }_{e}\right) =0.
\end{equation*}
As a result of 
\begin{equation*}
f\left( 1-x\right) \sim f\left( 1\right) -xf^{\prime }\left( 1\right) +\frac{%
1}{2}x^{2}f^{^{\prime \prime }}\left( 1\right) ,
\end{equation*}
we get the small survival probability estimate $\overline{\rho }_{e}\sim
2\left( \mu -1\right) /\sigma _{c}^{2}$ when the BGW process is nearly
supercritical. As a function of $\mu -1$, $\overline{\rho }_{e}$ is always
continuous at $0$ ($\overline{\rho }_{e}=0$ if $\mu -1\leq 0$), but with a
discontinuous slope at $\left( \mu -1\right) _{+}$, close to $2/\sigma
_{c}^{2}<\infty $. As $\mu \rightarrow \infty $ clearly $\overline{\rho }%
_{e}\rightarrow 1.$

A full power-series expansion of $\overline{\rho }_{e}$\ in terms of $\mu
-1>0$\ can also be obtained as follows: define $\overline{\phi }\left(
z\right) $\ by $\phi \left( z\right) =1+\mu \left( z-1\right) +\overline{%
\phi }\left( 1-z\right) $, so with $\overline{\phi }\left( 0\right) =0.$\
The equation $\rho _{e}=\phi \left( \rho _{e}\right) $\ becomes 
\begin{equation*}
\frac{\overline{\phi }\left( \overline{\rho }_{e}\right) }{\overline{\rho }%
_{e}}=\mu -1.
\end{equation*}
Lagrange inversion formula gives 
\begin{equation}
\bullet \text{ }\overline{\rho }_{e}=\sum_{k\geq 1}\rho _{k}\left( \mu
-1\right) ^{k}\text{, with}  \label{e4a}
\end{equation}
\begin{equation*}
\rho _{k}=\frac{1}{k}\left[ x^{k-1}\right] \left( \frac{\overline{\phi }%
\left( x\right) }{x^{2}}\right) ^{-k}.
\end{equation*}
Note $\rho _{1}=2/\phi ^{\prime \prime }\left( 1\right) $ with $\phi
^{\prime \prime }\left( 1\right) \sim \sigma _{c}^{2}$ when $\mu $ is
slightly above $1$. To the first order in $\mu -1$, we recover $\overline{%
\rho }_{e}\sim 2\left( \mu -1\right) /\sigma _{c}^{2}$. The second-order
coefficient is found to be $\rho _{2}=4/3\cdot \phi ^{\prime \prime \prime
}\left( 1\right) /\phi ^{\prime \prime }\left( 1\right) ^{3}.$ Let us check
these formulas on an explicit example. \newline

\textbf{Example:} If $\phi \left( z\right) =1/\left( 1+\mu \left( 1-z\right)
\right) $, with $\mu >1$ (the shifted geometric case), the fixed point $\rho
_{e}=1/\mu $ is explicitly found. Here $\overline{\phi }\left( x\right)
/x^{2}=\mu ^{2}/\left( 1+\mu x\right) $ with $\rho _{k}=\mu ^{-\left(
k+1\right) }.$ Thus, consistently, $\overline{\rho }_{e}=\sum_{k\geq 1}\rho
_{k}\left( \mu -1\right) ^{k}=1-1/\mu $ and, owing to $\phi ^{\prime \prime
}\left( 1\right) =2\mu ^{2}\sim \sigma _{c}^{2}=2$ as $\mu \rightarrow 1^{+}$
and $\phi ^{\prime \prime \prime }\left( 1\right) =6\mu ^{3}\sim 6$, $\rho
_{2}=\mu ^{-3}=4/3\cdot \phi ^{\prime \prime \prime }\left( 1\right) /\phi
^{\prime \prime }\left( 1\right) ^{3}\sim 1.$

For the general LF branching mechanism, $\sigma _{c}^{2}=2\pi _{0}/\pi $.
The latter example is the particular Geo$_{0}\left( \pi \right) $ case with $%
\pi _{0}=\pi $. $\triangleright $

\subsection{The joint laws of the current and cumulated sterile and prolific
individuals at and up to generation $n$}

Let $\left( \overline{N}_{n}^{0}\left( 1\right) ,\overline{N}_{n}^{1}\left(
1\right) \right) $\ be the cumulated number of (sterile, prolific)
individuals, descending from the founder up to generation $n$, so with $%
\overline{N}_{n}^{0}\left( 1\right) +\overline{N}_{n}^{1}\left( 1\right) =%
\overline{N}_{n}\left( 1\right) $. Let 
\begin{equation}
\Phi _{n}\left( z,z_{0},z_{1},\overline{z},\overline{z}_{0},\overline{z}%
_{1}\right) =\mathbf{E}\left( z^{N_{n}\left( 1\right)
}z_{0}^{N_{n}^{0}\left( 1\right) }z_{1}^{N_{n}^{1}\left( 1\right) }\overline{%
z}^{\overline{N}_{n}\left( 1\right) }\overline{z}_{0}^{\overline{N}%
_{n}^{0}\left( 1\right) }\overline{z}_{1}^{\overline{N}_{n}^{1}\left(
1\right) }\right) .  \label{f15}
\end{equation}
Clearly, with\textit{\ }$\Phi _{0}\left( z,z_{0},z_{1},\overline{z},%
\overline{z}_{0},\overline{z}_{1}\right) =z\overline{z}\left( \pi _{0}z_{0}%
\overline{z}_{0}+\overline{\pi }_{0}z_{1}\overline{z}_{1}\right) ,$\textit{\ 
}a recursion from the root yields 
\begin{equation}
\bullet \text{ }\Phi _{n+1}\left( z,z_{0},z_{1},\overline{z},\overline{z}%
_{0},\overline{z}_{1}\right) =\pi _{0}\overline{z}\left( \overline{z}_{0}-%
\overline{z}_{1}\right) +\overline{z}\overline{z}_{1}\phi \left( \Phi
_{n}\left( z,z_{0},z_{1},\overline{z},\overline{z}_{0},\overline{z}%
_{1}\right) \right) .  \label{f16}
\end{equation}

Three particular cases of interest are:

$\left( \mathbf{i}\right) $ Note $\Phi _{n}\left( 1,1,1,1,\overline{z}_{0},%
\overline{z}_{1}\right) =:\Phi _{n}\left( \overline{z}_{0},\overline{z}%
_{1}\right) $ gives the joint law of $\left( \overline{N}_{n}^{0}\left(
1\right) ,\overline{N}_{n}^{1}\left( 1\right) \right) .$ It obeys 
\begin{equation*}
\Phi _{n+1}\left( \overline{z}_{0},\overline{z}_{1}\right) =\overline{\phi }%
_{\overline{z}_{0},\overline{z}_{1}}\left( \Phi _{n}\left( \overline{z}_{0},%
\overline{z}_{1}\right) \right) ;\text{ }\Phi _{0}\left( \overline{z}_{0},%
\overline{z}_{1}\right) =\pi _{0}\overline{z}_{0}+\overline{\pi }_{0}%
\overline{z}_{1},
\end{equation*}
where $\overline{\phi }_{\overline{z}_{0},\overline{z}_{1}}\left( z\right)
:=\pi _{0}\left( \overline{z}_{0}-\overline{z}_{1}\right) +\overline{z}%
_{1}\phi \left( z\right) $, resulting in:

\begin{equation*}
\Phi _{n}\left( \overline{z}_{0},\overline{z}_{1}\right) =\overline{\phi }_{%
\overline{z}_{0},\overline{z}_{1}}^{\circ n}\left( \pi _{0}\overline{z}_{0}+%
\overline{\pi }_{0}\overline{z}_{1}\right) ,
\end{equation*}
the $n^{\text{th}}$-iterate of the `marked' generating function (g.f.): $%
\overline{\phi }_{\overline{z}_{0},\overline{z}_{1}}$ evaluated at $\Phi
_{0}\left( \overline{z}_{0},\overline{z}_{1}\right) $.\newline

$\bullet $ Defining $\Psi _{n}\left( z,z_{0},z_{1},\overline{z},\overline{z}%
_{0},\overline{z}_{1}\right) =\mathbf{E}\left( z^{N_{n}\left( 1\right)
}z_{0}^{N_{n}^{0}\left( 1\right) }z_{1}^{N_{n}^{1}\left( 1\right) }\overline{%
z}^{\overline{N}_{n-1}\left( 1\right) }\overline{z}_{0}^{\overline{N}%
_{n-1}^{0}\left( 1\right) }z_{1}^{\overline{N}_{n-1}^{1}\left( 1\right)
}\right) =\Phi _{n}\left( z/\overline{z},z_{0}/\overline{z}_{0},z_{1}/%
\overline{z}_{1},\overline{z},\overline{z}_{0},\overline{z}_{1}\right) $, $%
\Psi _{n}$\ obeys the same recurrence relation (\ref{f16}) than $\Phi _{n}$,
but now with the initial condition\textit{\ }$\Psi _{0}\left( z,z_{0},z_{1},%
\overline{z},\overline{z}_{0},\overline{z}_{1}\right) =z\left( \pi _{0}z_{0}+%
\overline{\pi }_{0}z_{1}\right) .$ We shall consider two other special cases:%
\newline

$\left( \mathbf{ii}\right) $%
\begin{eqnarray*}
\Psi _{n}\left( z,1,1,1,\overline{z}_{0},1\right) &=&:\Psi _{n}\left( z,%
\overline{z}_{0}\right) =\mathbf{E}\left( z^{N_{n}\left( 1\right) }\overline{%
z}_{0}^{\overline{N}_{n-1}^{0}\left( 1\right) }\right) \text{, }\Psi
_{0}\left( z,\overline{z}_{0}\right) =z, \\
\Psi _{n}\left( z,\overline{z}_{1}\right) &=&\phi _{\overline{z}_{0}}^{\circ
n}\left( z\right) \text{ where }\phi _{\overline{z}_{0}}\left( z\right) =\pi
_{0}\left( \overline{z}_{0}-1\right) +\phi \left( z\right) ,
\end{eqnarray*}
giving the joint law of the number of individuals alive at $n$ and the
cumulated number of sterile individuals up to generation $n-1$ (in view of
the forthcoming discussion: do the cumulated number of sterile individuals
in the past exceed (or not) the current population size? see [Howard, (2012)
and Avan et al. (2015)].

\begin{equation*}
\mathbf{E}\left( \overline{z}_{0}^{\overline{N}_{n-1}^{0}\left( 1\right)
}\mid N_{n}\left( 1\right) =k\right) =\frac{\left[ z^{k}\right] \Psi
_{n}\left( z,\overline{z}_{0}\right) }{\left[ z^{k}\right] \Psi _{n}\left(
z,1\right) }
\end{equation*}
gives the law of $\overline{N}_{n-1}^{0}\left( 1\right) $ given $N_{n}\left(
1\right) =k.$\newline

$\left( \mathbf{iii}\right) $

\begin{eqnarray*}
\Psi _{n}\left( z,1,1,1,1,\overline{z}_{1}\right) &=&:\Psi _{n}\left( z,%
\overline{z}_{1}\right) =\mathbf{E}\left( z^{N_{n}\left( 1\right) }\overline{%
z}_{1}^{\overline{N}_{n-1}^{1}\left( 1\right) }\right) ,\text{ }\Psi
_{0}\left( z,\overline{z}_{1}\right) =z, \\
\Psi _{n}\left( z,\overline{z}_{1}\right) &=&\phi _{\overline{z}_{1}}^{\circ
n}\left( z\right) \text{ where }\phi _{\overline{z}_{1}}\left( z\right)
:=\pi _{0}\left( 1-\overline{z}_{1}\right) +\overline{z}_{1}\phi \left(
z\right) ,
\end{eqnarray*}
giving the joint law of the number of individuals alive at $n$ and the
cumulated number of prolific individuals up to generation $n-1.$ This
situation is developed in Section $4.1$ where it appears in a discrete
version of the Lamperti's theorem.

Both cases $\left( \mathbf{ii}\right) $ and $\left( \mathbf{iii}\right) $
have initial condition $z$ and so, as a function of $z$, are iterates of
`marked' generating functions.

\subsection{Joint law of $\left( \overline{N}_{n}^{0}\left( 1\right) ,%
\overline{N}_{n}^{1}\left( 1\right) \right) $ in the linear-fractional case}

It is case $\left( \mathbf{i}\right) $. It requires the $n^{\text{th}}$%
-iteration of $\overline{\phi }_{\overline{z}_{0},\overline{z}_{1}}$ which
is a LF g.f. (not a p.g.f. because $\overline{\phi }_{\overline{z}_{0},%
\overline{z}_{1}}\left( 1\right) =\pi _{0}\overline{z}_{0}+\overline{\pi }%
_{0}\overline{z}_{1}\neq 1$).

$\overline{\phi }_{\overline{z}_{0},\overline{z}_{1}}$ can be put under the
form $\left( \alpha z+\beta \right) /\left( \gamma z+\delta \right) $ with:

\begin{eqnarray*}
\alpha &=&\overline{\pi }_{0}\pi \overline{z}_{1}-\overline{\pi }_{0}%
\overline{\pi }\overline{z}_{0};\text{ }\beta =\overline{\pi }_{0}\overline{z%
}_{0} \\
\gamma &=&-\pi ;\text{ }\delta =1
\end{eqnarray*}
and so: 
\begin{equation*}
\mathbf{E}\left( \overline{z}_{0}^{\overline{N}_{n}^{0}\left( 1\right) }%
\overline{z}_{1}^{\overline{N}_{n}^{1}\left( 1\right) }\right) =\overline{%
\phi }_{\overline{z}_{0},\overline{z}_{1}}^{\circ n}\left( \pi _{0}\overline{%
z}_{0}+\overline{\pi }_{0}\overline{z}_{1}\right) =\frac{\alpha _{n}\left(
\pi _{0}\overline{z}_{0}+\overline{\pi }_{0}\overline{z}_{1}\right) +\beta
_{n}}{\gamma _{n}\left( \pi _{0}\overline{z}_{0}+\overline{\pi }_{0}%
\overline{z}_{1}\right) +\delta _{n}}.
\end{equation*}
The fixed points of the transformation $\overline{\phi }_{\overline{z}_{0},%
\overline{z}_{1}}\left( z\right) :=\pi _{0}\left( \overline{z}_{0}-\overline{%
z}_{1}\right) +\overline{z}_{1}\phi \left( z\right) $ solving $\overline{%
\phi }_{\overline{z}_{0},\overline{z}_{1}}\left( z\right) =z$ are:

\begin{equation*}
z_{\pm }:=z_{\pm }\left( \overline{z}_{0},\overline{z}_{1}\right) =\frac{%
\left( 1-\overline{\pi }_{0}\pi \overline{z}_{1}+\pi _{0}\overline{\pi }%
\overline{z}_{0}\right) \pm \sqrt{\Delta }}{2\overline{\pi }},
\end{equation*}
where $\Delta =\left( 1-\overline{\pi }_{0}\pi \overline{z}_{1}+\pi _{0}%
\overline{\pi }\overline{z}_{0}\right) ^{2}-4\pi _{0}\overline{\pi }%
\overline{z}_{0}.$ From the conjugacy property stating that, with $T\left(
z\right) =\frac{z-z_{+}}{z-z_{-}},$

\begin{equation}
\bullet \text{ }u\frac{1-\overline{\pi }z_{-}}{1-\overline{\pi }z_{+}}%
=T\circ \overline{\phi }_{\overline{z}_{0},\overline{z}_{1}}\left( \cdot
\right) \circ T^{-1}\left( u\right) ,  \label{conj}
\end{equation}
is conjugate to $\overline{\phi }_{\overline{z}_{0},\overline{z}_{1}}\left(
\cdot \right) $, we get 
\begin{equation*}
\frac{\overline{\phi }_{\overline{z}_{0},\overline{z}_{1}}\left( z\right)
-z_{+}}{\overline{\phi }_{\overline{z}_{0},\overline{z}_{1}}-z_{-}}=\frac{1-%
\overline{\pi }z_{-}}{1-\overline{\pi }z_{+}}\frac{z-z_{+}}{z-z_{-}}.
\end{equation*}
Upon iteration, we get: 
\begin{eqnarray*}
\frac{\overline{\phi }_{\overline{z}_{0},\overline{z}_{1}}^{\circ n}\left(
z\right) -z_{+}}{\overline{\phi }_{\overline{z}_{0},\overline{z}_{1}}^{\circ
n}\left( z\right) -z_{-}} &=&\frac{\overline{\phi }_{\overline{z}_{0},%
\overline{z}_{1}}\left( \overline{\phi }_{\overline{z}_{0},\overline{z}%
_{1}}^{\circ n-1}\left( z\right) \right) -z_{+}}{\overline{\phi }_{\overline{%
z}_{0},\overline{z}_{1}}\left( \overline{\phi }_{\overline{z}_{0},\overline{z%
}_{1}}^{\circ n-1}\left( z\right) \right) -z_{-}} \\
&=&\frac{1-\overline{\pi }z_{-}}{1-\overline{\pi }z_{+}}\frac{\overline{\phi 
}_{\overline{z}_{0},\overline{z}_{1}}^{\circ n-1}\left( z\right) -z_{+}}{%
\overline{\phi }_{\overline{z}_{0},\overline{z}_{1}}^{\circ n-1}\left(
z\right) -z_{-}} \\
&=&...=\left( \frac{1-\overline{\pi }z_{-}}{1-\overline{\pi }z_{+}}\right)
^{n}\frac{z-z_{+}}{z-z_{-}},
\end{eqnarray*}
and so\textit{\ } 
\begin{eqnarray*}
\bullet \text{ }\mathbf{E}\left( \overline{z}_{0}^{\overline{N}%
_{n}^{0}\left( 1\right) }\overline{z}_{1}^{\overline{N}_{n}^{1}\left(
1\right) }\right) &=&\overline{\phi }_{\overline{z}_{0},\overline{z}%
_{1}}^{\circ n}\left( \pi _{0}\overline{z}_{0}+\overline{\pi }_{0}\overline{z%
}_{1}\right) \\
&=&z_{-}\left( \overline{z}_{0},\overline{z}_{1}\right) +\frac{z_{+}\left( 
\overline{z}_{0},\overline{z}_{1}\right) -z_{-}\left( \overline{z}_{0},%
\overline{z}_{1}\right) }{1-\left( \frac{1-\overline{\pi }z_{-}\left( 
\overline{z}_{0},\overline{z}_{1}\right) }{1-\overline{\pi }z_{+}\left( 
\overline{z}_{0},\overline{z}_{1}\right) }\right) ^{n}\frac{z-z_{+}\left( 
\overline{z}_{0},\overline{z}_{1}\right) }{z-z_{-}\left( \overline{z}_{0},%
\overline{z}_{1}\right) }}\mid _{z=\pi _{0}\overline{z}_{0}+\overline{\pi }%
_{0}\overline{z}_{1}}.
\end{eqnarray*}
Note that the joint p.g.f. of $\left( \overline{N}_{n}^{0}\left( 1\right) ,%
\overline{N}_{n}\left( 1\right) \right) $ is given by 
\begin{equation*}
\Phi _{n}\left( \overline{z}_{0},\overline{z}\right) :=\mathbf{E}\left( 
\overline{z}_{0}^{\overline{N}_{n}^{0}\left( 1\right) }\overline{z}^{%
\overline{N}_{n}\left( 1\right) }\right) =\mathbf{E}\left( \overline{z}_{0}^{%
\overline{N}_{n}^{0}\left( 1\right) }\left( \overline{z}_{0}\overline{z}%
\right) ^{\overline{N}_{n}^{1}\left( 1\right) }\right) .
\end{equation*}
In the subcritical case, it obeys $\Phi _{n+1}\left( \overline{z}_{0},%
\overline{z}\right) =\overline{z}\phi _{\overline{z}_{0}}\left( \Phi
_{n}\left( \overline{z}_{0},\overline{z}\right) \right) $, $\Phi _{0}\left( 
\overline{z}_{0},\overline{z}\right) =\overline{z}$, with $\phi _{\overline{z%
}_{0}}\left( z\right) =\pi _{0}\left( \overline{z}_{0}-1\right) +\phi \left(
z\right) $, with $\overline{z}_{0}$ viewed as a parameter. We refer to
Section $3.7$ for asymptotic results ($n\rightarrow \infty $) in the
subcritical case making use of this recurrence.

\subsection{Joint law of the number of individuals alive at $n$ and the
cumulated number of sterile individuals up to generation $n-1$}

In the case $\left( \mathbf{ii}\right) $, with $\overline{z}_{0}\in \left[
0,1\right] ,$ $\phi _{\overline{z}_{0}}\left( z\right) =\pi _{0}\left( 
\overline{z}_{0}-1\right) +\phi \left( z\right) $, with $\phi _{\overline{z}%
_{0}}\left( 0\right) =\pi _{0}\overline{z}_{0}<1$ and $1>\phi _{\overline{z}%
_{0}}\left( 1\right) =\pi _{0}\left( \overline{z}_{0}-1\right) +1>\phi _{%
\overline{z}_{0}}\left( 0\right) $. Here, $\Psi _{n}\left( z,\overline{z}%
_{0}\right) =\phi _{\overline{z}_{0}}^{\circ n}\left( z\right) ,$ with 
\begin{equation*}
A=\left[ 
\begin{array}{ll}
\pi -\pi _{0} & \pi _{0} \\ 
-\overline{\pi } & 1
\end{array}
\right] \rightarrow A_{\overline{z}_{0}}=\left[ 
\begin{array}{ll}
\overline{\pi }_{0}\pi -\pi _{0}\overline{\pi }\overline{z}_{0} & \pi _{0}%
\overline{z}_{0} \\ 
-\overline{\pi } & 1
\end{array}
\right] \text{.}
\end{equation*}
The matrix $A_{\overline{z}_{0}}^{n}$ could be computed to compute $\Psi
_{n}\left( z,\overline{z}_{0}\right) $ but we adopt a different point of
view, based on (\ref{conj}). The search for fixed points: $\phi _{\overline{z%
}_{0}}\left( z\right) =z$ yields: 
\begin{equation*}
z_{\pm }\left( \overline{z}_{0}\right) =\frac{1+\pi _{0}\overline{\pi }%
\overline{z}_{0}-\overline{\pi }_{0}\pi \pm \sqrt{\Delta }}{2\overline{\pi }}%
,
\end{equation*}
where 
\begin{equation*}
\Delta =\left( 1+\pi _{0}\overline{\pi }\overline{z}_{0}-\overline{\pi }%
_{0}\pi \right) ^{2}-4\pi _{0}\overline{\pi }\overline{z}_{0}>0
\end{equation*}
and $0<z_{-}\left( \overline{z}_{0}\right) \leq 1<z_{+}\left( \overline{z}%
_{0}\right) <z_{*}.$ In particular, $z_{-}\left( 1\right) =\pi _{0}/%
\overline{\pi }$, $z_{+}\left( 1\right) =1$ if $\pi _{0}<\overline{\pi }$
(subcritical case) or $z_{-}\left( 1\right) =1$, $z_{+}\left( 1\right) =\pi
_{0}/\overline{\pi }$ if $\pi _{0}>\overline{\pi }$ (supercritical case). If 
$\pi _{0}=\overline{\pi }$ (critical case), $z_{\pm }\left( 1\right) =1$.
From (\ref{conj}), for all $z\in \left[ 0,1\right] $ therefore 
\begin{eqnarray*}
\bullet \text{ }\Psi _{n}\left( z,\overline{z}_{0}\right) &=&\mathbf{E}%
\left( z^{N_{n}\left( 1\right) }\overline{z}_{0}^{\overline{N}%
_{n-1}^{0}\left( 1\right) }\right) =z_{-}\left( \overline{z}_{0}\right) +%
\frac{z_{+}\left( \overline{z}_{0}\right) -z_{-}\left( \overline{z}%
_{0}\right) }{1-\left( \frac{1-\overline{\pi }z_{-}\left( \overline{z}%
_{0}\right) }{1-\overline{\pi }z_{+}\left( \overline{z}_{0}\right) }\right)
^{n}\frac{z-z_{+}\left( \overline{z}_{0}\right) }{z-z_{-}\left( \overline{z}%
_{0}\right) }} \\
&\rightarrow &z_{-}\left( \overline{z}_{0}\right) \text{ geometrically fast
as }n\rightarrow \infty .
\end{eqnarray*}
and 
\begin{equation}
\mathbf{E}_{\overline{z}_{0}}\left( z^{N_{n}\left( 1\right) }\right) :=\frac{%
\Psi _{n}\left( z,\overline{z}_{0}\right) }{\Psi _{n}\left( 1,\overline{z}%
_{0}\right) }\rightarrow 1\text{ geometrically fast as }n\rightarrow \infty .
\label{f17}
\end{equation}
With $z_{-}\left( 0\right) =0$ and $z_{+}\left( 0\right) =\left( 1-\overline{%
\pi }_{0}\pi \right) /\overline{\pi }>1,$

\begin{equation*}
\mathbf{P}\left( \overline{N}_{n}^{0}\left( 1\right) =0\right) =\Psi
_{n+1}\left( 1,0\right) =\frac{z_{+}\left( 0\right) }{1+\left( \overline{\pi 
}_{0}\pi \right) ^{-\left( n+1\right) }\left( z_{+}\left( 0\right) -1\right) 
},
\end{equation*}
going geometrically fast to $0$.

Of interest is the conditional p.g.f. of the cumulated number of sterile
individuals $\overline{N}_{n-1}^{0}\left( 1\right) $ given the current
population size $N_{n}\left( 1\right) :$%
\begin{equation}
\bullet \text{ }\mathbf{E}\left( \overline{z}_{0}^{\overline{N}%
_{n-1}^{0}\left( 1\right) }\mid N_{n}\left( 1\right) =k\right) =\frac{\left[
z^{k}\right] \Psi _{n}\left( z,\overline{z}_{0}\right) }{\left[ z^{k}\right]
\Psi _{n}\left( z,1\right) }.  \label{f18}
\end{equation}
It can be found explicitly because $\Psi _{n}\left( z,\overline{z}%
_{0}\right) $ is under the form of a LF model in $z$. With $a\left( 
\overline{z}_{0}\right) :=\left( \frac{1-\overline{\pi }z_{-}\left( 
\overline{z}_{0}\right) }{1-\overline{\pi }z_{+}\left( \overline{z}%
_{0}\right) }\right) >1$, $\Psi _{n}\left( z,\overline{z}_{0}\right) $ can
indeed be put under the form $\left( \alpha _{n}z+\beta _{n}\right) /\left(
\gamma _{n}z+\delta _{n}\right) $ with 
\begin{eqnarray*}
\alpha _{n} &=&z_{+}\left( \overline{z}_{0}\right) -z_{-}\left( \overline{z}%
_{0}\right) -z_{-}\left( \overline{z}_{0}\right) \left( a\left( \overline{z}%
_{0}\right) ^{n}-1\right) \\
\beta _{n} &=&z_{-}\left( \overline{z}_{0}\right) z_{+}\left( \overline{z}%
_{0}\right) \left( a\left( \overline{z}_{0}\right) ^{n}-1\right) \\
\gamma _{n} &=&-\left( a\left( \overline{z}_{0}\right) ^{n}-1\right) \\
\delta _{n} &=&z_{+}\left( \overline{z}_{0}\right) -z_{-}\left( \overline{z}%
_{0}\right) +z_{+}\left( \overline{z}_{0}\right) \left( a\left( \overline{z}%
_{0}\right) ^{n}-1\right)
\end{eqnarray*}
\begin{equation*}
\alpha _{n}\delta _{n}-\beta _{n}\gamma _{n}=a\left( \overline{z}_{0}\right)
^{n}\left( z_{+}\left( \overline{z}_{0}\right) -z_{-}\left( \overline{z}%
_{0}\right) \right) ^{2}.
\end{equation*}
The $z^{k}-$coefficient of $\Psi _{n}\left( z,\overline{z}_{0}\right) $
[respectively $\Psi _{n}\left( z,1\right) $] are then given by (\ref{coeffs}%
) in terms of the fixed points $z_{\pm }\left( \overline{z}_{0}\right) $
[respectively $z_{\pm }\left( 1\right) $]. When $n$ is large ($n>>1$) 
\begin{equation}
\left[ z^{k}\right] \Psi _{n}\left( z,\overline{z}_{0}\right) \sim \left(
z_{+}\left( \overline{z}_{0}\right) -z_{-}\left( \overline{z}_{0}\right)
\right) ^{2}a\left( \overline{z}_{0}\right) ^{-n}z_{+}\left( \overline{z}%
_{0}\right) ^{-\left( k+1\right) }.  \label{f19}
\end{equation}
If in addition $k>>n>>1$%
\begin{equation}
\mathbf{E}\left( \overline{z}_{0}^{\overline{N}_{n-1}^{0}\left( 1\right)
}\mid N_{n}\left( 1\right) =k\right) ^{1/k}\sim \alpha \left( \overline{z}%
_{0}\right) :=\frac{z_{+}\left( 1\right) }{z_{+}\left( \overline{z}%
_{0}\right) }.  \label{f20}
\end{equation}
\textbf{Example:} (this is a very naive estimate). For a population whose
founder age is $2.10^{5}$\ years, considering the time elapsed between two
consecutive generations is about $20$\ years (this is questionable as this
time could vary with time in the past), the current number of generations
away from the founder is $n=10^{4}$. If the current population size is $%
k=8.10^{9}>>n$\ individuals, with $\rho \left( \pi _{0},\pi \right)
:=-F^{\prime }\left( 0\right) $\ where $F\left( \lambda \right) =\log \alpha
\left( e^{-\lambda }\right) ,$\ by Cram\'{e}r's theorem, [Cram\'{e}r,
(1938)], 
\begin{equation}
\frac{1}{k}\left( \overline{N}_{n-1}^{0}\left( 1\right) \mid N_{n}\left(
1\right) =k\right) \sim \rho \left( \pi _{0},\pi \right) :=\frac{%
z_{+}^{\prime }\left( 1\right) }{z_{+}\left( 1\right) }>0,\text{ a.s.,}
\label{f21}
\end{equation}
accrediting the fact that the cumulated number of sterile individuals could
be of the same order of magnitude $\rho $\ (several billions) than the
current population size, depending on the values of the independent
parameters $\left( \pi _{0},\pi \right) $\ of the LF branching mechanism.
More precisely, observing 
\begin{eqnarray*}
z_{+}^{\prime }\left( 1\right) &=&\frac{\pi _{0}\overline{\pi }_{0}}{\pi
_{0}-\overline{\pi }},\text{ }\rho \left( \pi _{0},\pi \right) =\frac{\pi
_{0}\overline{\pi }_{0}}{\pi _{0}-\overline{\pi }}\text{ if }\pi _{0}>%
\overline{\pi } \\
z_{+}^{\prime }\left( 1\right) &=&\frac{\pi _{0}\pi }{\overline{\pi }-\pi
_{0}},\text{ }\rho \left( \pi _{0},\pi \right) =\frac{\pi \overline{\pi }}{%
\overline{\pi }-\pi _{0}}\text{ if }\pi _{0}<\overline{\pi },
\end{eqnarray*}
respectively for the subcritical (supercritical) BGW process, we get 
\begin{equation*}
\bullet \text{ }\rho \left( \pi _{0},\pi \right) >1\text{ both }\left\{ 
\begin{array}{c}
\text{if }\sqrt{\overline{\pi }}>\pi _{0}>\overline{\pi } \\ 
\text{if }\sqrt{\pi _{0}}>\overline{\pi }>\pi _{0}.
\end{array}
\right.
\end{equation*}
clarifying the conditions on $\left( \pi _{0},\pi \right) $\ under which the
cumulated number of sterile individuals can exceed the current number of
prolific ones. If $\pi _{0}<\overline{\pi },$\ for $\rho \left( \pi _{0},\pi
\right) $\ to be of order say of few tens, $\pi _{0}$\ and $\overline{\pi }$%
\ both need to be quite close to one another. For example, $\pi _{0}=0.400$\
and $\overline{\pi }=0.405$\ yields $\rho \left( \pi _{0},\pi \right) =\frac{%
\pi \overline{\pi }}{\overline{\pi }-\pi _{0}}=48.$\ This corresponds to a
nearly supercritical BGW with $\mu =\overline{\pi }_{0}/\pi =1.008$. Note
though that the corresponding event $N_{n}\left( 1\right) =k$\ has an
extremely small probability to occur.

In the opposite direction, $\overline{\pi }<\pi _{0}^{2}<\pi _{0}$\ and also 
$\pi _{0}<\overline{\pi }^{2}<\overline{\pi }$\ are conditions for the
current number of prolific individuals to exceed the cumulated number of
sterile ones over the past ($\rho \left( \pi _{0},\pi \right) <1$). For
instance, $\pi _{0}=0.5$\ and $\overline{\pi }=0.2$\ yields $\rho \left( \pi
_{0},\pi \right) =\frac{\pi _{0}\overline{\pi }_{0}}{\pi _{0}-\overline{\pi }%
}=0.833<1,$\ ($\mu =0.625<1$)$.$\ And $\pi _{0}=0.4$\ and $\overline{\pi }%
=0.7$\ yields $\rho \left( \pi _{0},\pi \right) =\frac{\pi \overline{\pi }}{%
\overline{\pi }-\pi _{0}}=0.7<1,$\ ($\mu =2>1$)$.$\ For this last situation,
the event $N_{n}\left( 1\right) =k$\ has the largest (although still very
small) probability to occur. $\triangleright $

\subsection{Back to the total progeny (subcritical case)}

In this Section, the BGW process is assumed to be subcritical, so that $%
\overline{N}_{\infty }\left( 1\right) <\infty $ a.s.. For some rare specific
models for $\phi $ the limiting probabilities $\overline{\pi }_{k}=\mathbf{P}%
\left( \overline{N}_{\infty }\left( 1\right) =k\right) $ (or its p.g.f. $%
\Phi \left( \overline{z}\right) =\lim_{n\rightarrow \infty }\Phi _{n}\left( 
\overline{z}\right) $) can be explicitly computed. This is the case for the
LF $\phi $ for which $\Phi \left( \overline{z}\right) $ is the solution to
the quadratic equation $\Phi \left( \overline{z}\right) =\overline{z}\phi
\left( \Phi \left( \overline{z}\right) \right) ,$ showing an algebraic
dominant singularity of order $-1/2$ at some\emph{\ }$z_{c}>1$ obtained
while cancelling the discriminant$.$

For instance, assuming $\pi _{0}=\pi >1/2$ ($\phi \left( z\right) =\pi
/\left( 1-\overline{\pi }z\right) ,$ the Geo$_{0}\left( \pi \right) $
special case) yields the exact expression 
\begin{equation*}
\Phi \left( \overline{z}\right) =\frac{1}{2\overline{\pi }}\left( 1-\sqrt{%
1-4\pi \overline{\pi }\overline{z}}\right) ,
\end{equation*}
with an algebraic singularity of order $-1/2$ at $z_{c}=1/\phi ^{\prime
}\left( \tau \right) =1/\left( 4\pi \overline{\pi }\right) >1.$ Note $\Phi
\left( z_{c}\right) =\tau =\frac{1}{2\overline{\pi }}>1$ and $\phi \left(
\tau \right) =\tau \phi ^{\prime }\left( \tau \right) =2\pi >1$ (the
subcriticality condition)$.$ With $\left[ a\right] _{k}:=a\left( a+1\right)
...\left( a+k-1\right) $, we get (denoting $Pi$ :$=3.1415...$) from Lagrange
inversion formula, 
\begin{eqnarray*}
\overline{\pi }_{k} &=&\left[ \overline{z}^{k}\right] \Phi \left( \overline{z%
}\right) =\frac{\left( 4\pi \overline{\pi }\right) ^{k}}{4\overline{\pi }}%
\frac{\left[ 1/2\right] _{k-1}}{k!} \\
&=&\frac{\left( \pi \overline{\pi }\right) ^{k}}{\overline{\pi }}\frac{%
\left( 2k-2\right) !}{k!\left( k-1\right) !}\underset{k\rightarrow \infty }{%
\sim }\frac{1}{\sqrt{2Pi}4\overline{\pi }}k^{-3/2}z_{c}^{-k}.
\end{eqnarray*}
For a general (aperiodic and different from an affine function, see Remark
below) $\phi $ obeying: $\phi $ has convergence radius $z_{*}>1$ (possibly $%
z_{*}=\infty $) and $\pi _{0}>0$, a large $k$ estimate for $\overline{\pi }%
_{k}$ can be obtained in general. For such $\phi $'s indeed, the unique
positive real root to the equation 
\begin{equation}
\phi \left( \tau \right) -\tau \phi ^{\prime }\left( \tau \right) =0,
\label{e9}
\end{equation}
exists, with $\rho _{e}=1<\tau <z_{*}$ if $\mu <1$ (assuming the subcritical
case).

\textbf{Remark (}affinity)\textbf{:} When $\phi \left( z\right) =\overline{%
\alpha }+\alpha z$ is affine (pure death Bernoulli branching mechanism), the
number $\tau $ below is rejected at $\infty $ and the following analysis of
the corresponding $\Phi \left( \overline{z}\right) $ is invalid. This case
deserves a special treatment. $\blacksquare $\newline

The point $\left( \tau ,\phi \left( \tau \right) \right) $ is indeed the
tangency point to the curve $\phi \left( z\right) $ of a straight line
passing through the origin $\left( 0,0\right) $. Let then $z_{c}:=\tau /\phi
\left( \tau \right) =1/\phi ^{\prime }\left( \tau \right) \geq 1$. The
searched $\Phi \left( \overline{z}\right) $ solves $\psi \left( \Phi \left( 
\overline{z}\right) \right) =\overline{z},$ where $\psi \left( z\right)
=z/\phi \left( z\right) $ obeys $\psi \left( \tau \right) =z_{c},$ $\psi
^{\prime }\left( \tau \right) =0$ and $\psi ^{\prime \prime }\left( \tau
\right) =-\frac{\tau \phi ^{\prime \prime }\left( \tau \right) }{\phi \left(
\tau \right) ^{2}}>-\infty .$ Thus, $\psi \left( z\right) \sim z_{c}+\frac{1%
}{2}\psi ^{\prime \prime }\left( \tau \right) \left( z-\tau \right) ^{2}$
else $z\sim z_{c}+\frac{1}{2}\psi ^{\prime \prime }\left( \tau \right)
\left( \Phi \left( z\right) -\tau \right) ^{2}$ (a branch-point
singularity). It follows that $\Phi \left( \overline{z}\right) $ displays a
dominant power-singularity of order $-1/2$ at $z_{c}$ with $\Phi \left(
z_{c}\right) =\tau $ in the sense (recall $\sigma _{c}^{2}=\tau ^{2}\phi
^{\prime \prime }\left( \tau \right) /\phi \left( \tau \right) $) 
\begin{equation}
\bullet \text{ }\Phi \left( \overline{z}\right) \underset{\overline{z}%
\rightarrow z_{c}}{\sim }\tau -\sqrt{\frac{2\phi \left( \tau \right) }{\phi
^{\prime \prime }\left( \tau \right) }}\left( 1-\overline{z}/z_{c}\right)
^{1/2}=\tau \left( 1-\frac{\sqrt{2}}{\sigma _{c}}\left( 1-\overline{z}%
/z_{c}\right) ^{1/2}\right) .  \label{e10}
\end{equation}
By singularity analysis therefore [see Flajolet and Sedgewick (1993)], we
get [in agreement with Harris, 1963, Theorem 13.1, p. $32$] 
\begin{equation}
\bullet \text{ }\mathbf{P}\left( \overline{N}_{\infty }\left( 1\right)
=k\right) =\left[ \overline{z}^{k}\right] \Phi \left( \overline{z}\right) 
\underset{k\rightarrow \infty }{\sim }\sqrt{\frac{\phi \left( \tau \right) }{%
2Pi\phi ^{\prime \prime }\left( \tau \right) }}k^{-3/2}z_{c}^{-k}+O\left(
k^{-5/2}z_{c}^{-k}\right) ,  \label{e11}
\end{equation}
to the dominant order in $k$, with a geometric decay term at rate $%
z_{c}^{-1}=\phi ^{\prime }\left( \tau \right) <1$\ and a `universal'
power-law decay term $k^{-3/2}$\textit{.} When $\mu =\phi ^{\prime }\left(
1\right) \rightarrow 1$ (critical case) then both $\tau $ and $%
z_{c}\rightarrow 1$ and the above estimate boils down to a pure power-law
with $\left[ \overline{z}^{k}\right] \Phi \left( \overline{z}\right) 
\underset{k\rightarrow \infty }{\sim }\frac{1}{\sqrt{2\pi \phi ^{\prime
\prime }\left( 1\right) }}k^{-3/2}$. It can more precisely be checked that
when $\left| \mu -1\right| \ll 1$, $z_{c}^{-1}\sim 1-\left( \mu -1\right)
^{2}.$

Note finally that with $F\left( \lambda \right) =\log \phi \left(
e^{-\lambda }\right) $ the log-Laplace transform of $M$, $\tau >0$ is also
the solution to $F^{\prime }\left( -\log \tau \right) =1$.

Just like the computation of $\rho _{e}$ in the general case, the
computation of $\tau $, as a fixed point, can be quite involved (although
explicit in the LF case). A power-series expansion of $\tau $ and $z_{c}$ in
terms of the variable $\mu -1$ can formally be obtained, particularly useful
when the model is nearly critical. Define $\overline{\phi }\left( z\right) $%
\ by $\phi \left( z\right) =1+\mu \left( z-1\right) +\overline{\phi }\left(
1-z\right) $, so with $\overline{\phi }\left( 0\right) =0.$\ With $\overline{%
\tau }=1-\tau $, the equation $\phi \left( \tau \right) -\tau \phi ^{\prime
}\left( \tau \right) =0$\ giving $\tau $\ becomes 
\begin{equation*}
\delta \left( \overline{\tau }\right) :=\overline{\phi }\left( \overline{%
\tau }\right) +\left( 1-\overline{\tau }\right) \overline{\phi }^{\prime
}\left( \overline{\tau }\right) =\mu -1.
\end{equation*}
By Lagrange inversion formula [see Comtet, (1970)], we get:

1/ $\overline{\tau }=\overline{\tau }\left( \mu -1\right) =\sum_{n\geq
1}\tau _{n}\left( \mu -1\right) ^{n},$ where 
\begin{equation*}
\tau _{n}=\frac{1}{n}\left[ x^{n-1}\right] \left( \frac{\delta \left(
x\right) }{x}\right) ^{-n}.
\end{equation*}

2/ $z_{c}=1/\phi ^{\prime }\left( 1-\overline{\tau }\right) =z\left( 
\overline{\tau }\right) =z_{c}\left( \mu -1\right) =\sum_{n\geq
1}z_{n}\left( \mu -1\right) ^{n},$ where 
\begin{equation*}
z_{n}=\frac{1}{n}\left[ x^{n-1}\right] \left( z^{\prime }\left( x\right) 
\frac{\delta \left( x\right) }{x}\right) ^{-n}.
\end{equation*}
\newline

\textbf{Example:} Let $\phi \left( z\right) =\pi /\left( 1-\overline{\pi }%
z\right) $ [the Geo$_{0}\left( \pi \right) $ branching mechanism]. Let us
briefly work out this explicit Geo$_{0}\left( \pi \right) $ case, where $%
\phi \left( z\right) $ has convergence radius $z_{*}=1/\overline{\pi }$. The
r.v. $M$ has mean $\mu =\overline{\pi }/\pi $ and variance $\sigma ^{2}=%
\overline{\pi }/\pi ^{2}=\mu /\pi .$

If $\mu <1$ ($\overline{\pi }<1/2):$ $\rho _{e}=1<\tau =1/\left( 2\overline{%
\pi }\right) <z_{*}=1/\overline{\pi }.$ We have $\phi \left( \tau \right)
=2\pi $ and $z_{c}=1/\left( 4\overline{\pi }\pi \right) >1.$ Note $%
z_{c}<z_{*}$.

If $\mu =1$ ($\overline{\pi }=1/2):$ $\rho _{e}=1=\tau <z_{*}=2.$ We have $%
\phi \left( \tau \right) =1$ and $z_{c}=1.$

If $\mu >1$ ($\overline{\pi }>1/2):$ $\rho _{e}=\pi /\overline{\pi }<\tau
=1/\left( 2\overline{\pi }\right) <1<z_{*}=1/\overline{\pi }<2.$ We have $%
\phi \left( \tau \right) =2\overline{\pi }$ and $z_{c}=1/\left( 4\overline{%
\pi }\pi \right) >1.$ Note $z_{c}\lessgtr z_{*}$ if $\overline{\pi }\lessgtr
3/4$ and $\rho _{e}=1/\mu $ with $\overline{\rho }_{e}\sim 2\left( \mu
-1\right) /\sigma ^{2}$ as $\mu \rightarrow 1_{+}$ ($\overline{\pi }%
\rightarrow \left( 1/2\right) _{+}$)$.$ $\triangleright $

\subsection{\textbf{The pure power-law case (geometric tilting)}}

Define the tilted new p.g.f. $\Phi _{c}\left( \overline{z}\right) =\Phi
\left( \overline{z}z_{c}\right) /\Phi \left( z_{c}\right) $ and let $%
\overline{N}_{c}\left( 1\right) $ be the r.v. such that $\Phi _{c}\left( 
\overline{z}\right) =\mathbf{E}\left( \overline{z}^{\overline{N}_{c}\left(
1\right) }\right) .$ With $\widetilde{\phi }_{c}\left( z\right) =z_{c}\phi
\left( \tau z\right) /\tau =\phi \left( \tau z\right) /\phi \left( \tau
\right) $ defining the new rescaled branching p.g.f. encountered in Section $%
2$, we have 
\begin{equation}
\Phi _{c}\left( \overline{z}\right) =\overline{z}\widetilde{\phi }_{c}\left(
\Phi _{c}\left( \overline{z}\right) \right) .  \label{e12}
\end{equation}
Thus, $\overline{N}_{c}\left( 1\right) $ is the tree size of a BGW process
with one single founder when the generating branching mechanism is $%
\widetilde{\phi }_{c}\left( z\right) $. We note $\widetilde{\phi }_{c}\left(
1\right) =1$, $\widetilde{\phi }_{c}^{\prime }\left( 1\right) =z_{c}\phi
^{\prime }\left( \tau \right) =1$ (a critical case with extinction
probability $\rho _{e}=1,$ the smallest positive root of $\rho _{e}=%
\widetilde{\phi }_{c}\left( \rho _{e}\right) $) and the convergence radius
of $\widetilde{\phi }_{c}$ is $z_{*}/\tau >1$. As a result, $\Phi _{c}\left( 
\overline{z}\right) \underset{\overline{z}\rightarrow 1}{\sim }1-\tau ^{-1}%
\sqrt{\frac{2\phi \left( \tau \right) }{\phi ^{\prime \prime }\left( \tau
\right) }}\left( 1-\overline{z}\right) ^{1/2},$ with singularity displaced
to the left at $1,$ so that 
\begin{equation}
\bullet \text{ }\mathbf{P}\left( \overline{N}_{c}\left( 1\right) =k\right)
=\left[ \overline{z}^{k}\right] \Phi _{c}\left( \overline{z}\right) 
\underset{k\rightarrow \infty }{\sim }\tau ^{-1}\sqrt{\frac{\phi \left( \tau
\right) }{2Pi\phi ^{\prime \prime }\left( \tau \right) }}k^{-3/2}.
\label{e13}
\end{equation}
The geometric cutoff appearing in the probability mass of $\overline{N}%
_{\infty }\left( 1\right) $ has been removed and we are left with a pure
power-law case. This means that looking at the tree size p.g.f. $\Phi
_{c}\left( \overline{z}\right) $ generated by the critical branching
mechanism $\widetilde{\phi }_{c}\left( z\right) =z_{c}\phi \left( \tau
z\right) /\tau $, $\Phi _{c}\left( \overline{z}\right) $ exhibits a
power-singularity of order $-1/2$ at $z_{c}=1$ so that the new tree size
probability mass has pure power-law tails of order $1/2$. In particular $%
\mathbf{E}\left( \overline{N}_{c}\left( 1\right) \right) =\infty .$ In the
explicit Geo$_{0}\left( \beta \right) $ example above, where $\phi \left(
z\right) =\beta /\left( 1-\alpha z\right) $, it can be checked that $%
\widetilde{\phi }_{c}\left( z\right) =1/\left( 2-z\right) ;$ when dealing
with $\phi \left( z\right) =\left( \beta /\left( 1-\alpha z\right) \right)
^{\theta }$, $\widetilde{\phi }_{c}\left( z\right) =\left( \theta /\left(
\theta +1-z\right) \right) ^{\theta }$. Similarly, when dealing with the
binomial p.g.f. $\phi \left( z\right) =\left( 1-\alpha +\alpha z\right) ^{d}$%
, $\widetilde{\phi }_{c}\left( z\right) =\left( 1-1/d+z/d\right) ^{d}$ and
when dealing with the Poisson p.g.f. $\phi \left( z\right) =e^{-\mu \left(
1-z\right) }$, $\widetilde{\phi }_{c}\left( z\right) =e^{-\left( 1-z\right)
}.$ \newline

Consider the critical LF branching model for which $\phi \left( z\right) $
is given by (\ref{Mlaw}). Solving $\Phi \left( \overline{z}\right) =%
\overline{z}\phi \left( \Phi \left( \overline{z}\right) \right) $ for this $%
\phi $ yields: 
\begin{equation*}
\Phi \left( \overline{z}\right) =\frac{1}{2\overline{\pi }}\left( 1-%
\overline{z}\left( \pi -\overline{\pi }\right) -\sqrt{\left( 1-z\right)
\left( 1-\left( \pi -\overline{\pi }\right) ^{2}\overset{}{\overline{z}}%
\right) }\right) ,
\end{equation*}
with dominant singularity at $z_{c}=1$ ($\left( \pi -\overline{\pi }\right)
^{-2}>1$)$.$ Therefore 
\begin{eqnarray*}
&&\Phi \left( \overline{z}\right) \underset{z\rightarrow 1}{\sim }1-\sqrt{%
\frac{\pi }{\overline{\pi }}}\sqrt{1-z} \\
\mathbf{P}\left( \overline{N}\left( 1\right) =k\right) &=&\left[ \overline{z}%
^{k}\right] \Phi \left( \overline{z}\right) \underset{k\rightarrow \infty }{%
\sim }\sqrt{\frac{1}{2Pi\phi ^{\prime \prime }\left( 1\right) }}k^{-3/2}.
\end{eqnarray*}
Here both $\tau $ and $\phi \left( \tau \right) $ equal $1$ and $\phi
^{\prime \prime }\left( 1\right) =\sigma ^{2}=2\overline{\pi }/\pi $. 
\newline

Binomial (polynomial) and Poisson (exponential) models are examples of $\phi 
$ having convergence radius $z_{*}=\infty $. For the negative binomial
model, $\phi $ exhibits a power-singularity of positive order $\theta >0$ at 
$z_{*}=1/\alpha $ with $1<z_{*}<\infty $, so with $\phi \left( z_{*}\right)
=\infty .$ \newline

Here is now a family of $\phi $'s with a power-singularity of negative order 
$-\alpha $, $\alpha \in \left( 0,1\right) .$ Let $\alpha ,\lambda \in \left(
0,1\right) $ and $z_{*}>1$. Define the (damped) Sibuya p.g.f., [see Sibuya
(1979)], 
\begin{equation*}
h\left( z\right) =1-\lambda \left( 1-z/z_{*}\right) ^{\alpha }\text{ and }%
\phi \left( z\right) =\frac{h\left( z\right) }{h\left( 1\right) }.
\end{equation*}
It can be checked that this $\phi $ is a proper p.g.f. with convergence
radius $z_{*}$ and which is finite at $z=z_{*}>1$, with $\phi \left(
z_{*}\right) =\frac{1}{h\left( 1\right) }>1$. Note that for $k\geq 1,$

\begin{equation*}
\pi _{k}=\left[ z^{k}\right] \phi \left( z\right) =\frac{\lambda }{h\left(
1\right) }\left( -1\right) ^{k-1}\binom{\alpha }{k}z_{*}^{k}\underset{%
k\rightarrow \infty }{\sim }\frac{\lambda \alpha }{h\left( 1\right) }%
k^{-\left( \alpha +1\right) }z_{*}^{k}/\Gamma \left( 1-\alpha \right) .
\end{equation*}
The latter singularity expansion of $\Phi $ applies to this branching
mechanism $\phi $ as well.

\subsection{\textbf{Total number of leaves (sterile individuals) versus
total progeny}}

In the branching population models just discussed it is important to control
the number of leaves in the BGW tree with a single founder because leaves
are nodes (individuals) of the tree (population) that gave birth to no
offspring (the frontier of the tree as sterile individuals), so responsible
of its extinction. Leaves are nodes with outdegree zero, so let $\overline{N}%
^{0}\left( 1\right) $ be the number of leaves in a BGW tree with $\overline{N%
}\left( 1\right) $ nodes. With $\Phi \left( \overline{z}_{0},\overline{z}%
\right) =\mathbf{E}\left( \overline{z}_{0}^{\overline{N}^{0}\left( 1\right) }%
\overline{z}^{\overline{N}\left( 1\right) }\right) $ the joint p.g.f. of $%
\left( \overline{N}^{0}\left( 1\right) ,\overline{N}\left( 1\right) \right) $
solves the functional equation 
\begin{equation}
\bullet \text{ }\Phi \left( \overline{z}_{0},\overline{z}\right) =\overline{z%
}\left( \pi _{0}\left( \overline{z}_{0}-1\right) +\phi \left( \Phi \left( 
\overline{z}_{0},\overline{z}\right) \right) \right) .  \label{e14}
\end{equation}
With $\overline{N}^{0}\left( 1;k\right) :=\overline{N}^{0}\left( 1\right)
\mid \overline{N}\left( 1\right) =k$, we have 
\begin{equation*}
\mathbf{E}\left( \overline{z}_{0}^{\overline{N}^{0}\left( 1\right) }\mid 
\overline{N}\left( 1\right) =k\right) =\frac{\left[ \overline{z}^{k}\right]
\Phi \left( \overline{z}_{0},\overline{z}\right) }{\left[ \overline{z}%
^{k}\right] \Phi \left( 1,\overline{z}\right) },
\end{equation*}
where $\Phi \left( 1,\overline{z}\right) =\Phi \left( \overline{z}\right) $.
It is shown using this in [Drmota (2009), Th. $3.13$, page $84$] that, under
our assumptions on $\phi $, 
\begin{equation}
\begin{array}{l}
\frac{1}{k}\mathbf{E}\left( \overline{N}^{0}\left( 1;k\right) \right) 
\underset{k\rightarrow \infty }{\rightarrow }m_{0}=\frac{\pi _{0}}{\phi
\left( \tau \right) } \\ 
\frac{1}{k}\sigma ^{2}\left( \overline{N}^{0}\left( 1;k\right) \right) 
\underset{k\rightarrow \infty }{\rightarrow }\sigma _{0}^{2}=\frac{\pi _{0}}{%
\phi \left( \tau \right) }-\frac{\pi _{0}^{2}}{\phi \left( \tau \right) ^{2}}%
-\frac{\pi _{0}^{2}}{\tau ^{2}\phi \left( \tau \right) ^{2}\phi ^{\prime
\prime }\left( \tau \right) } \\ 
\frac{\overline{N}^{0}\left( 1;k\right) -m_{0}k}{\sigma _{0}\sqrt{k}}%
\underset{k\rightarrow \infty }{\overset{d}{\rightarrow }}\mathcal{N}\left(
0,1\right) .
\end{array}
\label{e15}
\end{equation}
As $k\rightarrow \infty $, $\frac{1}{k}\overline{N}^{0}\left( 1;k\right) $
converges in probability to $m_{0}<1$, the asymptotic fraction of nodes in a
size-$k$ tree which are leaves. For the Geo$_{0}\left( \pi _{0}\right) $
generated tree with $\phi \left( z\right) =\pi _{0}/\left( 1-\overline{\pi }%
_{0}z\right) $, it can be checked that $m_{0}=1/2$, whereas for the Poisson
generated tree with p.g.f. $\phi \left( z\right) =e^{\mu \left( z-1\right) }$%
, $m_{0}=e^{-1}.$ For the negative binomial tree generated by $\phi \left(
z\right) =\left( \beta /\left( 1-\alpha z\right) \right) ^{\theta }$, $%
m_{0}=\left( \theta /\left( \theta +1\right) \right) ^{\theta }$ and for the
Flory $d$-ary tree generated by the p.g.f. $\phi \left( z\right) =\left(
1-\alpha +\alpha z\right) ^{d}$, $m_{0}=\left( \left( d-1\right) /d\right)
^{d}.$\newline

\emph{Almost sure convergence and large deviations:} The functional equation
solving $\Phi \left( \overline{z}_{0},\overline{z}\right) $\ may be put
under the form

\begin{equation*}
\Phi \left( \overline{z}_{0},\overline{z}\right) =\overline{z}\phi _{%
\overline{z}_{0}}\left( \Phi \left( \overline{z}_{0},\overline{z}\right)
\right) ,
\end{equation*}
where $\phi _{\overline{z}_{0}}\left( z\right) =\pi _{0}\left( \overline{z}%
_{0}-1\right) +\phi \left( z\right) $, with $\overline{z}_{0}$\ viewed as a
parameter. Let $\tau \left( \overline{z}_{0}\right) $\ solve $\phi _{%
\overline{z}_{0}}\left( \tau \left( \overline{z}_{0}\right) \right) -\tau
\left( \overline{z}_{0}\right) \phi _{\overline{z}_{0}}^{\prime }\left( \tau
\left( \overline{z}_{0}\right) \right) =0$, else 
\begin{equation*}
\pi _{0}\left( \overline{z}_{0}-1\right) +\phi \left( \tau \left( \overline{z%
}_{0}\right) \right) -\tau \left( \overline{z}_{0}\right) \phi ^{\prime
}\left( \tau \left( \overline{z}_{0}\right) \right) =0
\end{equation*}
with $\tau \left( 1\right) =\tau $. We have $\left[ \overline{z}^{k}\right]
\Phi \left( 1,\overline{z}\right) ^{1/k}\rightarrow 1/z_{c}=\phi ^{\prime
}\left( \tau \right) $\ and $\left[ \overline{z}^{k}\right] \Phi \left( 
\overline{z}_{0},\overline{z}\right) ^{1/k}\rightarrow 1/z_{c}\left( 
\overline{z}_{0}\right) =\phi ^{\prime }\left( \tau \left( \overline{z}%
_{0}\right) \right) $, therefore 
\begin{equation*}
\mathbf{E}\left( \overline{z}_{0}^{\overline{N}^{0}\left( 1\right) }\mid 
\overline{N}\left( 1\right) =k\right) ^{1/k}=\left( \frac{\left[ \overline{z}%
^{k}\right] \Phi \left( \overline{\overline{z}}_{0},\overline{z}\right) }{%
\left[ \overline{z}^{k}\right] \Phi \left( 1,\overline{z}\right) }\right)
^{1/k}\rightarrow \alpha \left( \overline{z}_{0}\right) =\frac{\phi ^{\prime
}\left( \tau \left( 1\right) \right) }{\phi ^{\prime }\left( \tau \left( 
\overline{z}_{0}\right) \right) }.
\end{equation*}
By Cram\'{e}r's theorem, [Cram\'{e}r, (1938)], for $\rho \gtrless \rho
_{*}=F^{\prime }\left( 0\right) >0,$\ 
\begin{equation*}
\lim_{k\rightarrow \infty }\frac{1}{k}\log \mathbf{P}\left( \frac{1}{k}%
\overline{N}^{0}\left( 1;k\right) \gtrless \rho \right) =f\left( \rho
\right) ,
\end{equation*}
where, with $F\left( \lambda \right) =\log \alpha \left( e^{-\lambda
}\right) ,$\ 
\begin{equation*}
f\left( \rho \right) :=\inf_{\lambda \geq 0}\left( \lambda \rho -F\left(
\lambda \right) \right) \leq 0.
\end{equation*}
The function $f$\ is the large deviation rate function, as the Legendre
transform of the concave function $F$. The value $\rho _{*}$\ of $\rho $\
for which $f\left( \rho \right) =0$\ is $F^{\prime }\left( 0\right) $. We
conclude in particular that as $k\rightarrow \infty $%
\begin{equation}
\frac{1}{k}\left( \overline{N}^{0}\left( 1;k\right) \right) \overset{a.s.}{%
\rightarrow }\rho _{*}=F^{\prime }\left( 0\right) .  \label{e16}
\end{equation}

\textbf{Examples: }$\left( i\right) $ With $\overline{\pi }_{0}:=1-\pi _{0}$%
, let $\phi \left( z\right) =\pi _{0}/\left( 1-\overline{\pi }_{0}z\right) $
[the Geo$_{0}\left( \pi _{0}\right) $ branching mechanism]. Then, $\tau
\left( \overline{z}_{0}\right) =\frac{\overline{z}_{0}-\sqrt{\overline{z}_{0}%
}}{\overline{\pi }_{0}\left( \overline{z}_{0}-1\right) }$, so with 
\begin{equation*}
\phi ^{\prime }\left( \tau \left( \overline{z}_{0}\right) \right) =\pi _{0}%
\overline{\pi }_{0}\left( \sqrt{\overline{z}_{0}}+1\right) ^{2}\text{ and }%
\phi ^{\prime }\left( \tau \left( 1\right) \right) =4\pi _{0}\overline{\pi }%
_{0},
\end{equation*}
leading to $\alpha \left( \overline{z}_{0}\right) =4\left( \sqrt{\overline{z}%
_{0}}+1\right) ^{-2}$ and $F\left( \lambda \right) =\log 4-2\log \left(
1+e^{-\lambda /2}\right) $ with $F^{\prime }\left( 0\right) =1/2.$ So here $%
\frac{1}{k}\overline{N}_{k}^{0}\left( 1\right) \overset{a.s.}{\rightarrow }%
1/2$ (not only in probability). One can check 
\begin{equation*}
f\left( \rho \right) =-\log 4-2\left( \rho \log \rho +\left( 1-\rho \right)
\log \left( 1-\rho \right) \right) ,
\end{equation*}
is the Cram\'{e}r's large deviation rate function for this example.

Clearly, with $\overline{N}^{1}\left( 1;k\right) :=\overline{N}^{1}\left(
1;k\right) \mid \overline{N}\left( 1\right) =k$ the number of nodes in a
size-$k$ tree which are not leaves (prolific nodes), $\overline{N}^{0}\left(
1;k\right) +\overline{N}^{1}\left( 1;k\right) =k$ and $\frac{1}{k}\overline{N%
}^{1}\left( 1;k\right) \overset{a.s.}{\rightarrow }1-\rho _{*}.$

$\left( ii\right) $ With $\pi _{0}+\pi _{1}+\pi _{2}=1$, let $\phi \left(
z\right) =\pi _{0}+\pi _{1}z+\pi _{2}z^{2}$ [the binary branching
mechanism]. Then, $\tau \left( \overline{z}_{0}\right) =\sqrt{\left( \pi _{0}%
\overline{z}_{0}\right) /\pi _{2}}$, so with 
\begin{equation*}
\phi ^{\prime }\left( \tau \left( \overline{z}_{0}\right) \right) =\pi _{1}+2%
\sqrt{\pi _{0}\pi _{2}\overline{z}_{0}}\text{ and }\phi ^{\prime }\left(
\tau \left( 1\right) \right) =\pi _{1}+2\sqrt{\pi _{0}\pi _{2}},
\end{equation*}
leading to $\alpha \left( \overline{z}_{0}\right) =\left( \pi _{1}+2\sqrt{%
\pi _{0}\pi _{2}}\right) /\left( \pi _{1}+2\sqrt{\pi _{0}\pi _{2}\overline{z}%
_{0}}\right) $ and $F\left( \lambda \right) =\log \alpha \left( e^{-\lambda
}\right) $ with $F^{\prime }\left( 0\right) =\left( \pi _{0}\pi _{2}\right)
/\left( \pi _{1}\sqrt{\pi _{0}\pi _{2}}+2\pi _{0}\pi _{2}\right) <1.$ So
here $\frac{1}{k}\overline{N}^{0}\left( 1;k\right) \overset{a.s.}{%
\rightarrow }F^{\prime }\left( 0\right) $. If $\pi _{0}=\pi _{2}=1/4$ and $%
\pi _{1}=1/2$, $F^{\prime }\left( 0\right) =1/4.$ $\triangleright $

\subsection{\textbf{Forests of trees with a random number of trees}}

When there are more than one founders, we are left with a forest of
independent trees.

Let $\phi _{0}\left( z\right) $ (the p.g.f. of the random number of initial
founders) with $\mu _{0}=\phi _{0}^{\prime }\left( 1\right) <\infty $ be
such that $\Psi \left( \overline{z}\right) =\phi _{0}\left( \Phi \left( 
\overline{z}\right) \right) $, the p.g.f. of the total number of nodes of
the forest, has itself a dominant power-singularity at $z_{c}$ of order $%
-1/2 $, so with $\Psi \left( \overline{z}\right) \underset{\overline{z}%
\rightarrow z_{c}}{\sim }\phi _{0}\left( \tau \right) +\phi _{0}^{\prime
}\left( \tau \right) \sqrt{\frac{2\phi \left( \tau \right) }{\phi ^{\prime
\prime }\left( \tau \right) }}\left( 1-\overline{z}/z_{c}\right) ^{1/2}$.
Then, with $\overline{\mathcal{N}}$ the total number of nodes of the forest 
\begin{equation}
\mathbf{P}\left( \overline{\mathcal{N}}=k\right) =\left[ \overline{z}%
^{k}\right] \Psi \left( \overline{z}\right) \underset{k\rightarrow \infty }{%
\sim }\phi _{0}^{\prime }\left( \tau \right) \sqrt{\frac{\phi \left( \tau
\right) }{2Pi\phi ^{\prime \prime }\left( \tau \right) }}k^{-3/2}z_{c}^{-k}.
\label{e17}
\end{equation}
Note that $z_{c}=1$ when $\tau =1$ and $\phi ^{\prime }\left( 1\right) =1$
in which critical case $\left[ \overline{z}^{k}\right] \Psi \left( \overline{%
z}\right) \underset{k\rightarrow \infty }{\sim }\mu _{0}\sqrt{\frac{1}{2\pi
\phi ^{\prime \prime }\left( 1\right) }}k^{-3/2}$. For all branching
mechanism $\phi $ with convergence radius $z_{*}>1$ and all $\phi _{0}$ such
that $\Psi \left( \overline{z}\right) =\phi _{0}\left( \Phi \left( \overline{%
z}\right) \right) $ still has a dominant singularity at $z_{c}$, the law of
the size $\mathcal{N}$ of the forest of trees with a random number of
founders has a power-law factor which is $k^{-3/2}$, so independent of the
model's details (a universality property). Note that $z_{c}$ and the scaling
constant in front of $k^{-3/2}z_{c}^{-k}$ are model-dependent though, both
requiring the computation of $\tau $ which is known in the LF case.

\section{Relation to random walks (RW's)}

In this Section, we investigate two relations of BGW processes to RW's.

\subsection{\textbf{Relation to a skip-free to the left random walk}}

There is a natural connection between $\left\{ N_{n}\left( i\right) \right\} 
$ and $\left\{ \overline{N}_{n-1}\left( i\right) \right\} $ showing how the
full past determines the present of BGW trees.

Let $\psi _{n}\left( \lambda \right) :=-\log \phi _{n}\left( e^{-\lambda
}\right) $ and $\psi \left( \lambda \right) :=-\log \mathbf{E}\left(
e^{-\lambda \left( M-1\right) }\right) $. The recursion (\ref{1.1}) is also 
\begin{equation}
\psi _{n+1}\left( \lambda \right) -\psi _{n}\left( \lambda \right) =\psi
\left( \psi _{n}\left( \lambda \right) \right) =\psi _{n}\left( \psi \left(
\lambda \right) \right) \text{, }n\geq 0,\text{ }\psi _{0}\left( \lambda
\right) =\lambda ,  \label{1.1b}
\end{equation}
involving the increment of the log-Laplace transform of $N_{n}\left(
1\right) $. Let $M^{-}:=M-1,$ taking values in $\left\{ -1,0,1,2,...\right\} 
$. BGW processes are intimately related to a homogeneous random walk.
Consider indeed the (skip-free to the left) random walk $S_{n+1}\left(
i\right) =S_{n}\left( i\right) +M_{n+1}^{-}$, $S_{0}\left( i\right) =i\geq 1$%
, with the sequence $\left\{ M_{n}^{-}\right\} $ i.i.d. all distributed like 
$M^{-}.$ Consider $\mathcal{S}_{n}\left( i\right) :=S_{n\wedge \theta
_{i,0}}\left( i\right) $ where $\theta _{i,0}=\inf \left( n\geq
1:S_{n}\left( i\right) =0\right) $. Then, with $\overline{N}_{n}\left(
i\right) $ the cumulated number of offspring up to time $n$ of a
supercritical BGW process started with $i$ founders, $\mathcal{S}_{\overline{%
N}_{n+1}\left( i\right) }\left( i\right) -\mathcal{S}_{\overline{N}%
_{n}\left( i\right) }\left( i\right) =\sum_{i=\overline{N}_{n-1}\left(
i\right) +1}^{\overline{N}_{n}\left( i\right) }M_{i}^{-}\overset{d}{=}%
\sum_{k=1}^{N_{n}\left( i\right) }M_{k}^{-},$ showing, consistently with (%
\ref{1.1b}), that, 
\begin{equation}
\bullet \text{ }N_{n}\left( i\right) \overset{d}{=}\mathcal{S}_{\overline{N}%
_{n-1}\left( i\right) }\left( i\right) ,\text{ }n\geq 1.  \label{lamp}
\end{equation}
Therefore $N_{n}\left( i\right) $\ is a (Lamperti's) time-changed version of 
$\mathcal{S}_{n}\left( i\right) .$\ 

Given $\theta _{i,0}=\infty $ (an event with probability $1-\rho _{e}^{i}$,
see below), $N_{n}\left( i\right) \overset{d}{=}S_{\overline{N}_{n-1}\left(
i\right) }\left( i\right) $ and given $\theta _{i,0}<\infty $ (an event with
probability $\rho _{e}^{i}$), $N_{n}\left( i\right) \overset{d}{=}S_{%
\overline{N}_{n-1}\left( i\right) \wedge \theta _{i,0}}\left( i\right) 
\overset{d}{=}S_{\overline{N}_{n-1}\left( i\right) }\left( i\right) $,
because $\theta _{i,0}\overset{d}{=}\overline{N}_{\tau _{i,0}}\left(
i\right) >\overline{N}_{n-1}\left( i\right) $ (see below).

The relation (\ref{lamp}) is the discrete space-time version of a theorem by
[Lamperti, 1967] in the context of continuous-state branching processes.
Although probably in the folklore, we could not find a clear reference where
it is enounced.\newline

The one-step transition matrix $\Pi ^{\left( 0\right) }=\Pi ^{\left(
0\right) }\left( i,j\right) $ of the random walk $\left\{ S_{n}^{\left(
0\right) }\left( i\right) \right\} $ with state $\left\{ 0\right\} $
absorbing, is $\mathbf{P}\left( S_{1}^{\left( 0\right) }\left( 0\right)
=j\right) =P_{0,j}=\delta _{0,j}$ and

\begin{eqnarray*}
\mathbf{P}\left( S_{1}^{\left( 0\right) }\left( i\right) =j\right) &=&\Pi
^{\left( 0\right) }\left( i,j\right) =\left[ z^{j-i+1}\right] \phi \left(
z\right) =\pi _{j+i-1}\text{(}\underset{\text{LF}}{=}\overline{\pi }_{0}\pi 
\overline{\pi }^{j-i}\text{), }i\geq 1,\text{ }j\geq i \\
\mathbf{P}\left( S_{1}^{\left( 0\right) }\left( i\right) =i-1\right) &=&\Pi
^{\left( 0\right) }\left( i,i-1\right) =\pi _{0}\text{, }i\geq 1,\text{ }%
j=i-1
\end{eqnarray*}
It is an upper-Hessenberg type matrix, with state $\left\{ 0\right\} $
isolated. The harmonic function of $\left\{ S_{n}^{\left( 0\right) }\left(
i\right) \right\} $, say $\mathbf{h}^{\prime }\equiv \left( h\left( 1\right)
,h\left( 2\right) ,...\right) ^{\prime },$ is the smallest solution to 
\begin{equation*}
\sum_{j\geq i-1}\Pi ^{\left( 0\right) }\left( i,j\right) h\left( j\right)
=h\left( i\right) ,\text{ }i\geq 1,
\end{equation*}
with conventional boundary condition $h\left( 1\right) =1.$ It is then an
increasing sequence. With $\rho _{e},$ the smallest solution to $\phi \left(
z\right) =z$ , we get 
\begin{equation*}
h\left( i\right) =\frac{1-\rho _{e}^{i}}{1-\rho _{e}}\text{, }i\geq 1.
\end{equation*}
Indeed, for instance in the LF case, 
\begin{eqnarray*}
\sum_{j\geq i-1}\Pi ^{\left( 0\right) }\left( i,j\right) h\left( j\right)
&=&\sum_{k\geq 0}\Pi ^{\left( 0\right) }\left( i,k+i-1\right) h\left(
k+i-1\right) \\
&=&\frac{1}{1-\rho _{e}}\left( 1-\left( \pi _{0}\rho _{e}^{i-1}+\overline{%
\pi }_{0}\pi \sum_{k\geq 1}\overline{\pi }^{k}\rho _{e}^{k+i-1}\right)
\right) =\frac{1-\rho _{e}^{i}}{1-\rho _{e}},
\end{eqnarray*}
with $1/\overline{\rho }_{e}=h\left( \infty \right) $.

Because $\theta _{i,0}$ is also the first passage time to $0$ of the walk $%
S_{n}^{\left( 0\right) }\left( i\right) $ which can move downward at most $1$
at each step, $\theta _{i,0}$ is the sum of $i$ independent copies of $%
\theta _{1,0}$. Furthermore, by first-step analysis [Pitman, p. $124$], $%
\mathbf{E}\left( z^{\theta _{1,0}}\right) =\Phi \left( z\right) ,$ where $%
\Phi \left( z\right) $ solves the functional equation 
\begin{equation}
\Phi \left( z\right) =z\phi \left( \Phi \left( z\right) \right) ,
\label{feq}
\end{equation}
with $\Phi \left( 1\right) =\rho _{e}$. Consequently, $\mathbf{E}\left(
z^{\theta _{i,0}}\right) =\Phi \left( z\right) ^{i}$ (translating that $%
\theta _{i,0}\overset{d}{=}\overline{N}_{\tau _{i,0}}\left( i\right) $),
with $\mathbf{P}\left( \theta _{i,0}<\infty \right) =\rho _{e}^{i}.$ Note $%
\mathbf{E}\left( \theta _{i,0}\right) =\infty $ but $\mathbf{E}\left(
z^{\theta _{i,0}}\mid \theta _{i,0}<\infty \right) =\left( \Phi \left(
z\right) /\Phi \left( 1\right) \right) ^{i},$ leading to $\mathbf{E}\left(
\theta _{i,0}\mid \theta _{i,0}<\infty \right) =i\rho _{e}/\left( 1-\phi
^{\prime }\left( \rho _{e}\right) \right) $. We also have [see Norris], 
\begin{equation*}
h\left( i\right) =h\left( \infty \right) \left( 1-\mathbf{P}\left( \theta
_{i,0}<\infty \right) \right) ,
\end{equation*}
together with $\theta _{i,k}=\inf \left( n\geq 1:S_{n}^{\left( 0\right)
}\left( i\right) =k\right) ,$ $k>i$ and 
\begin{equation}
\mathbf{P}\left( \theta _{i,k}<\theta _{i,0}\right) =\frac{h\left( i\right) 
}{h\left( k\right) }=\frac{1-\rho _{e}^{i}}{1-\rho _{e}^{k}},  \label{hit}
\end{equation}
with $\mathbf{P}\left( \theta _{i,\infty }<\theta _{i,0}\right) =1-\rho
_{e}^{i}$, the probability of non-extinction given $S_{0}^{\left( 0\right)
}\left( i\right) =i$, as required.

We also observe that $\Phi \left( z\right) ^{i},$ as a solution to the
functional equation (\ref{feq}), is the p.g.f. of $\overline{N}_{\tau
_{i,0}}\left( i\right) $, the total limiting number of cumulated individuals
which appeared over time in the population (possibly infinite on the set of
explosion), when it has $i$ founders. It can be solved by the Lagrange
inversion formula, [see Comtet p. $159$]. Recalling $\mathbf{E}\left(
z^{\theta _{i,0}}\right) =\Phi \left( z\right) ^{i}$ and observing $\mathbf{E%
}\left( z^{S_{n}\left( i\right) }\right) =z^{i-n}\phi \left( z\right)
^{n}=z^{i}\mathbf{E}\left( z^{S_{n}\left( 0\right) }\right) $, by
Lagrange-B\"{u}rmann inversion formula, we get: 
\begin{equation*}
\left[ z^{n}\right] \Phi \left( z\right) ^{i}=\frac{i}{n}\left[
z^{n-i}\right] \phi \left( z\right) ^{n}=\frac{i}{n}\left[ z^{0}\right] 
\mathbf{E}\left( z^{S_{n}\left( i\right) }\right) .
\end{equation*}
Note we deal here with the `free' RW $S_{n}\left( i\right) $, not the one
absorbed at $0$.

The latter equality yields the Dwass-Kemperman formula [see Pitman, p. $124$%
] as 
\begin{equation}
\bullet \text{ }\mathbf{P}\left( \theta _{i,0}=n\right) =\frac{i}{n}\mathbf{P%
}\left( S_{n}\left( i\right) =0\right) ,\text{ }n\geq i,  \label{kemp}
\end{equation}
and more generally, while observing 
\begin{equation*}
\left[ z^{n+j}\right] \Phi \left( z\right) ^{i}=\frac{i}{n+j}\left[
z^{n+j-i}\right] \phi \left( z\right) ^{n}=\frac{i}{n+j}\left[ z^{j}\right] 
\mathbf{E}\left( z^{S_{n}\left( i\right) }\right) ,
\end{equation*}
\begin{equation}
\mathbf{P}\left( \theta _{i,0}=n+j\right) =\frac{i}{n+j}\mathbf{P}\left(
S_{n}\left( i\right) =j\right) ,\text{ }n\geq i-j.  \label{kempg}
\end{equation}
The prefactor $\frac{i}{n}$ in (\ref{kemp}) is thus the conditional
probability that $S_{n}\left( i\right) $ first hits $0$ at $n$, given $%
S_{n}\left( i\right) =0.$ Recalling the LF p.g.f. can be put under the form $%
\phi \left( z\right) =\frac{\pi _{0}-z\left( \overline{\pi }-\overline{\pi }%
_{0}\right) }{1-z\overline{\pi }},$%
\begin{equation*}
\mathbf{P}\left( S_{n}\left( i\right) =0\right) =\left[ z^{n-i}\right] \phi
\left( z\right) ^{n}=\pi _{0}^{n}\left[ z^{n-i}\right] \left( \frac{%
1+z\left( \overline{\pi }_{0}-\overline{\pi }\right) /\pi _{0}}{1-z\overline{%
\pi }}\right) ^{n}
\end{equation*}
where $\left( \frac{1+z\left( \overline{\pi }_{0}-\overline{\pi }\right)
/\pi _{0}}{1-z\overline{\pi }}\right) ^{n}$ is the product of $\left(
1+z\left( \overline{\pi }_{0}-\overline{\pi }\right) /\pi _{0}\right) ^{n}$
times $\left( 1-z\overline{\pi }\right) ^{-n}$ and so $\mathbf{P}\left(
S_{n}\left( i\right) =0\right) $ has the convolution explicit form 
\begin{equation}
\bullet \text{ }\mathbf{P}\left( S_{n}\left( i\right) =0\right) =\mathbf{P}%
\left( \overline{N}_{\tau _{i,0}}\left( i\right) =n\right) =\pi
_{0}^{n}\left( a*b\right) _{n-i}\text{ with}  \label{lawS0}
\end{equation}
\begin{equation*}
a_{k}=\binom{n}{k}\left[ \left( \overline{\pi }_{0}-\overline{\pi }\right)
/\pi _{0}\right] ^{k}\text{ and }b_{k}=\binom{n+k-1}{n-1}\overline{\pi }^{k}.
\end{equation*}
Multiplying this probability by $i/n$ also yields $\mathbf{P}\left( \theta
_{i,0}=n\right) $ explicitly for the LF model.

Finally, let $i,j\neq 1$\ and $j\geq \left( i-n\right) \wedge 0$. We have%
\emph{\ }

\begin{equation*}
S_{n}\left( i\right) =i+\sum_{k=1}^{n}M_{k}^{-}.
\end{equation*}
We therefore get, for $j\neq 0$ 
\begin{equation*}
\Pi ^{n}\left( i,j\right) :=\mathbf{P}\left( S_{n}\left( i\right) =j\right)
=\left[ z^{j-i+n}\right] \phi \left( z\right) ^{n}=\frac{n+j}{i}\mathbf{P}%
\left( \theta _{i,0}=n+j\right) ,
\end{equation*}
so that the resolvent of $\left\{ S_{n}\left( i\right) \right\} $ reads ($%
j\neq 0$) 
\begin{eqnarray*}
\bullet \text{ }g_{i,j}\left( u\right) &:&=\delta _{i,j}+\sum_{n\geq
1}u^{n}\Pi ^{n}\left( i,j\right) =\delta _{i,j}+\frac{1}{i}\sum_{n\geq
i-j}u^{n}\left( n+j\right) \mathbf{P}\left( \theta _{i,0}=n+j\right) \\
&=&\delta _{i,j}+\frac{u^{1-j}}{i}\sum_{k\geq i}u^{k-1}k\mathbf{P}\left(
\theta _{i,0}=k\right) =\delta _{i,j}+\frac{u^{1-j}}{i}\frac{d}{du}\left(
\Phi \left( u\right) ^{i}\right) .
\end{eqnarray*}
In particular, 
\begin{equation*}
g_{i,i}\left( u\right) =1+\frac{u^{1-i}}{i}\frac{d}{du}\left( \Phi \left(
u\right) ^{i}\right) .
\end{equation*}
This leads to the first return time $\theta _{i,i}$ of $\left\{ S_{n}\left(
i\right) \right\} $ to state $i$ p.g.f.: $\mathbf{E}\left( u^{\theta
_{i,i}}\right) =1-g_{i,i}\left( u\right) ^{-1}.$ In particular, $\mathbf{E}%
\left( u^{\theta _{1,1}}\right) =1-1/\left( 1+\Phi ^{\prime }\left( u\right)
\right) $ with, from (\ref{feq}), $\mathbf{P}\left( \theta _{1,1}<\infty
\right) =1-1/\left( 1+\Phi ^{\prime }\left( 1\right) \right) $ with $\Phi
^{\prime }\left( 1\right) =\rho _{e}/\left( 1-\phi ^{\prime }\left( \rho
_{e}\right) \right) $ and $\phi ^{\prime }\left( \rho _{e}\right) \in \left(
0,1\right) .$ The analysis makes use of the relation between $\theta _{i,0}$
and $\overline{N}_{\tau _{i,0}}\left( i\right) $. We refer to [Brown et al.,
(2010)] for related hitting time questions.\newline

We now come to the relation of BGW processes with the skip-free to the left
RW making use of (\ref{lamp})$.$

We first derive the scale function of the reflected skip-free to the right
RW $\left\{ R_{n}\left( i\right) =-S_{n}\left( -i\right) \right\} $, $i\geq
0 $, having moves up by one unit and arbitrary moves down [see Tak\'{a}cs
(1967), Marchal, (2001), Avram (2019)]$:$

With $k\geq i\geq 0$, define $\theta _{i,k}=\inf \left( n\geq 1:R_{n}\left(
i\right) \geq k\right) $ and $\theta _{i,-1}=\inf \left( n\geq 1:R_{n}\left(
i\right) \leq -1\right) .$ The scale function $w_{u}\left( k\right) $ of $%
\left\{ R_{n}\left( i\right) \right\} $ is defined by: 
\begin{equation*}
w_{u}\left( k\right) =\frac{1}{\mathbf{E}\left( u^{\theta _{.k}};\theta
_{.,k}<\theta _{.,-1}\right) },
\end{equation*}
where `$.$' is any initial state $i$. It has the scaling property ($k\geq
i\geq 1$) 
\begin{equation*}
\mathbf{E}\left( u^{\theta _{i,k}};\theta _{i,k}<\theta _{i,-1}\right) =%
\frac{w_{u}\left( i\right) }{w_{u}\left( k\right) }.
\end{equation*}
From the Markov property and the skip-free to the right property entailing $%
R_{\theta _{i,k}}\left( i\right) =k,$ by first-step analysis [see Marchal
(2001) Eq. 3.1], 
\begin{equation*}
w_{u}\left( i\right) =u\sum_{k=-1}^{i}\pi _{k+1}w_{u}\left( i-k\right) ,%
\text{ }i\geq 0.
\end{equation*}
This leads, if $\widehat{w}_{u}\left( z\right) =\sum_{i\geq 0}w_{u}\left(
i\right) z^{i}$, to 
\begin{equation*}
\widehat{w}_{u}\left( z\right) =\frac{1}{\phi \left( z\right) -z/u}\text{, }%
z\in \left( 0,\rho \left( u\right) \right) ,
\end{equation*}
where $\rho =\rho \left( u\right) $ solves $u\phi \left( \rho \right) =\rho $%
. The function $w\left( i\right) :=w_{1}\left( i\right) $ is called the
one-step scale function of $\left\{ R_{n}\left( i\right) \right\} $. It is
known from its generating function $\widehat{w}\left( z\right) :=\widehat{w}%
_{1}\left( z\right) =\frac{1}{\phi \left( z\right) -z}$, $z\in \left( 0,\rho
\left( 1\right) \right) .$

The scale function of $S_{n}\left( i\right) =-R_{n}\left( -i\right) $
coincides with the one of $R_{n}\left( i\right) ,$ but now, by symmetry,
with 
\begin{equation}
\mathbf{E}\left( u^{\theta _{i,-1}};\theta _{i,-1}<\theta _{i,k}\right) =%
\frac{w_{u}\left( k-i\right) }{w_{u}\left( k\right) },  \label{scale}
\end{equation}
where $\theta _{i,k}=\inf \left( n\geq 1:S_{n}\left( i\right) \geq k\right) $
and $\theta _{i,-1}=\inf \left( n\geq 1:S_{n}\left( i\right) \leq -1\right)
. $ When shifting up the initial condition of $S_{n}$ by one unit, we
recover the initial problem of when and how many times $S_{n}\left( i\right) 
$ passes through state $0$.\newline

Here, with $i\geq 1$, $\theta _{i,0}$ is the first hitting time of $0$ of
the non-increasing process $\min \left( S_{m}\left( i\right)
;m=0,...,n\right) $ with steps in $\left\{ 0,-1\right\} $. Let $\mathcal{S}%
_{n}^{*}\left( i\right) =\max \left( \mathcal{S}_{m}\left( i\right)
;m=0,...,n\right) <\infty $ on the event $\theta _{i,0}<\infty $, having
probability $\rho _{e}^{i}$.

The above solution of the two-sided exit problem in terms of the scale
function yields immediately the distribution of the overall maximum of $%
\mathcal{S}_{n}^{*}\left( i\right) ,$\ [see Bertoin, (2000)]. Putting $k=i+j$%
, for all $j\geq 0,$\ it indeed holds from (\ref{scale}) at $u=1$\ that : 
\begin{equation}
\bullet \text{ }\mathbf{P}\left( \mathcal{S}_{\theta _{i,0}}^{*}\left(
i\right) \leq i+j\right) =\frac{w\left( j\right) }{w\left( i+j\right) },
\label{Ber}
\end{equation}
since, if the event $\theta _{i,0}<\theta _{i,i+j}$\ is realized,
necessarily, $\max \left( \mathcal{S}_{m}\left( i\right) ;m=0,...,\theta
_{i,0}\right) \leq i+j.$

Based on the arguments developed in [Bennies and Kersting, (2000), p. $783$%
], given $\theta _{i,0}<\infty $, the time at which $S_{n}\left( i\right) $
attains its last maximum is uniformly distributed on $\left\{ i,...,\theta
_{i,0}\right\} .$

The process $N_{n}\left( i\right) $ being a time-changed version of $%
\mathcal{S}_{n}\left( i\right) $, the events ``$\mathcal{S}_{\theta
_{i,0}}^{*}\left( i\right) \leq i+j$'' and ``$N_{\tau _{i,0}}^{*}\left(
i\right) \leq i+j$'' coincide and we get 
\begin{equation}
\bullet \text{ }\mathbf{P}\left( N_{\tau _{i,0}}^{*}\leq i+j\right) =\frac{%
w\left( j\right) }{w\left( i+j\right) }\text{, }j\geq 0.  \label{width}
\end{equation}
The random variable $N_{\tau _{i,0}}^{*}$\ is the maximal value which the
process $\left\{ N_{n}\left( i\right) \right\} $\ can take in the course of
its history. It is known as the width of its profile (the apogee), the area
under the profile being $\left\{ \overline{N}_{\tau _{1,0}}\left( i\right)
\right\} $. Its distribution is given above, with possibly $N_{\tau
_{i,0}}^{*}\left( i\right) =\infty $\ in case of explosion ($\tau
_{i,0}=\infty $).\textit{\ }These results complete the ones of [Lindvall,
(1976)]. Note that $i^{\prime }>i\geq 1\Rightarrow N_{\tau _{i^{\prime
},0}}^{*}\left( i^{\prime }\right) \succ N_{\tau _{i,0}}^{*}\left( i\right) $%
, stochastically. Similarly, given $\tau _{i,0}<\infty $, the time at which $%
\left\{ N_{n}\left( i\right) \right\} $ attains its last maximum is
uniformly distributed on $\left\{ i,...,\tau _{i,0}\right\} .$\newline

\textbf{Examples:} $\left( i\right) $ For the LF p.g.f., one can check that 
\begin{eqnarray*}
\widehat{w}\left( z\right) &=&\frac{1}{\pi _{0}-\overline{\pi }}\frac{1-z%
\overline{\pi }}{\left( 1-z\right) \left( 1-z\overline{\pi }/\pi _{0}\right) 
}\text{ if }\mu \neq 1 \\
&=&\frac{1}{\overline{\pi }}\frac{1-z\overline{\pi }}{\left( 1-z\right) ^{2}}%
\text{ if }\mu =1.
\end{eqnarray*}
Therefore, with $w\left( 0\right) =1/\pi _{0}$, for all $j\geq 1$%
\begin{equation*}
w\left( j\right) =\left\{ 
\begin{array}{c}
\frac{\pi }{\pi _{0}-\overline{\pi }}\left[ 1-\left( \overline{\pi }_{0}/\pi
\right) \left( \pi _{0}/\overline{\pi }\right) ^{-\left( j+1\right) }\right]
>0\text{ if }\mu <1\text{ (subcritical)} \\ 
\frac{1}{\overline{\pi }}\left( 1+\pi j\right) >0\text{ if }\mu =1\text{
(critical)} \\ 
\frac{\pi }{\overline{\pi }-\pi _{0}}\left( \left( \overline{\pi }_{0}/\pi
\right) \left( \overline{\pi }/\pi _{0}\right) ^{j+1}-1\right) >0\text{ if }%
\mu >1\text{ (supercritical),}
\end{array}
\right.
\end{equation*}
increasing sequences in all cases, with $\rho =\pi _{0}/\overline{\pi }>1$ ($%
\mu =\overline{\pi }_{0}/\pi <1$) and $\rho =\rho _{e}=\pi _{0}/\overline{%
\pi }<1$ ($\mu >1$). The first sequence is bounded above, converging to $%
\frac{\pi }{\pi _{0}-\overline{\pi }}$.

From (\ref{width}), 
\begin{eqnarray*}
\mathbf{P}\left( N_{\tau _{i,0}}^{*}\left( i\right) <\infty \right) &=&1%
\text{ since }\frac{w\left( j\right) }{w\left( i+j\right) }\underset{%
j\rightarrow \infty }{\rightarrow }1\text{ if }\mu \leq 1\text{
((sub)-critical)} \\
\mathbf{P}\left( N_{\tau _{i,0}}^{*}\left( i\right) <\infty \right) &=&\rho
_{e}^{i}\text{ since }\frac{w\left( j\right) }{w\left( i+j\right) }\underset{%
j\rightarrow \infty }{\rightarrow }\rho _{e}^{i}\text{ if }\mu >1\text{
(supercritical).}
\end{eqnarray*}

- In the subcritical case ($\rho >1$), for large $j$, 
\begin{equation*}
\mathbf{P}\left( N_{\tau _{i,0}}^{*}\left( i\right) \leq i+j\right) =\frac{%
w\left( j\right) }{w\left( i+j\right) }=\frac{1-\left( \overline{\pi }%
_{0}/\pi \right) \left( \pi _{0}/\overline{\pi }\right) ^{-\left( j+1\right)
}}{1-\left( \overline{\pi }_{0}/\pi \right) \left( \pi _{0}/\overline{\pi }%
\right) ^{-\left( i+j+1\right) }}\sim 1-\mu \left( 1-\rho ^{-i}\right) \rho
^{-j},
\end{equation*}
with geometric decay.

- In the critical case, for large $j$, 
\begin{equation*}
\mathbf{P}\left( N_{\tau _{i,0}}^{*}\left( i\right) \leq i+j\right) =\frac{%
w\left( j\right) }{w\left( i+j\right) }=\frac{1+\pi j}{1+\pi \left(
i+j\right) }\sim 1-\frac{i}{i+j},
\end{equation*}
so that the width $N_{\tau _{i,0}}^{*}\left( i\right) -i$ of a BGW process
shows Pareto tails with index $1$ (just like its height $\tau _{i,0}$).

- In the supercritical case ($\rho _{e}<1$), for large $j$, 
\begin{equation*}
\mathbf{P}\left( N_{\tau _{i,0}}^{*}\left( i\right) \leq i+j\right) =\frac{%
w\left( j\right) }{w\left( i+j\right) }=\frac{\left( \overline{\pi }_{0}/\pi
\right) \left( \overline{\pi }/\pi _{0}\right) ^{j+1}-1}{\left( \overline{%
\pi }_{0}/\pi \right) \left( \overline{\pi }/\pi _{0}\right) ^{i+j+1}-1}\sim
\rho _{e}^{i}\left( 1-\frac{1}{\mu }\rho _{e}^{j}\right) ,
\end{equation*}
with geometric decay towards its limit.

Skip-free to the left RW's with LF upward jumps are important particular
cases in queueing and ruin theory, [Brown et al., (2010)]. In some gambling
game indeed, at each step there is a probability $\pi _{0}$ to lose one euro
and given a win phase (w.p. $\overline{\pi }_{0}$) the amount of the win is
Geo$_{0}\left( \pi \right) -$distributed, the shifted `time' till a first
failure: given an initial fortune $i\geq 1$, the first time of ruin is $%
\theta _{i,0}$ and $\mathcal{S}_{\theta _{i,0}}^{*}\left( i\right) $ the
largest amount of gains over this time window. For $\theta _{i,0}$ and $%
\mathcal{S}_{\theta _{i,0}}^{*}\left( i\right) $ to be finite w.p. $1$, it
is necessary that $\pi _{0}\geq \overline{\pi }$. In that case, ruin occurs
almost surely. For the gambler, the best situation is the critical case ($%
\pi _{0}=\overline{\pi }$) because then, the maximum of its gains is maximal
and also, the time till its eventual ruin is very long (Pareto tails). If
the maximal possible gain were known to him, a good strategy would be to
stop gambling when its gain has attained one of its maximum.\newline

$\left( ii\right) $ Another fundamental skip-free to the left RW (also
skip-free to the right) that will appear in the next section is when the
p.g.f. of $M$ is the one: $\phi _{b}\left( z\right) :=q+rz+pz^{2}$ ($q+r+p=1$%
) of a binary branching mechanism with holding probability $r$. In that
case, the p.g.f. of $M^{-}$ is the one of the simple RW (SRW): $z^{-1}\phi
_{b}\left( z\right) =qz^{-1}+r+pz.$ We have 
\begin{eqnarray*}
\frac{1}{\phi _{b}\left( z\right) -z} &=&\frac{1}{q\left( 1-z\right) \left(
1-pz/q\right) }=\frac{1}{q-p}\left( \frac{1}{1-z}-\frac{p/q}{1-pz/q}\right) 
\text{ if }p\neq q \\
&=&\frac{1}{p\left( 1-z\right) ^{2}}\text{ if }p=q,
\end{eqnarray*}
leading to the scale function 
\begin{equation}
w\left( j\right) =\left\{ 
\begin{array}{c}
\frac{1}{q-p}\left[ 1-\left( q/p\right) ^{-\left( j+1\right) }\right] >0%
\text{ if }p<q\text{ (subcritical)} \\ 
\frac{1}{p}\left( 1+j\right) \text{ if }p=q\text{ (critical)} \\ 
\frac{1}{p-q}\left( \left( p/q\right) ^{j+1}-1\right) >0\text{ if }p>q\text{
(supercritical).}
\end{array}
\right.  \label{scalebinary}
\end{equation}
It gives $\mathbf{P}\left( N_{\tau _{i,0}}^{*}\left( i\right) \leq
i+j\right) =\frac{w\left( j\right) }{w\left( i+j\right) },$ explicitly. This
distribution is also the one of the maximum of the SRW, started at $i$, till
it hits $0$ for the first time. $\triangleright $

\subsection{There is a linear-fractional BGW process nested inside the
profile of SRW's excursions}

The latter connection with a skip-free to the left RW is valid for any BGW
process and it takes an enlightening form when applied to the LF model. Here
is now a simple random walk (SRW) construction which is specific to the LF
model. In a special case, it is due to [Harris, (1952)]. An additional work
of us emphasizing the symmetries of the problem and fixing some technical
details is under preparation.

Let $\Sigma _{n}\left( 1\right) $ be the simple $\left\{ p,q,r\right\} -$%
nearest-neighbor random walk on the integers started at $\Sigma _{0}\left(
1\right) =1$ (so with increment or step $+1$ w.p. $p$, $-1$ w.p. $q$ and $0$
(stay alike) w.p. $r$, $p+q+r=1.$ We set $\Sigma _{-1}\left( 1\right) :=0$
and we assume that $\Sigma _{n}\left( 1\right) $ is stopped on its first
hitting time of $0$. With $\pm $ representing moves up and down, $0$ stay
alike moves and $s\cdot \cdot \cdot s$ concatenation of any length of $s-$%
steps moves, the profile of this SRW presents strings of the type: highlands 
$+0\cdot \cdot \cdot 0-$, valleys $-0\cdot \cdot \cdot 0+$ and terraces
either $+0\cdot \cdot \cdot 0+$ (left) or $-0\cdot \cdot \cdot 0-$ (right).
Before discussing the relationship of this SRW with a BGW tree, we start
with the following time-changed version of the SRW $\Sigma _{n}\left(
1\right) $: with $p_{0}=p/\left( 1-r\right) ,$ $q_{0}=q/\left( 1-r\right) $,
let $\Sigma _{n}\left( 1\right) $ be the simplest $\left\{
p_{0},q_{0},0\right\} -$nearest-neighbor random walk on the integers started
at $\Sigma _{0}\left( 1\right) =1$, absorbing the plateaus appearing in the
profile of $\Sigma _{n}\left( 1\right) .$ We also set \emph{$\Sigma $}$%
_{-1}\left( 1\right) =0$. The profile of this SRW presents maxima: $\wedge $%
, (else $+-$), minima: $\vee $ (else $-+$) and rise: $++$ or fall $--$
points. [Harris (1952)] considered this case. Define $\mathcal{\theta }%
_{1,0} $ to be the time of the first visit of $\Sigma _{n}\left( 1\right) $
to the origin (possibly $\infty $ with positive probability, not $1$). When
fixing at $0$ the height of the ground level, the upper part of this SRW
defines the profile of a non-trivial landscape (excursion) between times $-1$
and $\mathcal{\theta }_{1,0}$. For $h\geq 1,$ we define the additive
functional 
\begin{equation*}
\bullet \text{ }\mathcal{N}_{h}\left( 1\right) =\sum_{n=-1}^{\mathcal{\theta 
}_{1,0}-1}\mathbf{1}_{\left\{ \Sigma _{n}\left( 1\right) =h;\text{ }\Sigma
_{n+1}\left( 1\right) =h+1\right\} },\text{ }\mathcal{N}_{0}\left( 1\right)
=1,
\end{equation*}
counting the rise points at $h$ of the SRW. 
\begin{figure}[tph]
\begin{center}
\resizebox{22cm}{!}{\includegraphics{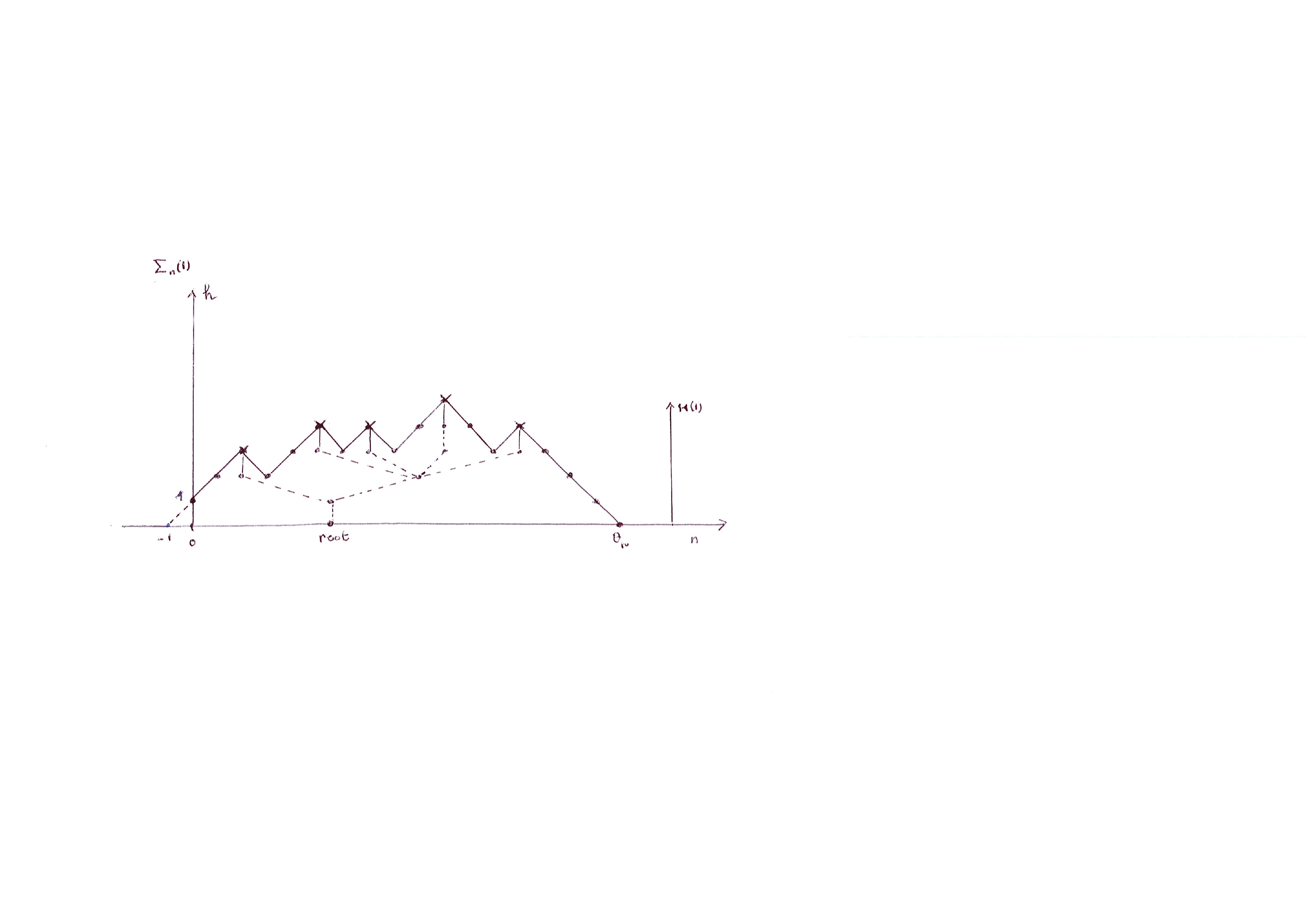}} %\label{SRWBGW}
\end{center}
\caption{A realization of the SRW $\Sigma _{n}\left( 1\right) $ and its
nested BGW (dotted lines). On top of its leaves, additional fictitious
branches were added to join the maxima of the SRW.}
\end{figure}

For instance, for the realization $0123234343454343210$ of a SRW excursion
starting at $0$ at $n=-1$, $\mathcal{N}_{h}\left( 1\right) =1,1,2,4,1$ for $%
h=0,...,4$, (see Figure $1$)$.$

Following [Harris, 1952], $\mathcal{N}_{h}\left( 1\right) $\ is a BGW
process with Geo$_{0}\left( q_{0}\right) $\ branching mechanism 
\begin{equation*}
\bullet \text{ }\mathcal{\phi }_{SRW}\left( z\right) =\text{ }q_{0}+\left(
1-q_{0}\right) \frac{\left( 1-p_{0}\right) z}{1-p_{0}z}=\frac{q_{0}}{1-p_{0}z%
}.
\end{equation*}
By first-step analysis indeed, for each $h\geq 1$ and given $m\geq 1$, the
probability $P\left( m\right) $ of any path event: $+\cdot \cdot \cdot
\left( m-1\right) +\cdot \cdot \cdot -$ coding the event that $m-1$ moves up
occurred at height $h$ between events of a rise point $+$ ($h-1\rightarrow h$%
) and the subsequent (random) fall point $-$ ($h\rightarrow h-1$) (without
the path visiting states below $h$ in between) obeys: 
\begin{equation*}
P\left( m\right) =\left\{ 
\begin{array}{c}
\delta _{m,0}\text{ w.p. }q_{0} \\ 
P\left( m-1\right) \text{ w.p. }p_{0}=1-q_{0}.
\end{array}
\right.
\end{equation*}
The p.g.f. of $P\left( m\right) $ is thus $\phi _{SRW}\left( z\right)
=\sum_{m\geq 0}P\left( m\right) z^{m}$ obeying $\phi _{SRW}\left( z\right)
=q_{0}+p_{0}z\phi _{SRW}\left( z\right) .$ Each such sub-excursion $+\cdot
\cdot \cdot \left( m-1\right) +\cdot \cdot \cdot -$ has the same law as the
one of the original full excursion when we shift its levels by a vertical
translation $-\left( h-1\right) $.

The range of $h$ is $\left\{ 0,...,\mathcal{H}\left( 1\right) -1\right\} $
where $\mathcal{H}\left( 1\right) $ is the height of any highest peak of the
SRW or the BGW's process lifetime shifted by one unit$.$ An individual in
the $(h-1)$-st generation has thus a probability $\left( 1-q_{0}\right)
p_{0}^{m-1}\left( 1-p_{0}\right) =q_{0}p_{0}^{m}$ of having exactly $m$
children, $m=0,1,..$. (The ancestor being the $0$ generation).

Clearly, there exists a similar quantity of $\mathcal{N}_{h}\left( 1\right) $
for the fall points $h\rightarrow h-1$ of the SRW for each $h\geq 1$. When
summing both over $h\geq 0$, these numbers are equal (when the SRW is a $%
\left( 0,0\right) -$ excursion) and the sum of the two is $1+\mathcal{\theta 
}_{1,0}$. Note that the local maxima of the walk are not taken into account
in $\mathcal{N}_{h}\left( 1\right) $.

The rise points of the SRW at $h$ are the offspring of the nodes of the BGW
tree at $h-1$. Starting from the root at $h=0$, the tree grows without
branching below the SRW's profile, until it meets its first $m\geq 1$ local
minima. The BGW process produces there $m+1$ offspring ($m$ of which
corresponding to finite sub-excursions of the SRW) and starting from the
right-most local first minimum, the process can be repeated along
independent pieces of the SRW's profile. The r.v. $\mathcal{N}_{h}\left(
1\right) -1$ is thus the number of local minima of an excursion of the SRW
at height $h$. For this BGW process, it is then easy to count the random
number of out-degree-$d$ nodes, $d=\left\{ 0,1,...\right\} $, for each $h$,
in particular the leaves ($d=0$) and the non-branching nodes ($d=1$). To
complete the picture, on top of each leaf of the tree a fictitious
additional branch (edge) can be added, ending up on the local maxima of the
SRW (see Figure $1$). The height of a highest leaf of the nested BGW tree
(its extinction time $\mathcal{\tau }_{1,0}$) corresponds to the one of a
highest peak $\mathcal{H}\left( 1\right) $ of the SRW, so with $\mathcal{H}%
\left( 1\right) =\mathcal{\tau }_{1,0}+1$. Clearly, $\overline{\mathcal{N}}_{%
\mathcal{\tau }_{1,0}}\left( 1\right) :=\sum_{h=0}^{\mathcal{\tau }_{1,0}}%
\mathcal{N}_{h}\left( 1\right) $ is the total number of individuals ever
born in the course of this Galton-Watson process and, with 
\begin{eqnarray*}
\mathcal{A}\left( 1\right) &=&\sum_{n=0}^{\mathcal{\theta }_{1,0}-1}\Sigma
_{n}\left( 1\right) \\
\mathcal{B}\left( 1\right) &=&\sum_{n=0}^{\mathcal{\theta }_{1,0}-1}\Sigma
_{n}\left( 1\right) -\sum_{n=0}^{\mathcal{\theta }_{1,0}-1}\left( \Sigma
_{n}\left( 1\right) -\Sigma _{n+1}\left( 1\right) \right) _{+},
\end{eqnarray*}
respectively the area and the restricted area under the (broken-line)
profile of the SRW, taking into account only the areas of the configurations 
$\square $ but not the ones $\bigtriangleup $ or pairs of half such
triangles: 
\begin{equation*}
\mathcal{B}\left( 1\right) =2\sum_{h=1}^{\mathcal{\tau }_{1,0}}h\mathcal{N}%
_{h}\left( 1\right) .
\end{equation*}
For the realization $0123234343454343210$ of the SRW excursion, $\mathcal{A}%
\left( 1\right) =51$, $\mathcal{B}\left( 1\right) =42$ and $\sum_{h=1}^{%
\mathcal{\tau }_{1,0}}h\mathcal{N}_{h}\left( 1\right) =21.$

Defining the p.g.f. $\emph{h}_{i}\left( z\right) :=\mathbf{E}\left( z^{%
\mathcal{\theta }_{i,0}}\right) $, upon conditioning on the first step 
\begin{equation*}
\emph{h}_{1}\left( z\right) =q_{0}z+p_{0}z\emph{h}_{2}\left( z\right)
\end{equation*}
with $\emph{h}_{2}\left( z\right) =\emph{h}_{1}\left( z\right) ^{2}.$ Thus, 
\begin{equation*}
\emph{h}\left( z\right) =\frac{1-\sqrt{1-4p_{0}q_{0}z^{2}}}{2p_{0}z},
\end{equation*}
and, with $k\geq 1,$%
\begin{equation*}
\mathbf{P}\left( \mathcal{\theta }_{1,0}=2k-1\right) =\frac{\left(
p_{0}q_{0}\right) ^{k}}{p_{0}}\frac{\binom{2k-2}{k-1}}{k!}.
\end{equation*}
In the critical case $p_{0}=q_{0}=1/2,$ $\mathcal{\theta }_{1,0}$ has Pareto
tails with index $1/2$.

We have $\emph{h}\left( 1\right) =\mathbf{P}\left( \mathcal{\theta }%
_{1,0}<\infty \right) =\frac{1-\left| p_{0}-q_{0}\right| }{2p_{0}}=1\wedge
q_{0}/p_{0}$, translating that when $p_{0}>q_{0}$, the SRW is transient. At
the same time, the nested BGW process is subcritical if $q_{0}/p_{0}>1$,
critical if $q_{0}/p_{0}=1$ and supercritical if $q_{0}/p_{0}<1.$ In this
last case, $\mathcal{\tau }_{1,0}=\inf \left( h\geq 1:\mathcal{N}_{h}\left(
1\right) =0\right) <\infty $ only with positive probability $\rho
_{e}=q_{0}/p_{0}$, the extinction probability of $\emph{N}_{h}\left(
1\right) $, whereas in the first two cases, $\mathcal{\tau }_{1,0}<\infty $
a.s.. As required, we have: $\mathbf{P}\left( \mathcal{\theta }_{1,0}<\infty
\right) =\mathbf{P}\left( \mathcal{\tau }_{1,0}<\infty \right) .$

The p.g.f. $\Phi \left( z\right) $ of the total number of nodes of a BGW
process with the geometric branching mechanism $\mathcal{\phi }_{SRW}\left(
z\right) $ solves $\Phi \left( z\right) =z\mathcal{\phi }_{SRW}\left( \Phi
\left( z\right) \right) $, leading to 
\begin{eqnarray*}
\Phi \left( z\right) &:&=\mathbf{E}\left( z^{\overline{\mathcal{N}}_{%
\mathcal{\tau }_{1,0}}\left( 1\right) }\right) =\frac{1-\sqrt{1-4p_{0}q_{0}z}%
}{2p_{0}} \\
\mathbf{P}\left( \overline{\mathcal{N}}_{\mathcal{\tau }_{1,0}}\left(
1\right) =k\right) &=&\frac{\left( p_{0}q_{0}\right) ^{k}}{p_{0}}\frac{%
\left( 2k-2\right) !}{k!\left( k-1\right) :},\text{ }k\geq 1.
\end{eqnarray*}
Hence: $\Phi \left( z^{2}\right) =z\emph{h}\left( z\right) ,$ translating
that $\mathcal{\theta }_{1,0}\overset{d}{=}2\overline{\mathcal{N}}_{\mathcal{%
\tau }_{1,0}}\left( 1\right) -1.$

In the critical case when $p_{0}=q_{0}=1/2$, it can be checked that $%
\overline{\mathcal{N}}_{\mathcal{\tau }_{1,0}}\left( 1\right) $ has Pareto
tails with index $1/2$ just like $\mathcal{\theta }_{1,0}$ therefore.\newline

\textbf{Remark:} the paths of the SRW excursion can be reconstructed from
the one of its nested BGW tree as follows (see Figure 1):

Start from any prolific node of the nested tree at height $h\geq 0$ and
consider the subtree rooted at this node. Browsing this subtree following
the contour process strategy yields the sub-excursion of the SRW above level 
$h+1.$ If the starting node has outdegree $d\geq 1$, there are $d$ returns
at level $h+1$ of the sub-excursion. In the processes, all the edges are
visited twice. [see Champagnat (2015), pages $33$ and $38$, for example]. $%
\blacksquare $\newline

Define now $\theta _{1,0}$ to be the time of the first visit to the origin
of the $\left\{ p,q,r\right\} -$SRW $\Sigma _{n}\left( 1\right) $. We define 
\begin{equation*}
A\left( 1\right) =\sum_{n=0}^{\theta _{1,0}-1}\Sigma _{n}\left( 1\right) ;%
\text{ }H\left( 1\right) =\underset{n=0,...,\theta _{1,0}-1}{\max }\Sigma
_{n}\left( 1\right) ,
\end{equation*}
respectively the area under the profile of $\Sigma _{n}\left( 1\right) $ and
its height. We also let, as before, 
\begin{equation*}
B\left( 1\right) =\sum_{n=0}^{\theta _{1,0}-1}\Sigma _{n}\left( 1\right)
-\sum_{n=0}^{\theta _{1,0}-1}\left( \Sigma _{n}\left( 1\right) -\Sigma
_{n+1}\left( 1\right) \right) _{+}
\end{equation*}
be the restricted area under the profile of $\Sigma _{n}\left( 1\right) $.

For each $h\geq 1,$ with $G_{n,h}$ an i.i.d. sequence of Geo$_{0}\left(
r\right) -$distributed r.v.'s define the additive functional 
\begin{equation*}
\bullet \text{ }N_{h}\left( 1\right) =\sum_{n=0}^{\theta _{1,0}-1}\left(
G_{n,h}\mathbf{1}_{\left\{ \Sigma _{n}\left( 1\right) =h;\text{ }\Sigma
_{n+1}\left( 1\right) =h+1\right\} }\right) ,\text{ }N_{0}\left( 1\right) =1.
\end{equation*}
In words, $N_{h}\left( 1\right) $ is the number of times that, before first
visiting $0,$ the random walk $\Sigma _{n}\left( 1\right) $ crosses from $h$
to $h+1$ (the rise points of the SRW) to which the sizes of each $h+1$%
-plateau forming the left terraces were attached. The local highlands of the
SRW are not taken into account in $N_{h}\left( 1\right) $. For instance, for
the realization $01(222)(33)23(4444)3434(55)4343210$ of a SRW excursion with
plateaus starting at $0$ at\emph{\ }$n=-1$\emph{, }$N_{h}\left( 1\right)
=1,3,3,7,2$\emph{\ }for\emph{\ }$h=0,...,4.$\emph{\ }Clearly there exists a
similar symmetric quantity for the fall points $h\rightarrow h-1$ of the SRW
for each $h\geq 1$, with plateaus now preceding the fall points (the right
terraces). When summing both over $h\geq 0$, these numbers are equal in
distribution (when the SRW is a $\left( 0,0\right) -$ excursion) and the sum
of the two is $1+\theta _{1,0}$ in law. If we define the width $\Sigma
_{\theta _{1,0}}^{**}\left( 1\right) $ of $\Sigma _{n}\left( 1\right) $ as
the largest size of the valleys in its profile, then 
\begin{equation*}
\Sigma _{\theta _{1,0}}^{**}\left( 1\right) =\max_{h=1,...,\tau
_{1,0}}N_{h}\left( 1\right) =:N_{\tau _{1,0}}^{*}\left( 1\right) ,
\end{equation*}
so the maximal value which the crossing process $N_{h}\left( 1\right) $ can
take.\newline

The sequence $N_{h}\left( 1\right) $\ is a BGW process with LF offspring
distribution (\ref{Mlaw}), with the correspondence: 
\begin{equation}
\overline{\pi }=p,\text{ }\pi _{0}=q\text{ and }\pi -\pi _{0}=r>0,
\label{corres}
\end{equation}
so with the restriction $\pi >\pi _{0}$\ ($r>0$). Equivalently, the
branching mechanism of the BGW process nested inside the SRW is 
\begin{mathletters}
\begin{equation}
\bullet \text{ }\phi _{SRW}\left( z\right) =q+\left( 1-q\right) \frac{\left(
1-p\right) z}{1-pz}=\frac{q+rz}{1-pz}.  \label{phisrw}
\end{equation}

\textit{Proof:} An individual in the $\left( h-1\right) $st generation has a
probability $\left( 1-q\right) p^{m-1}\left( 1-p\right) $ of having exactly $%
m$ children, $m=1,2,...$, given it is productive; its probability of having
no offspring being $q$ (The ancestor being at generation $0$). This
probability mass function enjoys the memory-less property of the geometric
distribution.

By first-step analysis, for each $h\geq 1$ and given $m\geq 1$, the
probability $P\left( m\right) $ of any path event $+\cdot \cdot \cdot \left(
m-1\right) +\cdot \cdot \cdot -$ coding the event that $m-1$ moves up
occurred at height $h$ between events of the type $+$ ($h-1\rightarrow h$)
and the subsequent $-$ ($h\rightarrow h-1$) (without the path visiting
states below $h$ in between and including no move steps between the
extremities) obeys: 
\end{mathletters}
\begin{equation*}
P\left( m\right) =\left\{ 
\begin{array}{c}
\delta _{m,0}\text{ w.p. }q \\ 
\delta _{m,1}\text{ w.p. }r \\ 
P\left( m-1\right) \text{ w.p. }p.
\end{array}
\right.
\end{equation*}
With probability $r$, the path configuration is $+0\cdot \cdot \cdot 0-$,
the one of a highland at height $h$, so with $m=1$ (no move up in between
the extremities). The p.g.f. of $P\left( m\right) $ is $\phi _{SRW}\left(
z\right) =\sum_{m\geq 0}P\left( m\right) z^{m}$ obeying $\phi _{SRW}\left(
z\right) =q+rz+pz\phi _{SRW}\left( z\right) .$ We call the linear-fractional
BGW process with branching mechanism $\phi _{SRW}\left( z\right) $ the
nested BGW process inside the $\left\{ p,q,r\right\} -$SRW.\newline

Note that the nested BGW process is subcritical if $\pi _{0}/\overline{\pi }%
=q/p>1$, critical if $\pi _{0}/\overline{\pi }=q/p=1$ and supercritical if $%
\pi _{0}/\overline{\pi }=q/p<1.$ In this last case, $\tau _{1,0}=\inf \left(
h\geq 1:N_{h}\left( 1\right) =0\right) <\infty $ only with positive
probability $\rho _{e}=q/p$, the extinction probability of $N_{h}\left(
1\right) $, whereas in the first two cases, $\theta _{1,0}<\infty $ a.s..
Concomitantly, $\mathbf{P}\left( \theta _{1,0}<\infty \right) =\frac{%
1-\left| p-q\right| }{2p}=1\wedge q/p$, translating that when $p>q$, the SRW
is transient at infinity. As required, in all cases, we have: $\mathbf{P}%
\left( \theta _{1,0}<\infty \right) =\mathbf{P}\left( \tau _{1,0}<\infty
\right) .$

As before, we let $\overline{N}_{\tau _{1,0}}\left( 1\right)
=\sum_{h=0}^{\tau _{1,0}}N_{h}\left( 1\right) $ the total number of nodes of
the BGW tree with branching mechanism (\ref{phisrw})$.$\newline

$\bullet $ By first-step analysis, we have the correspondences:

\begin{eqnarray*}
&&H\left( 1\right) \overset{d}{=}\tau _{1,0}+1 \\
&&\theta _{1,0}\overset{d}{=}2\overline{N}_{\tau _{1,0}}\left( 1\right) -1 \\
&&\Sigma _{\theta _{1,0}}^{**}\left( 1\right) \overset{d}{=}N_{\tau
_{1,0}}^{*}\left( 1\right) \\
&&B\left( 1\right) \overset{d}{=}2\sum_{h=1}^{H\left( 1\right)
-1}hN_{h}\left( 1\right)
\end{eqnarray*}

\textit{Proof:} We already mentioned the first one. We prove the second one,
the other ones being obtained similarly.

Defining the p.g.f. $h_{i}\left( z\right) :=\mathbf{E}\left( z^{\theta
_{i,0}}\right) $, upon conditioning on the first step 
\begin{equation*}
h_{1}\left( z\right) =qz+rzh_{1}\left( z\right) +pzh_{2}\left( z\right)
\end{equation*}
with $h_{2}\left( z\right) =h_{1}\left( z\right) ^{2}.$ Thus 
\begin{equation*}
h\left( z\right) :=h_{1}\left( z\right) =\mathbf{E}\left( z^{\theta
_{1,0}}\right) =\frac{1-rz-\sqrt{\left( 1-rz\right) ^{2}-4pqz^{2}}}{2pz}
\end{equation*}
with, as required while considering that the lengths of the plateaus are Geo$%
_{0}\left( r\right) -$distributed$,$%
\begin{equation*}
h\left( z\right) =\frac{1-\sqrt{1-4p_{0}q_{0}z^{2}}}{2p_{0}z}\mid
_{z\rightarrow \frac{\left( 1-r\right) z}{1-rz}}
\end{equation*}

Now the p.g.f. $\Phi \left( z\right) =\mathbf{E}\left( z^{\overline{N}_{\tau
_{1,0}}\left( 1\right) }\right) $ of the total number of nodes of a BGW
process with general LF branching mechanism $\phi \left( z\right) $ solves $%
\Phi \left( z\right) =z\phi \left( \Phi \left( z\right) \right) $, leading
to 
\begin{equation*}
\Phi \left( z\right) =\frac{1-\left( \pi -\pi _{0}\right) z-\sqrt{\left(
1-\left( \pi -\pi _{0}\right) z\right) ^{2}-4\pi _{0}\overline{\pi }z}}{2%
\overline{\pi }}.
\end{equation*}
When dealing with $\phi _{SRW}\left( z\right) $ with the correspondences (%
\ref{corres}), $\Phi \left( z\right) $ becomes $\Phi _{SRW}\left( z\right) $
obeying:

\begin{equation*}
\Phi _{SRW}\left( z^{2}\right) =zh\left( z\right) ,
\end{equation*}
translating that $\theta _{1,0}\overset{d}{=}2\overline{N}_{\tau
_{1,0}}\left( 1\right) -1.$ When $p=q$, $\theta _{1,0}$ has Pareto tails
with index $1/2$. This extends Theorem $5$ of [Harris, 1952].

- The height $H\left( 1\right) \geq 1$ of the highest maximum or peak of $%
\Sigma _{n}\left( 1\right) $ (as attained by the random walk before its
first return to the origin) can be identified to $\tau _{1,0}+1$ where $\tau
_{1,0}$ is the extinction time of the LF BGW tree associated to the SRW. Its
distribution is thus given recursively by 
\begin{equation*}
\mathbf{P}\left( H\left( 1\right) \leq h+1\right) =\phi _{SRW}\left( \mathbf{%
P}\left( H\left( 1\right) \leq h\right) \right) ,\text{ }h\geq 0.
\end{equation*}
We get 
\begin{eqnarray*}
\mathbf{P}\left( H\left( 1\right) >h\right) &=&\frac{1}{\left( \frac{1-p}{1-q%
}\right) ^{h}+\frac{p}{1-q}\left( 1+...+\left( \frac{1-p}{1-q}\right)
^{h-1}\right) } \\
&=&\frac{p-q}{p-q\left( \frac{1-p}{1-q}\right) ^{h}}\text{ if }p\neq q \\
&=&\frac{1-p}{1+p\left( h-1\right) }\text{ if }p=q.
\end{eqnarray*}
When $p=q$, $H\left( 1\right) $ has Pareto tails with index $1$.

- The height $L\left( 1\right) $ of the lowest minimum of $\Sigma _{n}\left(
1\right) $ (else, the height of the deepest valley of the SRW's profile) is
the first time at which a fictitious leaf appears in $\left\{ N_{h}\left(
1\right) \right\} $. Its distribution is Geo$\left( \pi _{1}\right) ,$ with $%
\pi _{1}=\phi _{SRW}^{^{\prime }}\left( 0\right) =\left( 1-p\right) \left(
1-q\right) ,$the probability that the associated BGW process generates a
single offspring, so 
\begin{equation*}
\mathbf{P}\left( L\left( 1\right) =h\right) =\overline{\pi }_{1}\pi
_{1}^{h-1},\text{ }h\geq 1.
\end{equation*}

- From (\ref{lawS0}), 
\begin{equation}
\mathbf{P}\left( \overline{N}_{\tau _{1,0}}\left( 1\right) =n\right)
=q^{n}\left( a*b\right) _{n-1}\text{ with}  \label{lawS1}
\end{equation}
\begin{equation*}
a_{k}=\binom{n}{k}\left[ r/q\right] ^{k}\text{ and }b_{k}=\binom{n+k-1}{n-1}%
p^{k}.
\end{equation*}
In the critical case when $p=q$, $\overline{N}_{\mathcal{\tau }_{1,0}}\left(
1\right) $ has Pareto tails with index $1/2$ and so does $\theta _{1,0}$
therefore.

- The law of $N_{\tau _{1,0}}^{*}\left( 1\right) $ and also of $\Sigma
_{\theta _{1,0}}^{**}\left( 1\right) $ is given by $\mathbf{P}\left( N_{\tau
_{1,0}}^{*}\left( 1\right) \leq 1+j\right) =\frac{w\left( j\right) }{w\left(
1+j\right) }$, with $w\left( j\right) $ the scale function defined in (\ref
{scalebinary}).\newline

The SRW $\Sigma _{n}\left( 1\right) $ looks very much like the profile of a
mountain chain. By considering the reflected SRW $-\Sigma _{n}\left(
1\right) $, the physical image is the one of a seabed profile. The highlands
of $\Sigma _{n}\left( 1\right) $ become the valleys of $-\Sigma _{n}\left(
1\right) .$

While considering the concatenation of $i\geq 2$ i.i.d. excursion landscapes
of the SRW (assuming state $0$ to be purely reflecting), the BGW to consider
is $N_{h}\left( i\right) ,$ being $i$ i.i.d. copies of $N_{h}\left( 1\right)
.$\newline

[Harris (1952)] observes that the connection of BGW processes and simple
SRWs remains valid if the transition probabilities $\left(
p_{h},q_{h}\right) ,$ $p_{h}+q_{h}=1,$ of the (non-homogeneous) SRW depend
on its current height $h$; the nested BGW process has then a corresponding
branching mechanism depending on the height: $\mathbf{P}\left(
M_{h}=m\right) =\left( 1-q_{h+1}\right) p_{h+1}^{m-1}\left( 1-p_{h+1}\right)
=q_{h+1}p_{h+1}^{m}$ ($p_{h+1}+q_{h+1}=1$), $h\geq 0,$ $m\geq 1$. The
iteration of variable height-dependent LF branching mechanisms is then
necessary.\newline

\begin{center}
\textbf{Acknowledgments:}
\end{center}

T. Huillet acknowledges partial support from the ``Chaire \textit{%
Mod\'{e}lisation math\'{e}matique et biodiversit\'{e}'' }of Veolia-Ecole
Polytechnique-MNHN-FondationX and support from the labex MME-DII Center of
Excellence (\textit{Mod\`{e}les math\'{e}matiques et \'{e}conomiques de la
dynamique, de l'incertitude et des interactions}, ANR-11-LABX-0023-01
project). This work was also funded by CY Initiative of Excellence (grant ``%
\textit{Investissements d'Avenir}''ANR- 16-IDEX-0008), Project ``EcoDep''
PSI-AAP2020-0000000013. S. Mart\'{i}nez was supported by the Center for
Mathematical Modeling ANID Basal FB210005.\newline

\begin{center}
\textbf{Declarations of interest}
\end{center}

The authors have no conflicts of interest associated with this paper.\newline

\begin{center}
\textbf{Data availability statement}
\end{center}

There are no data associated with this paper.

\end{document}